\documentclass[11pt]{amsart}

\usepackage[latin1]{inputenc}

\usepackage{subfigure, amsmath, verbatim, amsfonts, amssymb}

\usepackage{tikz}
\usepackage{xcolor}
\usetikzlibrary{arrows, shapes}
\usetikzlibrary{arrows.meta}
\allowdisplaybreaks

\newtheorem{theorem}{Theorem}
\newtheorem{proposition}[theorem]{Proposition}
\newtheorem{lemma}[theorem]{Lemma}
\newtheorem{corollary}[theorem]{Corollary}

\newtheorem*{theorem*}{Theorem}

\theoremstyle{definition}

\newtheorem{definition}[theorem]{Definition}

\newtheorem{remark}[theorem]{Remark}
\newtheorem{example}[theorem]{Example}
\newtheorem{examples}[theorem]{Examples}

\newcommand{\spe}[1]{\mathtt{#1}}

\DeclareMathOperator{\qsym}{QSym}

\DeclareMathOperator{\id}{id}
\DeclareMathOperator{\cano}{cano}

\DeclareMathOperator{\Hilb}{Hilb}
\DeclareMathOperator{\dHilb}{\mathbf{Hilb}}

\DeclareMathOperator{\ps}{ps}

\newcommand{\bdot}{\bullet}

\DeclareMathOperator{\dcone}{dcone}
\begin{document}

\title{Coloring Complexes and Combinatorial Hopf Monoids}

\author{Jacob A.~White}
\address{School of Mathematical and Statistical Sciences\\
University of Texas - Rio Grande Valley\\
Edinburg, TX 78539}
\keywords{Chromatic Polynomials, Quasisymmetric Functions, Combinatorial Species, Combinatorial Hopf Algebras, Balanced Simplicial Complexes} 

\subjclass{18D10, 18D35}
\date{\today}


\begin{abstract}

We generalize the notion of a coloring complex of a graph to linearized combinatorial Hopf monoids. We determine when a linearized combinatorial Hopf monoid has such a construction, and discover some inequalities that are satisfied by the quasisymmetric function invariants associated to the combinatorial Hopf monoid. We show that the collection of all such coloring complexes forms a linearized combinatorial Hopf monoid, which is the terminal object in the category of combinatorial Hopf monoids with convex characters. We also study several examples of combinatorial Hopf monoids.

\end{abstract}
\maketitle
\section{Introduction}

In their landmark paper, Aguiar, Bergeron and Sottile \cite{aguiar-bergeron-sottile} proved that $\qsym$, the Hopf algebra of quasisymmetric functions, forms the terminal combinatorial Hopf algebra. They showed how many well-known examples of quasisymmetric generating functions came from combinatorial Hopf algebras that had very simple characters. A fundamental example is the combinatorial Hopf algebra of graphs, where the resulting quasisymmetric function is Stanley's chromatic symmetric function \cite{stanley-coloring-1}. Other examples include quasisymmetric generating functions for $P$-partitions \cite{gessel}, and the Billera-Jia-Reiner quasisymmetric function associated to matroids \cite{billera-jia-reiner}. The goal of this paper is to show a similar type of result for linearized combinatorial Hopf monoids in species, and deduce some consequences of that result. Our terminal object is a species whose structures we refer to as \emph{abstract relative coloring complexes}. First, we discuss some background regarding coloring complexes of graphs, and applications of our new results to graph theory.

The chromatic polynomial can be studied using geometry: given a graph $G$, then up to a shift, the chromatic polynomial of $G$, denoted $\chi(G, x),$ is the Hilbert polynomial of an ideal \cite{steingrimsson} which arises from a relative simplicial complex $\Phi(G)$ which is the set difference of the Coxeter complex of type $A$ and a subcomplex $\Gamma(G)$, usually referred to as the \emph{coloring complex}. The coloring complex $\Gamma(G)$ has been used to find several new inequalities regarding the coefficients of the chromatic polynomial \cite{hersh-swartz,hultman}. For this paper, we refer to the relative simplicial complex $\Phi(G)$ as the coloring complex, and show how $\Phi(G)$ can be used to learn information about the chromatic symmetric function  of $G$, denoted $X(G, \mathbf{x})$. If $V$ is the vertex set of $G$, then the vertex set of $\Phi(G)$ is $2^V \setminus \{\emptyset, V\}$, the collection of all non-empty proper subsets of $V.$ A collection $\{V_1, \ldots, V_k\}$ is a face of $\Phi(G)$ if the following conditions are satisfied:
\begin{enumerate}
    \item We have $V_i \subsetneq V_{i+1}$ for all $i \in [k-1]$.
    \item We have $V_{i+1} \setminus V_i$ is an independent set of $G$ for all $i \in [k+1]$, where $V_{k+1} = V$ and $V_0 = \emptyset$.
\end{enumerate}

We briefly review the definitions of quasisymmetric functions, integer compositions, and monomial and fundamental quasisymmetric functions; see Section 7.19 of \cite{stanleyec2} for relevant definitions. 
Let $\mathbf{x} = x_1, x_2, \ldots $ be a sequence of commuting indeterminates. Recall that an \emph{integer composition} $\alpha$ of a positive integer $n$ is a sequence $(\alpha_1, \ldots, \alpha_k)$ of positive integers such that $\alpha_1+\cdots + \alpha_k = n$. We write $\ell(\alpha) = k$, and $\alpha \models n$. Let $n \in \mathbb{N}$ and let $F(\mathbf{x}) \in \mathbb{C}[[\mathbf{x}]]$ be a homogeneous formal power series in $\mathbf{x}$, where the degree of every monomial in $F(\mathbf{x})$ is $n$. Then $F(\mathbf{x})$ is a \emph{quasisymmetric function} if it satisfies the following property:
for every $\alpha \models n$ and every $i_1 < i_2 < \cdots < i_{\ell(\alpha)}$, we have $[\prod_{j=1}^{\ell(\alpha)} x_{i_j}^{\alpha_j}]F(\mathbf{x}) = [\prod_{j=1}^{\ell(\alpha)} x_j^{\alpha_j}]F(\mathbf{x})$. 

Given an integer composition $\alpha = (\alpha_1, \alpha_2, \ldots, \alpha_k)$ of $n$, we let \[M_{\alpha} = \sum\limits_{i_1 < \cdots < i_k} \prod_{j=1}^k x_{i_j}^{\alpha_j}.\] These are the \emph{monomial quasisymmetric functions}, which form a basis for the ring of quasisymmetric functions. We will often express elements of $\qsym$ as generating functions in $\mathbf{x}$, or in terms of the monomial quasisymmetric functions. Given a quasisymmetric function $f$ and an integer composition $\alpha$, let $[M_{\alpha}]f$ be the coefficient of $M_{\alpha}$ in the monomial quasisymmetric expansion of $f$.

The second basis we focus on is the basis of fundamental quasisymmetric functions, first introduced by Gessel \cite{gessel}.
The set of integer compositions is partially ordered by refinement: we write $\alpha \leq \beta$ if $\beta$ is a refinement of $\alpha$. Note that this order is dual to the typical order studied in the literature. With respect to this partial order, the set of integer compositions forms a lattice.
The \emph{fundamental quasisymmetric functions} $F_{\alpha}$ are defined by:
\[ F_{\alpha} = \sum\limits_{\beta \geq \alpha} M_{\beta}. \]

Let $G$ be a graph with vertex set $V$.
The complex $\Phi(G)$ is a \emph{balanced} relative simplicial complex of dimension $|V|-2$, which means that it is pure and comes with a coloring $\kappa: V(\Phi(G)) \to [|V|-1]$ such that every facet has exactly one vertex of each color. Given any balanced relative simplicial complex $\Phi$ on vertex set $N$ of dimension $d$ and $S \subseteq [d+1]$, we define \[f_S(\Phi) = |\{\sigma \in \Phi(G): \{\kappa(v): v \in \sigma \} = S \}|.\] This is the \emph{flag $f$-vector} of $\Phi$. If we write $S = \{s_1, s_2, \ldots, s_k \}$, with $s_1 < \cdots < s_k$, then we define $\alpha(S) = (s_1, s_2-s_1, \ldots, s_k-s_{k-1}, |V|-s_k)$. In Corollary \ref{cor:graphs}, we show that
\begin{align} X(G, \mathbf{x}) & = \sum_{\sigma \in \Phi} M_{\alpha(\kappa(\sigma))} \\ &= \sum_{S \subset [|V| - 1]} f_S(\Phi(G)) M_{\alpha(S)}.  \label{eq:flag} \end{align}
 Thus $f_S(\Phi(G)) = [M_{\alpha(S)}] X(G, \mathbf{x})$ and we can use the relative coloring complex to understand the chromatic symmetric function, similar to how the coloring complex was used to study the chromatic polynomial. 

In order to get new information using this perspective, we prove new results about balanced relative simplicial complexes.
In Section \ref{sec:geometry}, we study flag $f$-vectors of balanced relative simplicial complexes, as well as the $f$-vector of pure relative simplicial complexes. We also study colored relative simplicial complexes, which are relative simplicial complexes $\Phi$ equipped with a coloring function $\kappa: V(\Phi) \to [k]$ with the property that no two vertices of a face $\sigma \in \Phi$ receive the same color. In Subsection \ref{subsec:relsimple}, we review the definition of pure relative simplicial complexes, and $f$-vectors. 

\begin{definition}
Given a sequence $(a_0, \ldots, a_k)$, we say the sequence is \emph{strongly flawless} if the following two systems of inequalities are satisfied:
\begin{enumerate}
    \item We have $a_0 \leq a_1 \leq \cdots \leq a_{\lfloor k/2 \rfloor}$.
    \item We have $a_i \leq a_{k-i}$ for $0 \leq i \leq \lfloor k/2 \rfloor.$
\end{enumerate}
\end{definition}
Examples of strongly flawless sequences include pure O sequences \cite{hibi} and the $h$-vectors of broken circuit complexes of matroids \cite{kubitzke}. Recently, the Whitney numbers of the second kind of a realizable matroid have been shown to be strongly flawless \cite{huh-et-al}, solving the top-heavy conjecture \cite{top-heavy}.

We show that the $f$-vector of a pure relative simplicial complex is also strongly flawless. In fact, we show more is true. We say that a sequence $(a_0, \ldots, a_n)$ is \emph{super flawless} if, for all $i$, we have $(n-i)a_i \leq (i+1)a_{i+1}.$ We show in Proposition \ref{prop:flawlesspre} that every super flawless sequence is a strongly flawless sequence.

We present the following result.
\begin{proposition}
Let $(f_{-1}, f_0, \ldots, f_d)$ be the $f$-vector of a pure relative simplicial complex $\Phi$. Then $(f_{-1}, \ldots, f_d)$ is super flawless.
\label{prop:flawless}
\end{proposition}
The proof is immediately given after the proof of Proposition \ref{prop:flawlesspre}.
In Subsection \ref{subsec:balanced}, we study balanced relative simplicial complexes, and their flag $f$-vectors. The flag $f$-vector is a refinement of the $f$-vector, as we have \[f_{i-1}(\Phi) = \sum_{S \subseteq [d]: |S| = i} f_S(\Phi).\] We present the following result.
\begin{proposition}
Let $\Phi$ be a balanced relative simplicial complex of dimension $d$-1. Let $S \subseteq T \subseteq [d]$. Then $f_S(\Phi) \leq f_T(\Phi)$.
\label{prop:increasing}
\end{proposition}
The proof is given immediately after Examples \ref{ex:balancedcomplexes3}.
Since the chromatic symmetric function of a graph is related to the flag $f$-vector of $\Phi(G)$, we obtain inequalities about chromatic symmetric functions and chromatic polynomials. Moreover, since $\chi(G, x)$ can be obtained from $X(G, \mathbf{x})$ by principal specialization, the $f$-vector of $\Phi(G)$ can be used to obtain inequalities about $\chi(G, x)$ as well. Using Equation \eqref{eq:flag}, Proposition \ref{prop:flawless} and Proposition \ref{prop:increasing}, one can prove the following:
\begin{proposition}
Let $G$ be a graph on $n$ vertices.
For $\alpha \leq \beta \models n$ in the refinement order, we have 
\begin{equation} \label{eq:Mincreasing} [M_{\alpha}] X(G, \mathbf{x}) \leq [M_{\beta}] X(G, \mathbf{x}). \end{equation} 

Moreover, if we write $\chi(G, c) = \sum_{i=0}^n f_i \binom{c}{i}$, then $f_0 = 0$ and $(f_1, \ldots, f_n)$ is super flawless.
\label{cor:graphsintro}
\end{proposition}
We refer to a quasisymmetric function which satisfies the inequality in Equation \eqref{eq:Mincreasing} as being \emph{$M$-increasing}. In this paper, we use Proposition \ref{prop:flawless} and \ref{prop:increasing} in the proof of Theorem \ref{thm:main} below, from which Corollary \ref{cor:graphs} and Proposition \ref{cor:graphsintro} follows. 

Our primary goal in this paper is to generalize the construction of the coloring complex of a graph to other combinatorial objects and quasisymmetric functions. Given a collection $\mathcal{C}$ of combinatorial objects, and a quasisymmetric function $\Psi(O, \mathbf{x})$ for each $O \in \mathcal{C}$, we wish to find balanced relative simplicial complexes $\Phi(O)$ such that an analogue of Equation \eqref{eq:flag} holds, with $X(G,\mathbf{x})$ replaced by $\Psi(O,\mathbf{x})$. Whenever we can show that such a construction exists, then we know that the quasisymmetric functions are $M$-increasing, and that the corresponding polynomials are strongly flawless. 

We observe that the definition of the coloring complex of $G$ involves flags $V_1 \subsetneq V_2 \subsetneq \cdots \subsetneq V_k$ such that $G$ restricted to $V_i \setminus V_{i-1}$ is $1$-colorable. In general, we see that, given a combinatorial object $O$ with a vertex set $V$, and a notion of $1$-colorable, we want flags $V_1 \subsetneq V_2 \subsetneq \cdots \subsetneq V_k \subsetneq V$ such that restricting $O$ to $S_i \setminus S_{i-1}$ results in something $1$-colorable. One well-known setting for discussing combinatorial objects with vertex sets (or labels), is the setting of combinatorial species. It is also well-known that coalgebras are a useful way to study the notion of restriction. Hence, if we work with combinatorial Hopf monoids in species, we will have a general notion of restriction and label set.

Our focus will be on combinatorial objects that form a linearized combinatorial Hopf monoid, and where the quasisymmetric function is induced by a linearized character.
In Section \ref{sec:combhopf}, we define linearized combinatorial Hopf monoids in species. Roughly speaking, we start with a linear species $\spe{H}$, also known in the literature as a vector species or species. For each finite set $N$, we have a vector space $\spe{H}_N$. We require $\spe{H}$ to be a Hopf monoid object in the category of linear species equipped with a character $\varphi$, where each vector space $\spe{H}_N$ comes with a distinguished basis $\mathcal{H}_N$ (which we refer to as $\spe{H}$-structures), such that:
\begin{enumerate}
    \item The product of basis elements is a basis element.
    \item When the coproduct of a basis element is nonzero, then it is a simple tensor of basis elements.
    \item The character takes on values $0$ or $1$ on basis elements.
\end{enumerate}
Our notion of linearization is more general than in \cite{aguiar-mahajan-2,linearization}, because we are linearizing \emph{pointed set species}, and regarding the base points as being the zero vector. For example, the Hopf monoid of posets discussed in this paper is linearized in our sense, but not linearized from set species, because the coproduct of a poset might be zero. Our decision to study linearized characters is motivated by the fact that many examples of characters studied in the literature are linearized from pointed set species.
In Subsection \ref{subsec:linearized}, we discuss the linearized Hopf monoids of graph $\mathbb{K}\spe{G}^{\bdot}$, posets $\mathbb{K}\spe{P}^{\bdot}$, and matroids $\mathbb{K}\spe{M}^{\bdot}$. All three examples have appeared in \cite{aguiar-mahajan-1,aguiar-mahajan-2,aguiar-ardila}.

For a Hopf monoid $\spe{H}$ in species, given a finite set $N$, the coproduct map $\Delta$ decomposes as a direct sum $\bigoplus_{S \sqcup T} \Delta_{S,T}$ where $\Delta_{S,T}:\spe{H}_N \mapsto \spe{H}_S \otimes \spe{H}_T$ where $\sqcup$ denotes disjoint union. Informally, we will think of $\Delta_{S,T}(\spe{h})$ as decomposing $\spe{h}$ into disjoint label sets. The axioms of a Hopf monoid ensures that, given $N = C_1 \sqcup C_2 \sqcup \cdots \sqcup C_k$, the expression $\Delta_{C_1, C_2, \ldots, C_k}(\spe{h})$ is well-defined. We discuss this in more detail in Section \ref{sec:combhopf}.

We also discuss the $\varphi$-chromatic quasisymmetric function, $\Psi_{\varphi}(\spe{h}, \mathbf{x})$, which is an invariant associated to any $\spe{H}$-structure $\spe{h}$ in a linearized combinatorial Hopf monoid $\spe{H}$. The definition is given in Equation \eqref{eq:deffunction}. The existence of this invariant comes from the work of Aguiar and Mahajan \cite{aguiar-mahajan-1} and Aguiar, Bergeron, and Sottile \cite{aguiar-bergeron-sottile}, although in Theorem \ref{thm:coloring} we give an alternative combinatorial description of the invariant in terms of $\varphi$-proper colorings. Our primary focus will be on the \emph{chromatic character} $\chi$ which we introduce in this paper in Definition \ref{def:chromatic}. It is defined for any linearized Hopf monoid, and generalizes well-known characters for graphs, posets, matroids, and generalized permutohedra. We describe the character in detail for the Hopf monoid of graphs, posets, and matroids, and derive the resulting chromatic quasisymmetric function. For example, from our general theory, $X(G, \mathbf{x}) = \Psi_{\chi}(\spe{g}, \mathbf{x})$, where $\chi$ is the chromatic character on the Hopf monoid of graphs.

Given a linearized Hopf monoid $\spe{H}$ with character $\varphi$, we want to associate a balanced relative simplicial complex $\Phi_{\varphi}(\spe{h})$ to each $\spe{H}$-structure $\spe{h} \in \spe{H}_N$ such that \[\Psi_{\varphi}(\spe{h}, \mathbf{x}) = \sum_{S \subseteq [|N| - 1]} f_S(\Phi_{\varphi}(\spe{h})) M_{\alpha(S)}.\]

Given a pointed set species $\spe{H}$ such that the linearization $\mathbb{K}\spe{H}$ is a linearized Hopf monoid with linearized character $\varphi$, we wish to find a set species $\spe{K}$, and a natural transformation $\Phi: \spe{H} \to \spe{K}$, such that:
\begin{enumerate}
    \item Every non-zero $\spe{K}$-structure $\spe{k}_N$ is a colored relative simplicial complex on $N$, and 
    \item For every non-zero $\spe{H}$-structure $\spe{h}$, we have 
    \[\Psi_{\varphi}(\spe{h}, \mathbf{x}) = \sum_{S \subseteq [|V| - 1]} f_S(\Phi_{\varphi}(\spe{h})) M_{\alpha(S)}.\]
\end{enumerate}
We refer to $\Phi$ as the \emph{geometric realization} of $\spe{H}$.
In Section \ref{sec:realization}, we give sufficient conditions on the linearized character $\varphi$ of a linearized combinatorial Hopf monoid that ensures that $\spe{H}$ has a geometric realization. The conditions require notation from Section \ref{sec:combhopf}; see Definition \ref{def:convex}. Characters which satisfy the conditions listed in Definition \ref{def:convex} are called \emph{convex characters}. We define the $\varphi$-coloring complex $\Phi_{\varphi}(\spe{h})$ in Definition \ref{def:coloringcomplex}. We define \[\Phi_{\varphi}(\spe{h}) = \{ F_1 \subsetneq F_2 \subsetneq \cdots \subsetneq F_k: \varphi^{\otimes k+1} \circ \Delta_{F_1, F_2 \setminus F_1, \ldots, N \setminus F_k}(\spe{h}) = 1 \}\] where $\varphi^{\otimes k+1} = \varphi \otimes \varphi \otimes \cdots \otimes \varphi$ where there are $k+1$ terms.
Informally, we are looking at flags where `restricting' $\spe{h}$ to $F_i \setminus F_{i-1}$ is $1$-colorable for all $i$, where `restriction' comes from the coproduct, and saying $\spe{k}$ is `$1$-colorable' means that $\varphi(\spe{k}) = 1$.
We show that the coloring complexes do form a geometric realization for $\spe{H}$ whenever $\varphi$ is a convex character. 

We also give sufficient conditions on when a linearized combinatorial Hopf monoid $\mathbb{K}\spe{H}$ has a geometric realization $\Phi$ such that, for every non-trivial $\spe{H}$-structure $\spe{h} \in \spe{H}_N$, the complex $\Phi(\spe{h})$ is balanced of dimension $|N|-2$. We refer to the corresponding characters as \emph{balanced convex characters}. Let $\Delta$ be the coproduct.
We say $\varphi$ is a balanced convex character if it satisfies the following conditions for all $N$ and $\spe{h} \in \mathcal{H}_N$:
\begin{enumerate}
    \item If $|N| = 1$, then $\varphi(\spe{h}) = 1$.
    \item If $|N| > 1$, then $\Delta(\spe{h}) \neq 1 \otimes \spe{h} + \spe{h} \otimes 1.$
    \item If $\varphi(\spe{h}) = 1$, and we write $\Delta(\spe{h}) = \sum_{i=1}^k \spe{x}_i \otimes \spe{y}_i$ where $\Delta$ is the coproduct, then $\varphi(\spe{x}_i) = \varphi(\spe{y}_i) = 1$ for all $i.$
\end{enumerate} 
We present the following result, which is a combination of Theorem \ref{thm:geometricrealization} and Theorem \ref{thm:geometricbalanced}.
\begin{theorem}
Let $(\mathbb{K}\spe{H},  \mathbb{K}(\varphi))$ be a linearized combinatorial Hopf monoid with a convex character $\varphi$. Given an $\spe{H}$-structure $\spe{h}$, let $\Phi_{\varphi}(\spe{h})$ be the $\varphi$-coloring complex of $\spe{h}$. Then \[\Psi_{\varphi}(\spe{h}, \mathbf{x}) = \sum_{S \subseteq [|N| - 1]} f_S(\Phi_{\varphi}(\spe{h})) M_{\alpha(S)}, \] where $N$ is the vertex set of $\spe{h}$. 

If $\varphi$ is a balanced convex character, then $\Phi_{\varphi}(\spe{h})$ is balanced for all nontrivial $\spe{h}$. In that case, $\Psi_{\varphi}(\spe{h}, \mathbf{x})$ is $M$-increasing. If we write $\chi_{\varphi}(\spe{h},x) = \sum_{i=0}^{|N|} f_i\binom{x}{i}$, then $f_0 = 0$ and $(f_1, \ldots, f_{|N|})$ is super flawless.
\label{thm:main}
\end{theorem}
In Theorem \ref{thm:chiconvex}, we show that the chromatic character $\chi$ is always convex, and give conditions under which it is a balanced convex character. We show that these conditions are met for the Hopf monoids of graphs, posets, and matroids, and describe the corresponding geometric realizations. In the case of graphs, we show that $\Phi_{\chi}(\spe{g})$ is the relative coloring complex described in this introduction. We also show how Theorem \ref{thm:main} implies inequalities for the chromatic symmetric function of a graph, the $P$-partition enumerator of a poset, and the Billera-Jia-Reiner quasisymmetric function of a matroid. These inequalities can also be derived from the fact that the corresponding quasisymmetric functions are all $F$-positive. However, our approach generalizes to other invariants for mixed graphs and double posets, where $F$-positivity no longer holds.

Thus, we have a vast generalization of the coloring complex of a graph. A natural problem is to classify the possible coloring complexes that can arise from our construction. That is, given a relative simplicial complex $\Phi$, does it arise as $\Phi_{\varphi}(\spe{h})$ for \emph{some} combinatorial Hopf monoid $(\mathbb{K}\spe{H}, \mathbb{K}(\varphi))$ and $\spe{H}$-structure $\spe{h}$?
In Section \ref{sec:coloringcomplex}, we define an abstract relative coloring complex; see Definition \ref{def:abstractcoloring}. The motivation comes from the observing another property of coloring complexes of graphs $\Phi(G)$. Given two flags 
\[F_1 \subsetneq \cdots \subsetneq F_k \]
and \[G_1 \subsetneq \cdots \subsetneq G_m \] 
in $\Phi(G)$ with $F_i = G_j$ for some $i$ and $j$, then 
\[F_1 \subsetneq \cdots \subsetneq F_i \subsetneq G_{j+1} \subsetneq \cdots \subsetneq G_m \]  is also a face of $\Phi(G)$. We refer to this new face as an \emph{exchange}, and say that a complex that is closed under exchanges satisfies the flag exchange condition (see Definition \ref{def:flagexchange}). 
An abstract relative coloring complex on $N$ is a balanced relative subcomplex of the Coxeter complex $\Sigma_N$ of dimension $|N|-2$ satisfying the flag exchange condition. We show that, if $\spe{H}$ is a linearized combinatorial Hopf monoid with a convex character $\varphi$, then the $\varphi$-coloring complex of an $\spe{H}$-structure $\spe{h}$ is an abstract relative coloring complex. Then we show that the collection of abstract relative coloring complexes forms a pointed set species $\spe{C}$. Moreover, $\mathbb{K}\spe{C}$ is also a linearized combinatorial Hopf monoid. In Theorem \ref{thm:terminal}, we show that $\mathbb{K}\spe{C}$ is the terminal object in the category of linearized combinatorial Hopf monoids with convex characters. In particular, the map $\Phi: \spe{H} \to \spe{C}$ that sends $\spe{h}$ to $\Phi(\spe{h})$ is a geometric realization, and gives rise to a Hopf monoid homomorphism $\mathbb{K}(\Phi): \mathbb{K}\spe{H} \to \mathbb{K}\spe{C}$. 
This will be a foundation for future papers: in order to prove a geometric result about realizations of Hopf monoids, it suffices to prove the result for the Hopf monoid of coloring complexes, or some corresponding Hopf submonoid.

Throughout the paper we study familiar examples of graphs, posets, and matroids. 
In Section \ref{sec:applications}, we focus on several other examples. In Subsection \ref{subsec:graphs}, we study a new character for graphs that was recently studied by Aval, Bergeron, and Machacek \cite{aval-bergeron-machacek}. It forms the only prior example of a character we have found that is convex but not balanced. Then we classify all the balanced convex characters of graphs. 

Then we study examples of Hopf monoids related to mixed graphs, double posets, rooted connected graphs, antimatroids, and generalized permutohedra. 
In Subsection \ref{subsec:mixed}, we study a Hopf monoid related to mixed graphs, along with two balanced convex characters. Mixed graphs were introduced in \cite{beck-et-al}. In Subsection \ref{subsec:rooted}, we study a Hopf monoid related to rooted connected graphs. Here the invariant appears to be entirely new. In Subsection \ref{subsec:double}, we study a Hopf monoid related to double posets. The Hopf algebra of double posets was introduced by Malvenuto and Reutenauer \cite{malvenuto}, and one of our quasisymmetric functions was previously studied by Grinberg \cite{grinberg}. The examples of mixed graphs and double posets are interesting because the resulting $\varphi$-chromatic quasisymmetric functions are not always $F$-positive. Thus our results about quasisymmetric functions being $M$-increasing does yield non-trivial results in some cases. 
 In Subsection \ref{subsec:antimatroid}, we study a Hopf monoid related to antimatroids. Antimatroids were introduced in \cite{antimatroid}. 

Finally, in Subsection \ref{subsec:permutohedra}, we study the Hopf monoid of generalized permutohedra that was introduced by Aguiar and Ardila \cite{aguiar-ardila}. 
Many combinatorial Hopf monoids in the literature have morphisms to the Hopf monoid of extended generalized permutohedra (examples are discussed in detail in \cite{aguiar-ardila}). This is an interesting phenomenon that is addressed in \cite{ardila-other}. This motivates us to wonder how coloring complexes compare to generalized permutohedra. We show that the Hopf monoid of antimatroids \emph{do not} have injective linearized morphisms to the Hopf monoid of generalized permutohedra, despite having injective morphisms to the Hopf monoid of coloring complexes. Thus we see that coloring complexes are a more general object that contain interesting submonoids.

\section{Balanced Relative Simplicial Complexes}
\label{sec:geometry}
In this section, we discuss balanced relative simplicial complexes. We prove some new results about their flag $f$-vectors, and we show that the $f$-vector of a pure relative simplicial complex is strongly flawless.

\subsection{Relative Simplicial Complexes and face vectors}
\label{subsec:relsimple}

\begin{definition}
A \emph{relative simplicial complex} on a vertex set $S$ is a collection $\Phi$ of subsets of $S$ with the following property:
\begin{itemize}
    \item For every $\rho \subseteq \sigma \subseteq \tau$, if $\rho, \tau \in \Phi$, then $\sigma \in \Phi$.
\end{itemize}
\end{definition}
The name comes from the fact that there exists simplicial complexes
$(\Gamma, \Sigma)$ with $\Gamma \subseteq \Sigma$, and $\Phi = \Sigma \setminus \Gamma$. We do \emph{not} require that every vertex appear in some face. The \emph{void complex} on $S$ is $\Phi = \emptyset$, while the \emph{empty complex} on $S$ is $\Phi = \{\emptyset \}$. To be consistent with later algebraic structures, we shall denote the void complex by $0_S$ and the empty complex by $1_S$.

Given $\Phi$ with vertex set $S = \{s_1, \ldots, s_k\}$, we let $\mathbb{C}[x_{s_1}, \ldots, x_{s_k}]$ be the polynomial ring with indeterminates $x_{s_1}, \ldots, x_{s_k}$. The \emph{Stanley-Reisner ideal} for $\Sigma$, denoted $I_{\Sigma}$, is generated by $\langle \prod_{s \in \sigma} x_s : \sigma \not\in \Sigma \rangle$, and the \emph{Stanley-Reisner module} for $(\Gamma,\Sigma)$ is $I_{\Gamma} / I_{\Sigma}$. The module is graded by total degree, and its Hilbert function $\Hilb(\Phi, n)$ is the number of monomials of degree $n$ in the module. It is known that the Hilbert function is eventually a polynomial in $n$, and that it is an invariant of $\Phi$, despite being defined in terms of $(\Gamma, \Sigma)$: details can be found in \cite{stanley-algebra}. 

The dimension of a face $\sigma \in \Phi$ is $|\sigma| - 1$, and the dimension of $\Phi$ is the maximum dimension of a face of $\Phi$. By convention, the dimension of the void simplicial complex, $0_S$, is $-\infty$. A maximal face is called a \emph{facet}. A complex is \emph{pure} if all the facets have the same dimension. Note that if we write $\Phi = \Sigma \setminus \Gamma$ for a pair of simplicial complexes $\Gamma \subseteq \Sigma$, then the facets of $\Phi$ are the facets of $\Sigma$ that are \emph{not} contained in $\Gamma$.

\begin{examples}
For the following examples, we use $abc$ to denote the set $\{a,b,c\}$.
Let $\Phi$ be the relative simplicial complex on $\{a,b,c,d,e\}$ given by faces $\{e, ad, de, be, bc, ce, ae, ade, cde, bce, abe \}$. This relative simplicial complex appears on the left in Figure \ref{fig:balancedcomplexes}. In our figures for relative simplicial complexes, we use dashed lines to indicate subsets $\sigma$ of a face $\tau \in \Phi$ where $\sigma \not\in \Phi$. For instance, since $\{a,d\} \in \Phi$, there is a solid edge between $a$ and $d$. Since $\{a\} \subseteq \{a,d\}$ but $\{a\} \not\in \Phi$, there is a dashed circle around the vertex $a$.

If we let $\Xi$ be the relative simplicial complex on $\{a,b,c,d\}$ with faces $\{d, ad, bc, bd, cd, abd, acd, bcd \}$, then $\Xi$ appears on the right in Figure \ref{fig:balancedcomplexes}.
\label{ex:balancedcomplexes}
\end{examples}

\begin{figure}
\begin{center}
\begin{tabular}{cc}
\begin{tikzpicture}
\draw[color=white, fill=gray!20] (30:2cm) -- (330:2cm) -- (210:2cm) -- (150:2cm) -- cycle;
  \node[circle, draw=red, fill=white, dashed, thick] (b) at (150:2cm) {$a$};
  \node[circle, draw=red, fill=white, dashed, thick] (a) at (210:2cm) {$b$};
  \node[circle, draw=red, fill=white,dashed, thick] (e) at (330:2cm) {$c$};
  \node[circle, draw=red, fill=white,dashed, thick] (d) at (30:2cm) {$d$};
  \node[circle, draw=black, fill=white,thick] (or) at (0:0cm) {$e$};
  \draw[dashed, red, thick] (d) -- (e);
  \draw[dashed, red, thick] (a) -- (b);
  \draw[thick] (b) -- (d);
  \draw[thick] (a) -- (e);
  \draw[thick] (b) -- (or) -- (d);
  \draw[thick] (a) -- (or) -- (e);

\end{tikzpicture}
&

\begin{tikzpicture}
\draw[color=white, fill=gray!20] (90:2cm) -- (330:2cm) -- (210:2cm) -- cycle;
\node[circle, draw=red, fill=white, dashed, thick] (c) at (90:2cm) {$a$};
  \node[circle, draw=red, fill=white, dashed, thick] (a) at (210:2cm) {$b$};
  \node[circle, draw=red, fill=white,dashed, thick] (e) at (330:2cm) {$c$};
  \node[circle, draw=black, fill=white,thick] (or) at (0:0cm) {$d$};
  
  \draw[dashed, red, thick] (a) -- (c) -- (e);
  \draw[thick] (a) -- (e);
  \draw[thick] (c) -- (or);
  \draw[thick] (a) -- (or) -- (e);

\end{tikzpicture}
\end{tabular}
\end{center}
\caption{Two examples of pure relative simplicial complexes.}
\label{fig:balancedcomplexes}
\end{figure}

Now we discuss inequalities involving the $f$-vectors of relative simplicial complexes. Our proofs are similar to ones appearing in work of Hibi \cite{hibi}. Suppose $\Phi$ is pure and has dimension $d-1$, where $d$ is a nonnegative integer.

For $i \geq -1$, let $f_i(\Phi)$ be the number of faces of $\Phi$ of dimension $i$. The $f$-vector is defined by $f(\Phi) = (f_{-1}(\Phi), \ldots, f_{d-1}(\Phi))$. Recall that, if $\emptyset \in \Phi$, then $\emptyset$ has dimension $-1$, and this is why we start the $f$-vector with $f_{-1}(\Phi)$.
First, it is well-known (\cite{fvector}) that we have \[\Hilb(\Phi,n) = \begin{cases} f_{-1}(\Phi) & n = 0 \\ \sum\limits_{i=1}^{d} f_{i-1}(\Phi) \binom{n-1}{i-1} & n > 0. \end{cases}\]
Thus, determining the $f$-vector determines the Hilbert function.

Recall that the $f$-vector is \emph{strongly flawless} if the following two systems of inequalities are satisfied:
\begin{enumerate}
    \item $f_{-1} \leq f_0 \leq \cdots \leq f_{\lfloor (d-1) / 2 \rfloor}$.
    \item $f_{i-1} \leq f_{d-i-1}$ for $0 \leq i \leq \lfloor d / 2 \rfloor$.
\end{enumerate}
This is just a translation of the definition of strongly flawless sequence to the sequence $(f_{-1}, \ldots, f_d)$, to accommodate the fact that we start the indexing at $-1.$
A sequence that only satisfies the second set of inequalities is called \emph{flawless}.
\begin{example}
Let $\Phi$ be the relative simplicial complex in Example \ref{ex:balancedcomplexes}. Then $f(\Phi) = (0, 1, 6, 4).$ 

If we let $\Xi$ be the relative simplicial complex in Example \ref{ex:balancedcomplexes}, then $f(\Xi) = (0, 1, 4, 3).$ Both $f$-vectors are strongly flawless. However, if we modified $\Xi$ by adding four new disjoint, isolated vertices, the new $f$-vector would be $(0,5,4,3)$, which is no longer strongly flawless.
\end{example}

Similarly, a sequence $(f_{-1}, \ldots, f_{d-1})$, where the indexing starts at $-1$, is \emph{super flawless} if for all $i$ we have \[(d-i)f_{i-1} \leq (i+1)f_i.\]
\begin{proposition}
Let $(f_{-1}, \ldots, f_{d-1})$ be super flawless. Then $(f_{-1}, \ldots, f_{d-1})$ is strongly flawless.
\label{prop:flawlesspre}
\end{proposition}
\begin{proof}
First, we show that $(f_{-1}, \ldots, f_{d-1})$ is flawless. Since our sequence is super flawless, we see that, for $j \geq 1$, we have $(d-i)(d-i-1)\cdots(d-i-j+1) f_{i-1} \leq (i+1)(i+2)\cdots(i+j)f_{i+j-1}$. Hence $\binom{d-i}{j} f_{i-1} \leq \binom{i+j}{j} f_{i+j-1}$. If we set $j = d-2i$, we obtain $\binom{d-i}{d-2i}f_{i-1} \leq \binom{d-i}{d-2i}f_{d-i-1}$, and thus $f_{i-1} \leq f_{d-i-1}$. Hence the sequence is flawless.

Now we show that it is strongly flawless.
Let $i \leq \lfloor (d-1)/2 \rfloor.$ Then $d \geq 2i+1$, so $(d-i) \geq (i+1)$. Thus $(i+1)f_{i-1} \leq (d-i)f_{i-1} \leq (i+1)f_i$, and thus $f_{i-1} \leq f_i$. \end{proof}

\begin{proof}[Proof of Proposition \ref{prop:flawless}]
 Let $\{F_1, \ldots, F_m \}$ be the set of facets of $\Phi$, and let $V$ be the vertex set of $\Phi$. If $m=0$, then there are no facets, and $\Phi = \Sigma \setminus \Gamma$ for a pair of simplicial complexes where $\Sigma = \Gamma$, and thus $f_i(\Phi) = 0$ for all $i$. Then the $f$-vector is trivially super flawless. So we suppose $m \geq 1$. Let $\Sigma = \langle F_1, \ldots, F_m \rangle$, and let $\Gamma = \Sigma \setminus \Phi$. We prove that $f(\Phi)$ is super flawless by induction on $m$.
 
 First, suppose $m=1$. Then $\Sigma$ is a simplex. Let $i \leq d$. We see that $(d-i)f_{i-1}(\Phi)$ is the cardinality of the set $A_i(\Phi) = \{(S, x): S \in \Phi, |S| = i, x \in V \setminus S \}$. Meanwhile, $(i+1)f_i(\Phi)$ is the cardinality of the set $B_i(\Phi) = \{(T, x): T \in \Phi, |T| = i+1, x \in T \}$. The function $\psi_i: A_i \to B_i$ given by $\psi_i(S, x) = (S \cup \{x\}, x)$ is well-defined because we know that $S \cup \{x\} \in \Phi$ whenever $S \in \Phi$ and $x \in V \setminus S$, because $S \subseteq S \cup \{x\} \subseteq F_1 \in \Phi$. Clearly $\psi_i$ is injective. Thus $(d-i)f_{i-1}(\Phi) \leq (i+1)f_i(\Phi).$ Hence $f(\Phi)$ is super flawless.

 Now suppose $m > 1$. Let $\Phi_1 = \langle F_m \rangle \cap \Phi$, and let $\Phi_2 = \Phi \setminus \Phi_1$. Clearly $\Phi_1$ is a pure relative simplicial complex. Every face in $\Phi$ that is contained in $F_m$ is a face of $\Phi_1$. Hence every face in $\Phi_2$ must be contained in some facet in $\{F_1, \ldots, F_{m-1} \}$, as otherwise it would be a face of $\Phi_1$. Thus $\Phi_2$ is a pure relative simplicial complex. We see that $f_i(\Phi) = f_i(\Phi_1) + f_i(\Phi_2)$. By induction, both $f(\Phi_1)$ and $f(\Phi_2)$ are both super flawless, and the sum of two super flawless sequences of the same dimension is also super flawless. Thus, $f(\Phi)$ is super flawless. \end{proof}

\subsection{Balanced and Colored Relative Simplicial Complexes}
\label{subsec:balanced}

\begin{figure}
\begin{center}
\begin{tabular}{cc}
\begin{tikzpicture}
\draw[color=white, fill=gray!20] (30:2cm) -- (330:2cm) -- (210:2cm) -- (150:2cm) -- cycle;
  \node[circle, draw=red, fill=white, dashed, thick] (b) at (150:2cm) {$1$};
  \node[circle, draw=red, fill=white, dashed, thick] (a) at (210:2cm) {$2$};
  \node[circle, draw=red, fill=white,dashed, thick] (e) at (330:2cm) {$1$};
  \node[circle, draw=red, fill=white,dashed, thick] (d) at (30:2cm) {$2$};
  \node[circle, draw=black, fill=white,thick] (or) at (0:0cm) {$3$};
  \draw[dashed, red, thick] (d) -- (e);
  \draw[dashed, red, thick] (a) -- (b);
  \draw[thick] (b) -- (d);
  \draw[thick] (a) -- (e);
  \draw[thick] (b) -- (or) -- (d);
  \draw[thick] (a) -- (or) -- (e);

\end{tikzpicture}
&

\begin{tikzpicture}
\draw[color=white, fill=gray!20] (90:2cm) -- (330:2cm) -- (210:2cm) -- cycle;
\node[circle, draw=red, fill=white, dashed, thick] (c) at (90:2cm) {$1$};
  \node[circle, draw=red, fill=white, dashed, thick] (a) at (210:2cm) {$2$};
  \node[circle, draw=red, fill=white,dashed, thick] (e) at (330:2cm) {$3$};
  \node[circle, draw=black, fill=white,thick] (or) at (0:0cm) {$4$};
  
  \draw[dashed, red, thick] (a) -- (c) -- (e);
  \draw[thick] (a) -- (e);
  \draw[thick] (c) -- (or);
  \draw[thick] (a) -- (or) -- (e);

\end{tikzpicture}
\end{tabular}
\end{center}
\caption{Two examples of colored relative simplicial complexes.}
\label{fig:balancedcomplexes2}
\end{figure}

A \emph{colored relative simplicial complex} is a triple $(\Phi, k, \kappa)$, where $\Phi$ is a relative simplicial complex, $k \in \mathbb{N}$, and $\kappa: S \to [k]$ is a coloring of vertices of $\Phi$ with the property that no two vertices of a face of $\Phi$ receive the same color. We shall call such a function $\kappa$ a \emph{geometric coloring}, to distinguish from other uses of the word coloring that will appear in this paper.

For a colored relative simplicial complex $\Phi$ with a fixed geometric coloring $\kappa$, we can define the \emph{flag $f$-vector}. Given $S \subseteq [k]$, we let $f_S(\Phi, \kappa)$ denote the number of faces $\{i_1, \ldots, i_m \}$ such that $\{\kappa(i_1), \ldots, \kappa(i_m) \} = S$. Often we will write $f_S(\Phi)$ instead of $f_S(\Phi, \kappa)$, as we usually will be working with various relative simplicial complexes with respect to a fixed vertex set $V$ and geometric coloring $\kappa: V \to [k]$. Given $S = \{i_1, \ldots, i_m \} \subseteq [k]$, with $i_1 < i_2 < \cdots < i_m$, we define the integer composition $\alpha(S) \models (k+1)$ by $\alpha(S) = (i_1, i_2-i_1, \ldots, i_m - i_{m-1}, k+1-i_m)$. We encode the flag $f$-vector with a quasisymmetric function of degree $k+1$, by defining
\[ \dHilb(\Phi, \mathbf{x})= \sum_{S \subseteq [k]} f_S(\Phi) M_{\alpha(S)}.\]

A flag $f$-vector is \emph{weakly increasing}
 if $f_S(\Phi) \leq f_T(\Phi)$ for all $S \subseteq T \subseteq [k]$. Equivalently, a flag $f$-vector is weakly increasing if $\dHilb(\Phi, \mathbf{x})$ is $M$-increasing.
 
 \begin{example}
 Let $P$ be a graded poset with bottom element $\hat{0}$ and top element $\hat{1}$, and let $\Delta(\bar{P})$ be the order complex of $P \setminus \{\hat{0}, \hat{1} \}$, which consists of sets $\{x_1, \ldots, x_k\}$ such that $x_1 < x_2 < \cdots < x_k$ in $P$. Then $\Delta(\bar{P})$ is a balanced simplicial complex, where we define $\kappa(x)$ to be its rank in $P$.
 
 Then \[\dHilb(\Delta(\bar{P}), \mathbf{x}) = \sum_{x_1 < \cdots < x_k} M_{\alpha(\{\kappa(x_1), \kappa(x_2), \ldots, \kappa(x_k)\})}.\]
 
 This is the $F$-quasisymmetric function of a graded poset introduced by Ehrenborg \cite{ehrenborg}. As explained in Example \ref{ex:complexposet}, the $F$-quasisymmetric function of a graded poset is a generalization of the $P$-partition enumerator.
 \label{ex:ehrenborg}
 \end{example}
 
 \begin{example}
 Let $\Phi$ and $\Xi$ be the relative simplicial complexes given in Example \ref{ex:balancedcomplexes}. Define geometric colorings for both complexes as indicated in Figure \ref{fig:balancedcomplexes2}. For example, we have $\kappa(e) = 3$ for $\Phi$. Then we see that both complexes are colored relative simplicial complexes with respect to their respective geometric colorings. 
 
 We calculate $\Hilb(\Phi, \mathbf{x}) = M_{3,1}+2(M_{1,1,2}+M_{1,2,1}+M_{2,1,1})+4M_{1,1,1,1},$ which is $M$-increasing. On the other hand, $f_{\{4\}}(\Xi) = 1$, but $f_{\{1,2,3,4 \}}(\Xi) = 0$. Thus, the flag $f$-vector of $\Xi$ is not weakly increasing.
 \label{ex:balancedcomplexes2}
 \end{example}
 
 Suppose that $\Phi = 0_V$ or $\Phi = 1_V$, the void or empty complexes. Then regardless of the vertex set $V$, or the function $\kappa:V \to [k]$, we see that $(\Phi, \kappa)$ is a colored relative simplicial complex.
 
 Naturally, there is also a flag $h$-vector, given by:
 \[f_S(\Phi) = \sum_{T \subseteq S} h_T(\Phi) \mbox{ or } h_S(\Phi) = \sum_{T \subseteq S} (-1)^{|S| - |T|} f_T(\Phi).\]
 The flag $h$-vector also depends upon the geometric coloring function $\kappa$.
 It holds that $\dHilb(\Phi, \mathbf{x}) = \sum_{S \subseteq [k]} h_S(\Phi) F_{\alpha(S)}$, so the flag $h$-vector expresses $\Phi$ in the basis of fundamental quasisymmetric functions.
 
 We say a flag $h$-vector is \emph{positive} if $h_S(\Phi) \geq 0$ for all $S \subseteq 0$. Suppose that the flag $h$-vector of $\Phi$ is positive, and let $S \subseteq T \subseteq [k]$. Then \[f_S(\Phi) = \sum_{R \subseteq S} h_R(\Phi) \leq \sum_{R \subseteq T} h_R(\Phi) = f_T(\Phi).\] Hence, if the flag $h$-vector of $\Phi$ is positive, then the flag $f$-vector of $\Phi$ is weakly increasing. The converse does not necessarily hold. For example, for $\Phi$ in Example \ref{ex:balancedcomplexes2}, we have $h_{\{1,2,3\}}(\Phi) = -1$, even though the flag $f$-vector is weakly increasing.

A colored relative simplicial complex $\Phi$ with geometric coloring $\kappa: V(\Phi) \to [k]$ is \emph{balanced} if $\Phi$ is pure, and $k=d$, where $d-1$ is the dimension of the complex. The definition of balanced we give here is equivalent to the definition of totally balanced in \cite{fvector}. However, it has become more common to use the phrase \emph{balanced} in place of totally balanced. We show that balanced relative simplicial complexes form a class of colored simplicial complexes whose flag $f$-vectors are weakly increasing.

\begin{examples}
We regard the empty complex $1_{\emptyset}$ with coloring $\kappa:\emptyset \to [0]$ as a balanced simplicial complex over $\emptyset$ of dimension $-1$. We regard the void complex $0_N$ as not being balanced on any ground set $N$. The colored simplicial complex $\Phi$ appearing in Example \ref{ex:balancedcomplexes2} is balanced. The colored complex $\Xi$ in that same example is not. In fact, since every pair of vertices of $\Xi$ is contained in some facet of $\Xi$, we see that any geometric coloring of $\Xi$ must use four colors, so there is no geometric coloring $\kappa$ for which $\Xi$ is balanced.
\label{ex:balancedcomplexes3}
\end{examples}


\begin{proof}[Proof of Proposition \ref{prop:increasing}]

We prove that the flag $f$-vector of $\Phi$ is weakly increasing by induction on the number of facets of $\Phi$. Let $F_1, F_2, \ldots, F_m$ be the facets of $\Phi$. 

First, suppose $m=1$. $\Phi$ has exactly one facet $\sigma$. Then $f_{[d]}(\Phi) = 1$. In fact, by replacing vertex $i$ with $f(i)$, we can assume that our simplicial complex has vertex set $[d]$, and $\sigma = [d]$. We also see that $f_S \leq 1$ for all $S \subseteq [d]$.

Let $S \subseteq T$. If $f_S = 0$, then $f_S \leq f_T$. So suppose $f_S = 1$. Then $S \in \Phi$. Since $[d]$ is also an element of $\Phi$ and $\Phi$ is a relative simplicial complex, we conclude that $T \in \Phi$. Hence $f_T = 1$.

Now assume $m > 1$. Let $\Phi_1 = \Phi \cap \langle F_m \rangle$, and let $\Phi_2 = \Phi \setminus \Phi_1$. Then the flag $f$-vector satisfies the following:
\[f_S(\Phi) = f_S(\Phi_1) + f_S(\Phi_2) \]
for every $S \subseteq [d]$. Moreover, both $\Phi_1$ and $\Phi_2$ are balanced. By induction, the flag $f$-vectors of $\Phi_1$ and $\Phi_2$ are both weakly-increasing. Thus the flag $f$-vector of $\Phi$ is also weakly increasing. \end{proof}
As a corollary, the $F$-quasisymmetric function of a graded poset is $M$-increasing. The fact that the flag $f$-vector of a graded poset is weakly increasing can also be proven using Theorem 2.1 of \cite{billera-hetyei}.

Finally, we mention the relationship between $\dHilb(\Phi, \mathbf{x})$ and $\Hilb(\Phi, n)$. Recall that the \emph{double cone} of $\Phi$ is the relative simplicial complex $\dcone(\Phi)$ with vertex set $V(\Phi) \cup \{-\infty, \infty \}$ given by:
\[\dcone(\Phi) = \{\sigma \cup \tau: \sigma \in \Phi, \tau \subseteq \{-\infty, \infty \} \} \]
Given a quasisymmetric function $f$, the \emph{principal specialization} $\ps_n(f)$ is obtained by setting $x_i = 1$ for $i \leq n$, and $x_i = 0$ otherwise. The invariant $\ps_n(f)$ is a polynomial in the variable $n$.
\begin{lemma}
Let $\Phi$ be a colored relative simplicial complex. Then
\[\ps_{n+1} \dHilb(\Phi, \mathbf{x}) = \Hilb(\dcone(\Phi), n).\]
\label{lem:doublecone}
\end{lemma}
\begin{proof}
Observe that $\ps_{n+1} M_{\alpha} = \binom{n+1}{\ell(\alpha)}$. Thus, 
\begin{align*}\ps_{n+1} \dHilb(\Phi, \mathbf{x}) & = \sum_{S \subseteq [d]} f_S(\Phi) \binom{n+1}{|S|+1} \\ & = \sum_{k=0}^d f_{k-1}(\Phi) \binom{n+1}{k+1}.\end{align*}
In particular, $\ps_{n+1} \dHilb(\Phi, \mathbf{x})$ counts pairs $(\sigma, S)$ where $\sigma \in \Phi$ and $S \subseteq [n+1]$ with $|S| = |\sigma|+1$. We write $\dcone(\Phi) = \Sigma \setminus \Gamma$ for simplicial complexes $\Gamma \subseteq \Sigma$. We give a bijection between the pairs $(\sigma, S)$ and the monomials of degree $n$ in $I_{\Gamma} / I_{\Sigma}$.

Given $\sigma \in \Phi$, we write $\sigma = \{v_1, \ldots, v_k \}$ where $\kappa(v_i) < \kappa(v_j)$ whenever $i < j$. Similarly, we write $S = \{s_1, \ldots, s_{k+1} \}$ with $1 \leq s_1 < s_2 < \cdots < s_{k+1} \leq n+1$. We define $F(\sigma, S) = x_{-\infty}^{s_1-1} x_{\infty}^{n+1-s_{k+1}} \prod_{i=1}^k x_{v_i}^{s_{i+1}-s_{i}}$. 
We see that the result is a monomial of degree $n$. Moreover, since $\sigma \in \Phi$, then for any $T \subseteq \{-\infty, \infty \}$, we have $\sigma \cup T \not \in \Gamma$, and thus $F(\sigma, S) \in I_{\Gamma}$. We also see that $F(\sigma, S) \not\in I_{\Sigma}$, so $F(\sigma, S)$ corresponds to a monomial of degree $n$ in $I_{\Gamma}/I_{\Sigma}$.

Conversely, given a monomial $m$ in $I_{\Gamma}/I_{\Sigma}$ of degree $n$, we may write it in the form $m = x_{-\infty}^{e_0} x_{\infty}^{e_{k+1}} \prod_{i=1}^k x_{v_i}^{e_i}$ where:
\begin{enumerate}
\item $e_0 + \cdots + e_{k+1} = n$, 
\item $e_i > 0$ for $1 \leq i \leq k$,
\item $\{v_1, \ldots, v_k \} \in \Phi$, and 
\item $\kappa(v_i) < \kappa(v_j)$ for $i < j$.
\end{enumerate}
Note that we are allowing $e_0 = 0$ or $e_{k+1} = 0$ as possibilities. Then let $G(m) = (\{v_1, \ldots, v_k \}, \{e_0+1, e_0+e_1+1, \ldots, e_0+\cdots+e_k+1 \})$. Then $F$ and $G$ are inverse bijections. \end{proof}

\subsection{Coxeter Complex of type A and Convex Albums}

In this section, we discuss relative simplicial complexes $\Phi$ that are subcomplexes of the Coxeter complex of type $A$. These form a key class of relative complexes in this paper.

Given a finite set $N$, let $P(N)$ be the collection of proper subsets of $N$. We define the Coxeter complex of the type $A$, denoted $\Sigma_N$, as follows: 
\[\Sigma_N = \{ \{F_1, F_2, \ldots, F_k \}: \emptyset \subset F_1 \subset F_2 \subset \cdots \subset F_k \subset N \} \]
The faces of $\Sigma_N$ are \emph{flags} of proper subsets of $N$. We denote the faces by $F_{\bdot}$ and write $F_{\bdot} : F_1 \subset F_2 \subset \cdots \subset F_k$ to mean $F_{\bdot} = \{F_1, \ldots, F_k \}$ and $F_1 \subset F_2 \subset \cdots \subset F_k$. Given a flag $F_{\bdot} : F_1 \subset F_2 \subset \cdots \subset F_k$, we write $\ell(F_{\bdot}) = k+1$.

There is another combinatorial description of a face of the Coxeter complex. Given a finite set $N$, a \emph{set composition} of $N$ is a sequence $(C_1, \ldots, C_k)$ of disjoint non-empty subsets whose union is $N$. We denote set compositions as $C_1|C_2|\cdots|C_k$, and refer to the sets $C_i$ as \emph{blocks}. We let $C \models N$ to denote that $C$ is a set composition of $N$, and let $\ell(C) = k$ be the length of the composition.

To every set composition $C \models N$, there is an associated flag $F(C) = \{F_1, \ldots, F_{\ell(C) - 1} \}$. We define $F_i = \bigcup_{j=1}^i C_j$. Note that $F(C) \in \Sigma_N$. Similarly, if $F_{\bdot} \in \Sigma_N$, and $\ell(F_{\bdot}) = k$, then there is an associated set composition $C(F_{\bdot})$, defined by:
\begin{enumerate}
    \item $C_1 = F_1$,
    \item $C_i = F_i \setminus F_{i-1}$ for $2 \leq i \leq k-1$, and
    \item $C_k = N \setminus F_{k-1}$.
\end{enumerate}

For example, to the set composition $13|2|45$, the associated edge in the complex $\Sigma_{\{1,2,3,4,5 \}}$ is the flag $\{1,3 \} \subset \{1,2,3 \}$. Hence, we can denote faces of the Coxeter complex by flags of subsets or by set compositions. 

Since there is a partial order on the faces of $\Sigma_N$ given by containment, we also have a partial order on set compositions: we say $C \leq C'$ if every block of $C'$ is a subset of a block of $C$, and, for each block $C_i$ of $C$, the blocks of $C'$ that are contained in $C_i$ appear consecutively in $C'$. This partial order is known as \emph{refinement}. Moreover, the maps $F$ and $C$ are order-preserving bijections.

An \emph{album} $\mathcal{A}$ is a collection of set compositions.
To every album $\mathcal{A}$, we can define a collection of flags $F(\mathcal{A}) = \{F(C): C \in \mathcal{A} \}$. Similarly, given a collection $\mathcal{F}$ of faces of the Coxeter complex $\Sigma_N$, we can associate an album $A(\mathcal{F}) = \{C(F_{\bdot}): F_{\bdot} \in \mathcal{F} \}$. This defines a correspondence between albums and collections of flags. 

\begin{example}
As an example, if we let $\mathcal{F}$ be the collection of \emph{all} flags in $\Sigma_{\{a,b,c\}}$, then 
\begin{align*}
A(\mathcal{F}) &  = \{abc, a|b|c, ab|c, b|a|c, b|ac, b|c|a, bc|a, c|b|a, c|ab, c|a|b, ac|b, a|c|b, a|bc \}. \end{align*}
\end{example}

Now we discuss how to determine if $F(\mathcal{A})$ is a relative simplicial complex.
An album is \emph{convex} if, whenever $C, C'' \in \mathcal{A}$ and there exists $C'$ with $C \leq C' \leq C''$, we have $C' \in \mathcal{A}$. Since the functions $F$ and $A$ are order-preserving, we arrive at the following:

\begin{proposition}
\label{prop:convex}

An album $\mathcal{A}$ is convex if and only if $F(\mathcal{A})$ is a relative simplicial complex.
\end{proposition}

\begin{example}
As an example, let 
\begin{align*} \mathcal{A} & = & \{  a|c|b|d, c|a|b|d, c|a|d|b, a|c|d|b, a|c|bd, ac|b|d, \\ & &   c|a|bd, ac|d|b, a|bc|d, c|ad|b, ac|bd \}.\end{align*}
The reader can check that this is a convex album. Hence $F(\mathcal{A})$ is a relative simplicial complex, which is given on the left in Figure \ref{fig:convexalbum}. In that image, we label the vertex correspond to $\{a,c\}$ as $ac$, and so on. Let $\mathcal{A}' = \mathcal{A} \cup \{abc|d \}$. Then $\mathcal{A}'$ is no longer convex, as $abc|d \leq c|ab|d \leq c|a|b|d$, but $c|ab|d \not \in \mathcal{A}'$. The set $F(\mathcal{A}')$ appears on the right in Figure \ref{fig:convexalbum}. We see that there is a vertex $\{a,b,c\}$ contained in a triangle, and yet contained in a dashed edge. Hence $F(\mathcal{A}')$ is not a relative simplicial complex.
\end{example}

\begin{figure}
\begin{center}
\begin{tabular}{cc}
\begin{tikzpicture}
\draw[color=white, fill=gray!20] (30:2cm) -- (330:2cm) -- (210:2cm) -- (150:2cm) -- cycle;
  \node[circle, draw=red, fill=white, dashed, thick] (b) at (150:2cm) {$abc$};
  \node[circle, draw=red, fill=white, dashed, thick] (a) at (210:2cm) {$c$};
  \node[circle, draw=red, fill=white,dashed, thick] (e) at (330:2cm) {$acd$};
  \node[circle, draw=red, fill=white,dashed, thick] (d) at (30:2cm) {$a$};
  \node[circle, draw=black, fill=white,thick] (or) at (0:0cm) {$ac$};
  \draw[dashed, red, thick] (d) -- (e);
  \draw[dashed, red, thick] (a) -- (b);
  \draw[thick] (b) -- (d);
  \draw[thick] (a) -- (e);
  \draw[thick] (b) -- (or) -- (d);
  \draw[thick] (a) -- (or) -- (e);

\end{tikzpicture}
&

\begin{tikzpicture}
\draw[color=white, fill=gray!20] (30:2cm) -- (330:2cm) -- (210:2cm) -- (150:2cm) -- cycle;
  \node[circle, draw=black, fill=white, thick] (b) at (150:2cm) {$abc$};
  \node[circle, draw=red, fill=white, dashed, thick] (a) at (210:2cm) {$c$};
  \node[circle, draw=red, fill=white,dashed, thick] (e) at (330:2cm) {$acd$};
  \node[circle, draw=red, fill=white,dashed, thick] (d) at (30:2cm) {$a$};
  \node[circle, draw=black, fill=white,thick] (or) at (0:0cm) {$ac$};
  \draw[dashed, red, thick] (d) -- (e);
  \draw[dashed, red, thick] (a) -- (b);
  \draw[thick] (b) -- (d);
  \draw[thick] (a) -- (e);
  \draw[thick] (b) -- (or) -- (d);
  \draw[thick] (a) -- (or) -- (e);

\end{tikzpicture}
\end{tabular}
\end{center}
\caption{Two examples of $F(\mathcal{A})$.}
\label{fig:convexalbum}
\end{figure}

\section{Combinatorial Hopf monoids}
\label{sec:combhopf}
In this section, we define combinatorial Hopf monoids. First, we define the notion of Hopf monoid in the category of linear species.
A \emph{linear species} is a functor $\spe{F}: Set \to Vec$ from the category of finite sets with bijections to the category of finite dimensional vector spaces over a field $\mathbb{K}$ and linear transformations. Sometimes linear species are called vector species, or just species. In this paper, we also use set species and pointed set species. We adopt the following notation: for a given (linear, pointed set, set) species $\spe{Q}$, and a finite set $N$, we let $\spe{Q}_N$ denote the corresponding vector space, pointed set, or set.

\subsection{Linearized Hopf Monoids}
\label{subsec:linearized}

Given any braided monoidal category $C$, there is a notion of a Hopf monoid object $H$ in $C$, which is an object $H$ together with several morphisms satisfying several axioms. The notion of braided monoidal category is studied in Section 1.1 of \cite{aguiar-mahajan-1}, and the full list of axioms for a Hopf monoid object appears in Section 1.2 of loc. cit.

Let $L$ be the category of linear species, with natural transformations as morphisms. Then $L$ has the structure of a braided monoidal category with respect to the Cauchy product. Given two species $\spe{Q}$ and $\spe{R}$, and a finite set $N$, we define \[(\spe{Q} \cdot \spe{R})_N = \bigoplus_{N = S \sqcup T} \spe{Q}_S \otimes \spe{R}_T\] where $\sqcup$ denotes disjoint union. Given a bijection $\sigma:M \to N$, and $S \subseteq M$, we let $\sigma|_S: S \to \sigma(S)$ be the restriction of $\sigma$ to $S.$ 
We define \[(\spe{Q} \cdot \spe{R})_{\sigma} = \bigoplus_{N \ S \sqcup T} \spe{Q}_{\sigma|_S} \otimes \spe{R}_{\sigma|_T}.
\]
Thus $\spe{Q} \cdot \spe{R}$ is a species, known as the Cauchy product of $\spe{Q}$ and $\spe{R}$. More
terminology about species and the Cauchy product are found in Section 8.1 of \cite{aguiar-mahajan-1}. 

For the category of linear species $L$, Aguiar and Mahajan have described the axioms for a Hopf monoid object in terms of an equivalent set of axioms in Section 8.3 of \cite{aguiar-mahajan-1}. We give some of the structural definition and axioms related to associativity and compatibility for a Hopf monoid object in the category of linear species with respect to the Cauchy product, and defer to Secion 8.3 of \cite{aguiar-mahajan-1}, Section 2 of \cite{aguiar-mahajan-2} or Section 2 of \cite{aguiar-ardila} for the full definition and axioms.
For every pair of disjoint finite sets $M, N$, there are linear transformations $\mu_{M,N}: \spe{H}_M \otimes \spe{H}_N \to \spe{H}_{M \sqcup N}$ and $\Delta_{M,N}: \spe{H}_{M \sqcup N} \to \spe{H}_M \otimes \spe{H}_N$. We refer to $\mu$ as \emph{multiplication} and $\Delta$ as \emph{comultiplication}. We focus only on \emph{connected} species, where $\dim \spe{H}_{\emptyset} = 1$. Recall that $\id_X: \spe{H}_X \to \spe{H}_X$ is the identity function.
We require several axioms, including:
\begin{enumerate}
    \item (associativity) For every triple $L, M$, and $N$ of disjoint sets, we have $\mu_{L, M \sqcup N} \circ (\id_L \otimes \mu_{M, N}) = \mu_{L \sqcup M, N} \circ (\mu_{L,M} \otimes \id_N)$.
    \item (identity) There exists $1 \in \spe{H}_{\emptyset}$ such that, for every finite set $N$, and $x \in \spe{H}_N$, we have $\mu_{\emptyset, N}(1 \otimes x) = x = \mu_{N, \emptyset}(x \otimes 1)$.  
    \item (coassociativity) For every triple $L, M$, and $N$ of disjoint sets, we have $(\id_L \otimes \Delta_{M, N}) \circ \Delta_{L, M \sqcup N} = (\Delta_{L,M} \otimes \id_N) \circ \Delta_{L \sqcup M, N}$.
    \item (counit) For every finite set $N$, and $x \in \spe{H}_N$, we have $\Delta_{\emptyset, N}(x) = 1 \otimes x$ and $\Delta_{N, \emptyset}(x) = x \otimes 1$.
    \item (compatibility) for every quadruple $A,B, C$ and $D$ of disjoint sets, we have \[\Delta_{A \sqcup C, B \sqcup D} \circ \mu_{A \sqcup B, C \sqcup D} = (\mu_{A,C} \otimes \mu_{B,D}) \circ (\id_A \otimes \tau_{B,C} \otimes \id_D) \circ (\Delta_{A,B} \otimes \Delta_{C,D}),\] where $\tau_{B,C}: \spe{H}_B \otimes \spe{H}_C \to \spe{H}_C \otimes \spe{H}_B$ is the linear map satisfying $\tau_{B,C}(x \otimes y) = y \otimes x$.
\end{enumerate}
Note that these are equalities of functions.
We let $\Delta_{L, M, N} = \id_L \otimes \Delta_{M,N} \circ \Delta_{L,M \sqcup N}$, and $\spe{x} \cdot \spe{y} = \mu_{M,N}(\spe{x} \otimes \spe{y}) $. Finally, many naturality conditions are required for the unit, counit, multiplication, and comultiplication maps. We use the naturality condition for comultiplication for some proofs. Given $N = S \sqcup T$, where $\sqcup$ denotes disjoint union, and given a bijection $\sigma: M \to N$, we require \[\Delta_{S,T} \circ \spe{H}_{\sigma} = \spe{H}_{\sigma|_{\sigma^{-1}(S)}} \otimes \spe{H}_{\sigma|_{\sigma^{-1}(T)}} \circ \Delta_{\sigma^{-1}(S), \sigma^{-1}(T)}.\] 

\begin{remark}
Our counit axiom is not the usual one presented in \cite{aguiar-mahajan-1}, but for connected species, our axiom is equivalent. We assume the reader is familiar with the usual notion of Hopf monoid in the category of linear species, and with the Cauchy product. 
If we define $\spe{I}$ by \[\spe{I}_N = \begin{cases} \{\vec{0}\} & N \neq \emptyset \\ \mathbb{K} & N = \emptyset \end{cases} \]
then we obtain the unit for the monoidal category $L$. For every species $\spe{H}$ and finite set $N$, we let \[\rho_N: \spe{H}_N \to \bigoplus_{N = S \sqcup T} \spe{I}_S \otimes \spe{H}_T\] be given by $\rho_N(x) = 1 \otimes x$, viewed as an element of $\spe{I}_{\emptyset} \otimes \spe{H}_N$. Then $\rho_N$ is an isomorphism, and $\rho$ is an ismorphism. Similarly, we let \[\lambda_N: \spe{H}_N \to \bigoplus_{N = S \sqcup T} \spe{H}_S \otimes \spe{I}_T\] be given by $\lambda_N(x) = x \otimes 1$, viewed as an element of $\spe{H}_{N} \otimes \spe{I}_{\emptyset}$.Then $\lambda_N$ is an isomorphism, and $\lambda$ is an ismorphism.

The unit of a Hopf monoid is a morphism $\iota: \spe{I} \to \spe{H}$ while the counit is a morphism $\epsilon: \spe{H} \to \spe{I}$. We refer to $\iota_{\emptyset}(1)$ as $1 \in \spe{H}_{\emptyset}.$ One of the axioms for a Hopf monoid object in a monoidal category is that $\iota \epsilon = \id_{\spe{I}}.$ Since $\spe{H}$ is connected, it follows from this axiom that
\[\epsilon_N(x) = \begin{cases} 0 & N \neq \emptyset \\ x & N = \emptyset \end{cases} \] for every finite set $N$ and every $x \in \spe{H}_N$.

The usual counit axiom is that \[ (\epsilon \cdot \id) \circ \Delta = \rho \]
and \[ (\id \cdot \epsilon) \circ \Delta = \lambda. \]

For a given finite set $N,$ and $\spe{h} \in \spe{H}_N$, we see that $(\epsilon \cdot \id) \circ \Delta (\spe{h}) = \Delta_{\emptyset, N}(\spe{h})$, while $\rho(\spe{h}) = 1 \otimes \spe{h}$. Thus $\Delta_{\emptyset, N}(\spe{h}) = 1 \otimes \spe{h}$. A similar argument shows that $\Delta_{N, \emptyset}(\spe{h}) = \spe{h} \otimes 1.$
\end{remark}

Now we discuss the notion of a pointed set species and its linearization. Most of our Hopf monoids come from this construction. 
A \emph{pointed set} consists of a pair $(A, x)$ where $x$ is an element of a set $A$. A morphism $\varphi: (A,x) \to (B, y)$ between pointed sets consists of a function $\varphi: A \to B$ such that $\varphi(x) = y$. An isomorphism is a bijective morphism. Unless we need to be careful, we will refer to the base point as $0$. So all of our pointed sets can be denoted by $(A,0)$. 

Given two pointed sets $(A, 0)$ and $(B, 0)$, their \emph{wedge product}, denoted by $A \wedge B$, is the pointed set  $(A \times B / \sim , (0,0))$, where $\sim$ is the equivalence relation induced by requiring $(a,0) \sim (0,0) \sim (0,b)$ for all $a \in A, b \in B$. 
A \emph{pointed set species} is a functor $\spe{F}: Set \to PSet$ from the category of finite sets with bijections, to the category of pointed finite sets with isomorphisms. A morphism between two pointed set species $\spe{Q}$ and $\spe{R}$ is a natural transformation $\varphi: \spe{Q} \to \spe{R}$.

Given a pointed set species $\spe{F}$, there is an associated linear species $\mathbb{K}\spe{F}$ called the \emph{linearization}: we define $(\mathbb{K}\spe{F})_N$ to be the vector space with basis $\spe{F}_N \setminus \{0\}$. We refer to $\spe{f} \neq 0$ as an $\spe{F}$-structure if there exists a finite set $N$ such that $\spe{f} \in \spe{F}_N \setminus \{0\}$. 

Given a morphism of pointed set species $\varphi: \spe{F} \to \spe{G}$, the \emph{linearization} $\mathbb{K}(\varphi): \mathbb{K}(\spe{F}) \to \mathbb{K}(\spe{G})$ is the natural transformation obtained by extending $\varphi_N: \spe{F}_N \to \spe{G}_N$ by linearity to induce a linear transformation $\mathbb{K}(\varphi)_N: \mathbb{K}(\spe{F})_N \to \mathbb{K}(\spe{G})_N$. We refer to a morphism of linear species as \emph{linearized} if it is the linearization of a morphism of pointed set species. Thus linearization is a functor from the category of pointed set species and natural transformations to $L$.

Let $\spe{H}$ be a pointed set species.
We say that $\mathbb{K}\spe{H}$ is a \emph{linearized} Hopf monoid if $\mathbb{K}\spe{H}$ is a Hopf monoid in the category of linear species, and:
\begin{enumerate}
    \item For every pair of disjoint finite sets $M, N$, the multiplication map $\mu_{M,N}$ on $\mathbb{K}\spe{H}$ is linearized from a morphism $\mu_{M,N}: \spe{H}_M \wedge \spe{H}_N \to \spe{H}_{M \sqcup N}$. This means that for every $\spe{x} \in \spe{H}_M, \spe{y} \in \spe{H}_N$, we have $\spe{x} \cdot \spe{y} \in \spe{H}_{M \sqcup N}$. If $\spe{x} = 0$ or $\spe{y} = 0$, then $\spe{x} \cdot \spe{y} = 0.$
    \item For every pair of disjoint finite sets $M, N$, the comultiplication map $\Delta_{M,N}$ on $\mathbb{K}\spe{H}$ is linearized from a morphism $\Delta_{M,N}: \spe{H}_{M \sqcup N} \to \spe{H}_M \wedge \spe{H}_N$. Thus, for every $\spe{x} \in \spe{H}_{M \sqcup N}$, if $\Delta_{M,N}(\spe{x}) \neq 0$ then there exists $\spe{x}|_M \in \spe{H}_M \setminus \{0\}$ and $\spe{x}/M \in \spe{H}_N \setminus \{0\}$ with $\Delta_{M,N}(\spe{x}) = \spe{x}|_M \otimes \spe{x}/M$.
\end{enumerate}
Our notion of linearization is more general than the usual one found in the literature, because we are working with pointed set species. The advantage is that we view the species of posets as a linearized species, despite the fact that the coproduct is sometimes zero. When we write `let $\mathbb{K}\spe{H}$ be a linearized Hopf monoid', then $\spe{H}$ is the underlying pointed set species.

\begin{proposition}
Let $\spe{H}$ be a pointed set species. Suppose that $\mathbb{K}\spe{H}$ is a linearized Hopf monoid, and that $\mathbb{K}\spe{H}$ is connected. Let $M$ and $N$ be disjoint finite sets with $\spe{x} \in \spe{H}_M$ and $\spe{y} \in \spe{H}_N$. Suppose that $\spe{x} \cdot \spe{y} = 0$. Then $\spe{x} = 0$ or $\spe{y} = 0.$
\end{proposition}
\begin{proof}
Consider $\alpha_{M,N} := \Delta_{M,N} \circ \mu_{M,N}: \spe{H}_M \wedge \spe{H}_N \to \spe{H}_M \wedge \spe{H}_N$. From the compatibility axiom for products and coproducts, we know that \[\Delta_{M,N} \circ \mu_{M,N} = (\mu_{M, \emptyset} \otimes \mu_{\emptyset, N}) \circ (\id_M \otimes \tau_{\emptyset, \emptyset} \otimes \id_N) \circ (\Delta_{M, \emptyset} \otimes \Delta_{\emptyset, N})\] is an equality of functions for the linearization $\mathbb{K}(\alpha_{M,N})$. Hence, we have $\mathbb{K}(\alpha_{M,N})(\spe{x} \otimes \spe{y}) = (\mu_{M, \emptyset} \otimes \mu_{\emptyset, N})(\spe{x} \otimes 1 \otimes 1 \otimes \spe{y}) = \spe{x} \otimes \spe{y}$.

Thus $\alpha_{M,N}(\spe{x}, \spe{y}) = (\spe{x}, \spe{y})$. On the other hand, $\Delta_{M, N}(\spe{x} \cdot \spe{y}) = \Delta_{M,N}(0) = 0$, since $\Delta_{M,N}$ is a function of pointed sets, and hence must preserve basepoints. Thus $(\spe{x},\spe{y}) = 0$, which by definition of wedge product, means that $\spe{x} = 0$ or $\spe{y} = 0$. \end{proof}

Most Hopf monoids that have been studied in the literature are linearized Hopf monoids, the duals of linearized Hopf monoids, or Hadamard products of such Hopf monoids. The Hadamard product is another product on species defined in Section 8.1 of \cite{aguiar-mahajan-1}, and also covered in Section 3 of \cite{aguiar-mahajan-2}. Examples of linearized Hopf monoids include the Hopf monoid of graphs, posets, matroids, hypergraphs, set partitions, linear orders, and generalized permutohedra. In fact, almost every Hopf monoid studied in \cite{aguiar-ardila} is a linearized Hopf monoid. 

Recall that a set species is an endofunctor $\spe{Q}$ on the category of finite sets with bijections. A morphism between two set species $\spe{Q}$ and $\spe{R}$ is a natural transformation.
There is a functor from the category of set species with natural transformations to the category of pointed set species. Given a set species $\spe{Q}$, and a finite set $N$, we let $\spe{Q}^{\bdot}_N = \spe{Q}_N \sqcup \{0\}$. Given a bijection $\sigma:M \to N$, we let $\spe{Q}^{\bdot}_{\sigma}: \spe{Q}^{\bdot}_M \to \spe{Q}^{\bdot}_N$ be given by: \[ \spe{Q}^{\bdot}_{\sigma}(x) = \begin{cases} 0 & x = 0 \\ \spe{Q}_{\sigma}(x) & x \neq 0 \end{cases} \]
We refer to $\spe{Q}^{\bdot}$ as the \emph{pointing} of $\spe{Q}$.
Most of our linearized species will be of two forms:
\begin{enumerate}
    \item They will be of the form $\mathbb{K}\spe{H}^{\bdot}$, where $\spe{H}$ is a set species.
    \item They will be of the form $\mathbb{K}\spe{K}$, where the elements of $\spe{K}_N$ are relative simplicial complexes, and the base point of $\spe{K}_N$ is the void complex $0_N$. Examples of pointed set species of this form appear in Section \ref{sec:realization}.
\end{enumerate}
Given a natural transformation of set species $\varphi: \spe{Q} \to \spe{R}$, we can also define a morphism $\varphi^{\bdot}: \spe{Q}^{\bdot} \to \spe{R}^{\bdot}$ by \[\varphi^{\bdot}_N(x) = \begin{cases} 0 & x = 0 \\ \varphi(x) & \mbox{ otherwise} \end{cases} \]
for every finite set $N$ and every $x \in \spe{Q}^{\bdot}_N$. Thus we have a functor from set species and morphisms to pointed set species and morphisms.
We now list several examples of linearized Hopf monoids. These examples have all appeared in \cite{aguiar-mahajan-1}.

\begin{example}
Given a finite set $N$, let $\spe{E}_N = \{1 \}$. This gives rise to a set species, called the \emph{exponential species}. Then $\mathbb{K}\spe{E}^{\bdot}$ is a linearized Hopf monoid, with multiplication maps $\mu_{M,N}(1 \otimes 1) = 1$, and comultiplication maps $\Delta_{M,N}(1) = 1 \otimes 1$.
\end{example}

\begin{example}

Given a finite set $N$, we let $\spe{S}_N$ denote the collection of set compositions with ground set $N$. This forms the set species of set compositions. Given a set decomposition $N = S \sqcup T$, if $C \models S$ is given by $C_1|\cdots|C_k$ and $C' \models T$ is given by $C'_1|\cdots|C'_r$, then their product $C \cdot C'$ is the set composition $C_1|\cdots|C_k|C'_1|\cdots|C'_r$. 
Given a set composition $C \models M \sqcup N$, we let $C|_M = C_1 \cap M | C_2 \cap M | \cdots | C_k \cap M$ and $C / M = C_1 \cap N | \cdots | C_k \cap N$, where it is understood that we remove any empty blocks from the composition. 
For example, $\Delta_{\{1,3,5\}, \{2,4,6\}}(12|35|46) = 1|35 \otimes 2|46$.

Then $\mathbb{K}\spe{S}^{\bdot}$ is a linearized Hopf monoid.
\end{example}

\begin{example}
Given a finite set $N$, let $\spe{G}_N$ denote the collection of graphs with vertex set $N$. Given a bijection $\sigma:M \to N$, and a graph $\spe{g} \in \spe{G}_M$, define $\spe{G}_{\sigma}(\spe{g})$ to be the graph on $N$ with edges $ij$ if and only if $\sigma^{-1}(i)\sigma^{-1}(j)$ is an edge of $\spe{g}$. Then this gives rise to the set species of graphs $\spe{G}$. 

The species of graphs $\mathbb{K}\spe{G}^{\bdot}$ is a linearized Hopf monoid. The product is given by $\spe{g} \cdot \spe{h} = \spe{g} \sqcup \spe{h}$, the disjoint union of graphs. Given a graph $\spe{g}$, and $S \subseteq N$, $\spe{g}|_S$ is the induced subgraph on $S$, and $\spe{g} / S$ is the induced subgraph on $N - S$. 
\end{example}

\begin{example}
Given a finite set $N$, let $\spe{P}_N$ denote all partial orders on $N$. Given a bijection $\sigma:M \to N$, and a partial order $\spe{p} \in \spe{P}_M$, define $\spe{P}_{\sigma}(\spe{g})$ to be the partial order on $N$ given by $x \leq y$ if and only if $\sigma^{-1}(x) \leq_{\spe{p}} \sigma^{-1}(y)$. Then this gives rise to a set species $\spe{P}$, the species of posets. The species of posets $\mathbb{K}\spe{P}^{\bdot}$ is a linearized Hopf monoid. The product is given by $\spe{p} \cdot \spe{R} = \spe{p} \sqcup \spe{R}$, the disjoint union of partial orders.  Given a partial order $\spe{p}$, and $S \subseteq N$, we define $\spe{p}|_S = \spe{p} / S = 0$ if $S$ is not an order ideal of $\spe{p}$. If $S$ is an order ideal, then $\spe{p}|_S$ and $\spe{p} / S$ are the induced subposets on $S$ and $N - S$ respectively. 
\end{example}

\begin{example}
Given a finite set $N$, let $\spe{M}_N$ denote the collection of matroids with ground set $N$. Given a bijection $\sigma:M \to N$, and a matroid $\spe{m} \in \spe{M}_M$, define $\spe{M}_{\sigma}(\spe{m})$ to be the matroid on $N$ where a set $S$ is a basis if and only if $\sigma^{-1}(S)$ is a basis of $\spe{m}$. Then this gives rise to a set species $\spe{M}$, the species of matroids. The species of matroids $\mathbb{K}\spe{M}^{\bdot}$ forms a linearized Hopf monoid. The product is given by the direct sum operation. Given a matroid $\spe{m}$, and $S \subset N$, we define $\spe{m}|_S$ to be the restriction, and $\spe{m} / S$ to be the contraction of matroids. 
\end{example}

\subsection{Characters and Combinatorial Hopf Monoids}
\label{subsec:characters}

Recall that a morphism of (set, pointed set, linear) species $\varphi: \spe{F} \to \spe{G}$ is a natural transformation between $\spe{F}$ and $\spe{G}$, where $\spe{F}$ and $\spe{G}$ are (set, pointed set, linear) species.

A morphism $\varphi: \spe{H} \to \mathbb{K}\spe{E}^{\bdot}$, where $\spe{H}$ is a Hopf monoid in species, is a \emph{character} if for all disjoint finite sets $M$ and $N$, and all $x \in \spe{H}_M$ and all $y \in \spe{H}_Y$, we have $\varphi_N(x) \cdot \varphi_M(y) = \varphi_{M \sqcup N}(x \cdot y)$. By an abuse of notation, we will write $\varphi(\spe{h})$ in place of $\varphi_{N}(\spe{h})$, when no confusion will arise. Given a a linearized Hopf monoid $\mathbb{K}(\spe{H})$, and a morphism $\varphi: \spe{H} \to \spe{E}^{\bdot}$ such that $\mathbb{K}(\varphi)$ is a character, then we shall refer to $\mathbb{K}(\varphi)$ as a linearized character.

Now we define a linearized combinatorial Hopf monoid. A \emph{linearized combinatorial Hopf monoid} is a linearized Hopf monoid $\mathbb{K}\spe{H}$ with a \emph{linearized} character $\mathbb{K}(\varphi)$. This means that $\varphi(\spe{h}) = 0$ or $\varphi(\spe{h}) = 1$ for every $\spe{H}$-structure.
The motivation is that many combinatorial Hopf algebras are studied where the character $\varphi$ only takes on the values $0$ and $1$.

First, we mention two examples of characters that are defined for every linearized Hopf monoid.
\begin{example}
Let $\mathbb{K}\spe{H}$ be a linearized Hopf monoid. Given a finite set $N$, and $\spe{h} \in \spe{H}_N$, we define \[\zeta_N(\spe{h}) = \begin{cases} 0 & \spe{h} = 0 \\ 1 & \mbox{ otherwise } \end{cases} \] for every $\spe{H}$-structure $\spe{h}$. Then $\mathbb{K}(\zeta)$ is a linearized character, which we call the \emph{zeta character}. Thus, for every linearized Hopf monoid $\mathbb{K}\spe{H}$, we see that $(\mathbb{K}\spe{H}, \mathbb{K}(\zeta))$ is a linearized combinatorial Hopf monoid.
\label{def:zeta}
\end{example}

\begin{definition}
Let $\mathbb{K}\spe{H}$ be a linearized Hopf monoid. Given a finite set $N$, we say an $\spe{H}$-structure $\spe{h} \in \spe{H}_N \setminus \{0\}$ is \emph{totally reducible} if $|N| = 1$, or there exists a nontrivial decomposition $N = S \sqcup T$, and totally reducible elements $\spe{x} \in \spe{H}_S$ and $\spe{y} \in \spe{H}_T$ such that $\spe{h} = \spe{x} \cdot \spe{y}$. We define 
    \[ \chi_N(\spe{h}) = 
\begin{cases}
1 & \text{if $\spe{h}$ is totally reducible} \\
0 & \text{ otherwise.}
\end{cases}
\]
We call $\mathbb{K}(\chi)$ the \emph{chromatic} character. Thus, for every linearized Hopf monoid $\mathbb{K}\spe{H}$, we see that $(\mathbb{K}\spe{H}, \mathbb{K}(\chi))$ is a linearized combinatorial Hopf monoid.
\label{def:chromatic}
\end{definition}
 For instance, if we let $\spe{H} = \spe{G}$, then a graph $\spe{g}$ is totally reducible if and only if it is edgeless. If we let $\spe{H} = \spe{P}$, then a poset $\spe{p}$ is totally reducible if and only if it is an antichain. Finally, if $\spe{H} = \spe{M}$, then a matroid $\spe{m}$ is totally reducible if and only if is a direct sum of loops and coloops, which means $\spe{m}$ has a unique basis. These characters were studied in context of Hopf algebras in \cite{aguiar-bergeron-sottile}, and in the context of Hopf monoids in \cite{aguiar-ardila}. These characters give rise to the chromatic symmetric function, Gessel's $P$-partition enumerator, and the Billera-Jia-Reiner invariant of a matroid, as we discuss below.

\subsection{Hopf submonoids}

Finally, given a species $\spe{Q}$, a species $\spe{R}$ is a \emph{subspecies} of $\spe{Q}$ if the following two conditions are satisfied:
\begin{enumerate}
    \item For every finite set $N$, we have $\spe{R}_N \subseteq \spe{Q}_N$.
    \item For every bijection $\sigma:M \to N$, and every $\spe{r} \in \spe{R}_N$, we have $\spe{R}_{\sigma}(\spe{r}) = \spe{Q}_{\sigma}(\spe{r})$.
\end{enumerate}

Suppose that $\mathbb{K}\spe{Q}$ is a linearized Hopf monoid. Then $\mathbb{K}\spe{R}$ is a linearized Hopf submonoid if $\spe{R}$ is subspecies of $\spe{Q}$ such that $\mathbb{K}\spe{R}$ is a linearized Hopf monoid, where the product and coproduct on $\mathbb{K}\spe{R}$ agrees with the corresponding structures on $\mathbb{K}\spe{Q}$. This is equivalent to requiring the following two conditions:
\begin{enumerate}
    \item For every pair of disjoint finite sets $M$ and $N$, and structures $\spe{p} \in \spe{R}_M$ and $\spe{r} \in \spe{R}_N$, we have $\spe{p} \cdot \spe{r} \in \spe{R}_{M \sqcup N}$.
    \item For every pair of disjoint finite sets $M$ and $N$, and every $\spe{R}$ structure $\spe{q} \in \spe{R}_{M \sqcup N}$, if $\Delta_{M,N}(\spe{q}) \neq 0$, then $\spe{q}|_S$ and $\spe{q}/S$ are both non-zero $\spe{R}$ structures.
\end{enumerate}

If $\mathbb{K}\spe{Q}$ is a linearized combinatorial Hopf monoid with linearized character $\mathbb{K}(\varphi)$, and $\mathbb{K}\spe{R}$ is a linearized Hopf submonoid, then $\mathbb{K}\spe{R}$ is a combinatorial Hopf submonoid with character $\mathbb{K}(\varphi)$ as well.

\subsection{The chromatic quasisymmetric function}
\label{subsec:chromaticquasisymmetric}

Recall that a \emph{combinatorial Hopf algebra} is a graded connected Hopf algebra $H$ with finite graded dimension, equipped with with a linear multiplicative map $\varphi: H \to \mathbb{K}$, which is called a character. Aguiar, Bergeron, and Sottile \cite{aguiar-bergeron-sottile} showed that there was a unique Hopf algebra homomorphism from any combinatorial Hopf algebra $H$ to the combinatorial Hopf algebra of quasisymmetric functions. In this section we review a construction that associates a Hopf algebra to a Hopf monoid. The end result is that, given a linearized combinatorial Hopf monoid $\spe{H}$, and an $\spe{H}$-structure $\spe{h}$, there is a quasisymmetric function invariant $\Psi_{\varphi}(\spe{h}, \mathbf{x})$ associated to $\spe{h}$. The advantage of working with species is that we are able to describe this invariant as enumerating \emph{colorings}. We shall refer to this invariant as the $\varphi$-chromatic quasisymmetric function of $\spe{H}$.

In Chapter 15, Section 1 of \cite{aguiar-mahajan-1}, Aguiar and Mahajan define the Full Fock functor $\mathcal{F}$, which sends linear species to graded vector spaces. Given a linear species $\spe{Q},$ we define $\mathcal{F}(\spe{Q}) = \bigoplus_{n \geq 0} \spe{Q}_{[n]}$. This defines a functor.

Let $k$ and $n$ be positive integers. We define $[k+1,k+n] = \{k+i: i \in [n] \}.$
Given two linearly ordered sets $S = \{s_1, \ldots, s_n \}$ and $T = \{t_1, \ldots, t_n \}$ with $s_1 < \cdots < s_n$ and $t_1 < \cdots < t_n$, we define $\cano_S^T: S \to T$ by $\cano(s_i) = t_i.$

Given two species $\spe{Q}$ and $\spe{R}$, we note that $\mathcal{F}(\spe{Q}) \cdot \mathcal{F}(\spe{R})$ is isomorphic to \[\bigoplus_{n \geq 0} \bigoplus_{n = s+t} \spe{Q}_{[s]} \otimes \spe{R}_{[t]}. \] We define \[\varphi_{\spe{Q},\spe{R}}: \mathcal{F}(\spe{Q}) \cdot \mathcal{F}(\spe{R}) \to \mathcal{F}(\spe{Q}\cdot \spe{R}), \] by 
\[\varphi_{\spe{Q}, \spe{R}} = \bigoplus_{n \geq 0} \bigoplus_{s+t=n} \spe{Q}_{\id} \otimes \spe{R}_{\cano_{[t]}^{[s+1,n]}}.\]
We also define 
 \[\psi_{\spe{Q},\spe{R}}: \mathcal{F}(\spe{Q}\cdot \spe{R}) \to \mathcal{F}(\spe{Q}) \cdot \mathcal{F}(\spe{R}),  \]
 by 
 \[\psi_{\spe{Q}, \spe{R}} = \bigoplus_{n \geq 0} \bigoplus_{S \sqcup T = [n]} \spe{Q}_{\cano_{S}^{[|S|]}} \otimes \spe{R}_{\cano_T^{[|T|]}}. \]

 In Proposition 15.25 of \cite{aguiar-mahajan-1}, Aguiar and Mahajan show that, if $\spe{H}$ is a Hopf monoid, then $\mathcal{F}(\spe{H})$ is a Hopf algebra, with multiplication given by $\mathcal{F}(\mu) \circ \varphi_{\spe{H}, \spe{H}}$, and comultiplication given by $\psi_{\spe{H}, \spe{H}} \circ \mathcal{F}(\Delta).$

Suppose that we have a character $\mathbb{K}(\varphi): \mathbb{K}\spe{H} \to \mathbb{K}\spe{E}$. Given a $\spe{H}$-structure $h \in \spe{H}_{[n]}$, we can define a character $\hat{\varphi}: \mathcal{F}(\mathbb{K}\spe{H}) \to \mathbb{K}$ by $\hat{\varphi}(h) = \varphi(h)$. Thus $(\mathcal{F}(\mathbb{K}\spe{H}), \hat{\varphi})$ is a combinatorial Hopf algebra.

We recall the quasisymmetric function associated to a character on a combinatorial Hopf algebra $(H, \varphi)$. Given $h \in H_n$, with $n \geq 1$, \begin{equation} \Psi_{\varphi}(h, \mathbf{x}) = \sum_{k = 1}^n \sum_{\alpha_1 + \alpha_2 + \cdots + \alpha_k = n} (\varphi^{\otimes k} \circ (\rho_{\alpha_1} \otimes \rho_{\alpha_2} \otimes \cdots \otimes \rho_{\alpha_k}) \circ \Delta^{k-1})(h) M_{\alpha} \label{eq:deffunction} \end{equation}  
where $\rho_n$ is the projection map from $H$ to $H_n$, and the summation is over all integer compositions $\alpha$ of $n$. For $h \in H_0$, we define $\Psi_{\varphi}(h, \mathbf{x}) = \epsilon(h)M_{\emptyset}$, where $\epsilon$ is the counit of $H$, and $\emptyset$ is the empty integer composition.
    
Thus, given a linearized combinatorial Hopf monoid $\mathbb{K}\spe{H}$ with linearized character $\mathbb{K}(\varphi)$, there is an associated combinatorial Hopf algebra $(\mathcal{F}(\mathbb{K}\spe{H}), \hat{\varphi})$. For each finite set $N$, let $\sigma: N \to [n]$ be any bijection. To every $\spe{H}$-structure $\spe{h} \in \spe{H}_N$ we define $\Psi_{\varphi}(\spe{h}, \mathbf{x}) = \Psi_{\hat{\varphi}}(\spe{H}_{\sigma}(\spe{h}), \mathbf{x})$, which we call the \emph{$\varphi$-chromatic quasisymmetric function}. We show below that the resulting quasisymmetric function does not depend on the choice of $\sigma$. We will abuse notation and denote the $\varphi$-chromatic quasisymmetric function by $\Psi_{\varphi}(\spe{h}, \mathbf{x})$.

We also define $\chi_{\varphi}(\spe{h}, n) = \ps_n \Psi_{\varphi}(\spe{h}, \mathbf{x})$, where $\ps$ is the principal specialization. Then $\chi_{\varphi}(\spe{h}, n)$ is a polynomial in $n$, which we call the \emph{$\varphi$-chromatic polynomial.}

We give two alternative formulas for $\Psi_{\varphi}(\spe{h}, \mathbf{x})$, both of which demonstrate that $\Psi_{\varphi}(\spe{h}, \mathbf{x})$ is well-defined. Let $\mathbb{K}\spe{H}$ be a linearized Hopf monoid with linearized character $\mathbb{K}(\varphi)$.
 Given a set composition $C \models N$, we define $\Delta_C = (\id_{C_1} \otimes \Delta_{C/C_1}) \circ \Delta_{C_1, N \setminus C_1}$. We also let $\varphi_C = (\bigotimes_{i=1}^k \varphi_{C_i}) \circ \Delta_C$. We use this notation for various proofs. 
 
  Let $N$ be a finite set, and let $\spe{h} \in \spe{H}_N$. We say a set composition is \emph{$\varphi$-proper} for $\spe{h}$ if $\varphi_C(\spe{h}) = 1$. Let $f: N \to \mathbb{N}$ be a function. For $i \in \mathbb{N}$, let $N_i = \{v \in N: f(v) \leq i \}$.
We call $f$ a \emph{$\varphi$-proper coloring} of $\spe{h}$ if $\varphi(\spe{h}|_{N_{i+1}}/ N_{i}) = 1$ for all $i$.

\begin{theorem}
Let $\mathbb{K}\spe{H}$ be a linearized combinatorial Hopf monoid with linearized character $\mathbb{K}(\varphi)$. Fix a finite set $N$, and $\spe{h} \in \spe{H}_N$. 
Then  \[\Psi_{\varphi}(\spe{h}, \mathbf{x}) =  \sum_{\mbox{$C$ is  $\varphi$-proper}} M_{\alpha(C)} = \sum_{\mbox{$f$ is $\varphi$-proper}} \prod_{n \in N} x_{f(n)}. \]

Similarly, $\chi_{\varphi}(\spe{h}, n)$ is the number of $\varphi$-proper colorings $f$ such that $f(N) \subseteq [n]$.
\label{thm:coloring}
\end{theorem}
For example, consider the combinatorial Hopf monoid $(\mathbb{K}\spe{G}, \mathbb{K}(\chi))$. Then the resulting invariant $\Psi_{\varphi}(\spe{g}, \mathbf{x}) = X(G, \mathbf{x}),$ the chromatic symmetric function. The invariant $\Psi_{\zeta}(\spe{p}, \mathbf{x})$ for the combinatorial Hopf monoid $(\mathbb{K}\spe{P}, \mathbb{K}(\zeta))$  enumerates $\spe{p}$-partitions. Finally, for the combinatorial Hopf monoid $(\mathbb{K}\spe{M}, \mathbb{K}(\chi))$, the invariant $\Psi_{\chi}(\spe{m}, \mathbf{x})$ is the Billera-Jia-Reiner quasisymmetric function of a matroid \cite{billera-jia-reiner}.

\begin{proof}

Let $N$ be a finite set, $n = |N|$, and let $\sigma: N \to [n]$ be a bijection. For any $S \subseteq N$, we define $\sigma_S^{[|S|]} = \cano_{\sigma(S)}^{[|S|]} \circ \sigma|_S$. We have \begin{align*} \Delta^{k-1}(\spe{H}_{\sigma}(\spe{h})) & =  \sum_{C_1 \sqcup  C_2 \sqcup \cdots \sqcup C_{k} = [n]}  \left( \bigotimes_{i=1}^{k} \spe{H}_{\cano_{C_i}^{[|C_i|]}} \right)\circ \Delta_{C_1|C_2|\cdots |C_{k}}(\spe{H}_{\sigma}(\spe{h})) \\
&= \sum_{C \models [n]: \ell(C) = k}  \left(\bigotimes_{i=1}^{k} \spe{H}_{\cano_{C_i}^{[|C_i|]}} \circ \spe{H}_{\sigma|_{\sigma^{-1}(C_i)}} \right) \circ \Delta_{\sigma^{-1}(C)}(\spe{h})  \\
&= \sum_{C \models [n]: \ell(C) = k}  \left( \bigotimes_{i=1}^{k} \spe{H}_{\sigma_{\sigma^{-1}(C_i)}^{[|C_i|]}} \right) \circ \Delta_{\sigma^{-1}(C)}(\spe{h}). \end{align*} 
where the first equality is the definition of the coproduct of $\mathcal{F}(\mathbb{K}\spe{H})$ on $\spe{H}_{\sigma}(\spe{h}).$ The second equality involves reindexing the summation with set compositions of length $k$, and applying the naturality of $\Delta_{C}$ to $\spe{H}_{\sigma}.$ The third equality comes from the definition of $\sigma|_{S}^{[|S|]}$ and the fact that $\spe{H}$ is a functor.
Thus we have
\begin{align*}
    \varphi^{\ell(\alpha)} \circ \rho^{\alpha} \circ \Delta^{k-1}(\spe{H}_{\sigma}(\spe{h})) & = 
    \sum_{C \models [n]: \alpha(C) = \alpha}  \left(\bigotimes_{i=1}^{\ell(\alpha)} \varphi_{[|C_i|]} \circ \spe{H}_{\sigma_{\sigma^{-1}(C_i)}^{[|C_i|]}} \right) \circ \Delta_{\sigma^{-1}(C)}(\spe{h}) \\
    &= \sum_{C \models [n]: \alpha(C) = \alpha}  \left( \bigotimes_{i=1}^{\ell(\alpha)} \varphi_{\sigma^{-1}(C)} \right) \circ \Delta_{\sigma^{-1}(C)}(\spe{h}) \\
    &= \sum_{C \models N: \alpha(C) = \alpha}  \varphi_C(\spe{h})
\end{align*}
where the second equality comes from the naturality of $\varphi$, and the third equality involves replacing $C$ with $\sigma^{-1}(C)$.
We see that $\varphi_C(\spe{h}) \neq 0$ only when $C$ is $\varphi$-proper, in which case $\varphi_C(\spe{h}) = 1$. Thus $[M_{\alpha}]\Psi_{\varphi}(\spe{h}, \mathbf{x})$ is the number of $\varphi$-proper compositions of type $\alpha$. The first equality follows.

To prove the second equality, first let $f: N \to \mathbb{N}$ be an arbitrary function. We let $B(f) = \{i \in \mathbb{N}: f^{-1}(i) \neq \emptyset \}$. Then $|B(f)| < \infty$, and we write $B(f) = \{i_1, \cdots, i_k \}$ with $i_1 < i_2 < \cdots < i_k$. Define the set composition $C(f) = C_1 | C_2 | \cdots | C_k$, where $C_j = f^{-1}(i_j)$ and $i_j \in B(f)$. Then we see that \[\mathbf{x}^f = \prod_{v \in N} x_{f(v)} = \prod_{j=1}^k x_{i_j}^{|C_j|}.\] Thus $\mathbf{x}^f$ appears as a term in $M_{\alpha(C)}$. We claim that \[M_{\alpha(C)} = \sum_{f: C(f) = C} \mathbf{x}^f. \] To see this, let $i_1 < \cdots < i_k$. Let $g_{i_1, \ldots, i_k}:N \to \mathbb{N}$ be defined by $g_{i_1, \ldots, i_k}(u) = i_j$ where $u \in C_j$. Then $g_{i_1, \ldots, i_k}$ is a function with $C(g_{i_1, \ldots, i_k}) = C$ and $\mathbf{x}^{g_{i_1, \ldots, i_k}} = \prod_{j=1}^k x_{i_j}^{\alpha_j}$. Thus every term in $M_{\alpha(C)}$ corresponds to a unique function.

Let $f: N \to \mathbb{N}$. We show that $f$ is $\varphi$-proper if and only if $C(f)$ is $\varphi$-proper. Suppose that $f$ is $\varphi$-proper. Then for all $i \in \mathbb{N}$, we see that $\varphi(\spe{h}|_{N_i}/N_{i-1}) = 1$. Let $C(f) = C_1|\cdots|C_k$. Since $C_1 \cup \cdots \cup C_j = N_i$ for some $i$, it follows that $C(f)$ is $\varphi$-proper. Similarly, if $\varphi(\spe{h}|_{N_i}/N_{i-1}) = 0$ for some $i$, then $N_{i-1} \neq N_i$, and if $N_{i-1} = C_1 \cup \cdots \cup C_j$, it follows that $\varphi_C(\spe{h}) = 0$, and $C$ is not $\varphi$-proper. \end{proof}

\section{Geometric Realization}
\label{sec:realization}
Now we introduce \emph{geometric realizations} for a linearized combinatorial Hopf monoid. Given a linearized combinatorial Hopf monoid $(\mathbb{K}\spe{H},  \mathbb{K}(\varphi))$, a geometric realization consists of a pair $(\spe{K}, \Sigma_{\varphi})$, where $\spe{K}$ is a pointed set species, and $\Sigma_{\varphi}: \spe{H} \to \spe{K}$ is a natural transformation such that:
\begin{enumerate}
    \item The structures of the pointed set species $\spe{K}$ are colored relative simplicial complexes $(\Phi, \kappa)$.
    \item For each $N$, the base point of $\spe{K}_N$ is the void complex $0_N.$
    \item The natural transformation $\Sigma_{\varphi}: \spe{H} \to \spe{K}$ satisfies the additional property that $\Psi_{\varphi}(\spe{h}, \mathbf{x}) = \dHilb(\Sigma_{\varphi}(\spe{h}), \mathbf{x})$ for every $\spe{H}$-structure $\spe{h}$.
\end{enumerate}
We briefly mention the motivation behind making the void complex the base point.
Given $0 \in \spe{H}_N$, we see that $\Delta_C(0) = 0$ for all $C$. Thus $\Psi_{\varphi}(0, \mathbf{x}) = 0$ and $\dHilb(\Sigma_{\varphi}(0), \mathbf{x}) = 0$. This forces $\Sigma_{\varphi}(0)$ to be the void complex. Thus the void complex naturally appears to play the same role as the zero vector. 
Our focus will be to describe a geometric realization that comes from a certain construction which relies on the fact that we already have a description of $\Psi_{\varphi}(\spe{h}, \mathbf{x})$ as a summation over set compositions.

\subsection{Convex characters}
\label{subsec:convex}

Given a linearized combinatorial Hopf monoid $(\mathbb{K}\spe{H},  \mathbb{K}(\varphi))$, and an $\spe{H}$-structure $\spe{h} \in \spe{H}_N$, we let \[\mathcal{A}_{\varphi}(\spe{h}) = \{C \models N: \varphi_C(h) = 1 \} \] be the \emph{album associated to} $\spe{h}$. We see that this is a natural album which arises when studying $\Psi_{\varphi}(\spe{h}, \mathbf{x})$. In fact, we have \[\Psi_{\varphi}(\spe{h}, \mathbf{x}) = \sum_{C \in \mathcal{A}_{\varphi}(\spe{h})} M_{\alpha(C)} = \sum_{\substack{F_{\bdot} \in F(\mathcal{A}_{\varphi}(\spe{h})) \\ \mbox{ flags in } N}} M_{\alpha(F_{\bdot})}. \]

Moreover, $F(\mathcal{A}_{\varphi}(\spe{h}))$ is a collection of sets, and if we define $\kappa:2^N \to [|N|-1]$ via $\kappa(S) = |S|$, then it is not hard to show that $\kappa$ is a geometric coloring for $F(\mathcal{A}_{\varphi}(\spe{h}))$. So if $\mathcal{A}_{\varphi}(\spe{h})$ is convex, then $F(\mathcal{A}_{\varphi}(\spe{h}))$ is a colored relative complex. Our next goal is to determine precisely when this happens.

\begin{definition}
Let $\mathbb{K}\spe{H}$ be a linearized Hopf monoid, and let $\mathbb{K}(\varphi)$ be a linearized character.
We say that $\varphi$ is a \emph{convex character} if, for every $\spe{H}$-structure $\spe{h}$ such that $\varphi(\spe{h}) = 1$, and any two set compositions $C \leq C'$ such that $C, C' \models N$, if $\varphi_{C'}(\spe{h}) = 1$, then $\varphi_C(\spe{h}) = 1$.
\label{def:convex}
\end{definition}

The word convex comes from the following proposition.
\begin{proposition}
Let $(\mathbb{K}\spe{H},  \mathbb{K}(\varphi))$ be a linearized combinatorial Hopf monoid. Then the albums $\mathcal{A}_{\varphi}(\spe{h})$ are convex for all $\spe{H}$-structures $\spe{h}$ if and only if $\varphi$ is a convex character. 
\label{prop:convexequiv}
\end{proposition}
\begin{proof}
First, suppose that $\mathcal{A}_{\varphi}(\spe{h})$ is convex for every $\spe{H}$-structure $\spe{h}$. Let $\spe{h}$ be an $\spe{H}$-structure such that $\varphi(\spe{h}) = 1$. Let $C \leq C'$ be two set compositions such that $C, C' \models N$ and $\varphi_{C'}(\spe{h}) = 1$. Then $C' \in \mathcal{A}_{\varphi}(\spe{h}).$ Since $\varphi(\spe{h}) = 1$, $N \in \mathcal{A}_{\varphi}(\spe{h}).$ Since $N \leq C \leq C'$, it follows that $C \in \mathcal{A}_{\varphi}(\spe{h})$, and hence $\varphi_{C}(\spe{h}) = 1$. Therefore $\varphi$ is a convex character.

Now suppose that $\varphi$ is a convex character. Let $N$ be a finite set and let $\spe{h} \in \spe{H}_N$. Let $C \leq C' \leq C'' \models N$, with $C, C'' \in \mathcal{A}_{\varphi}(\spe{h})$. Let $F(C) = S_1 \subset S_2 \subset \cdots \subset S_k$. Since $\varphi_C(\spe{h}) = 1$, we know $\varphi(\spe{h}|_{S_{i-1}}/ S_i) = 1$ for all $i$. We observe that $\varphi_{C''|_{S_{i-1}}/ S_i}(\spe{h}|_{S_{i-1}}/S_i) = 1$ for all $i$. Since $C'|_{S_{i-1}}/ S_i \leq C''|_{S_{i-1}}/ S_i$ and $\varphi$ is a convex character, we see that $\varphi_{C'|_{S_{i-1}}/ S_i}(\spe{h}|_{S_{i-1}}/S_i) = 1$ for all $i$. Then $\varphi_{C'}(\spe{h}) = 1$, so $C' \in \mathcal{A}_{\varphi}(\spe{h})$. \end{proof}

\begin{definition}
Let $(\mathbb{K}\spe{H}, \mathbb{K}(\varphi))$ be a linearized combinatorial Hopf monoid. We say that $\varphi$ is a \emph{balanced convex character} if, for every finite set $N$ and every non-zero $\spe{H}$-structure $\spe{h} \in \spe{H}_N$, the following conditions are satisfied:
\begin{enumerate}
    \item If $|N| = 1$, then $\varphi(\spe{h}) = 1$.
    \item If $|N| > 1$, then there exists non-empty $S \subsetneq N$ such that \[\Delta_{S, N \setminus S}(\spe{h}) \neq 0. \]
    \item If $\varphi(\spe{h}) = 1$, then $\varphi(\spe{h}|_S) = \varphi(\spe{h}/S) = 1$ for any $S \subsetneq N$ such that $\Delta_{S, N \setminus S}(\spe{h}) \neq 0$.
\end{enumerate}
\label{def:balancedconvex}
\end{definition}

\begin{proposition}
Let $(\mathbb{K}\spe{H},  \mathbb{K}(\varphi))$ be a linearized combinatorial Hopf monoid. If $\varphi$ is a balanced convex character, then $\varphi$ is a convex character.
\label{prop:criterion}

\end{proposition}
\begin{proof}
Fix a finite set $N$ with $|N| \geq 1$. We prove the following statement by induction on $|N|$: 

\begin{center} for every $\spe{h} \in \spe{H}_N$ such that $\varphi(\spe{h}) = 1$, and any two set compositions $C \leq C'$ such that $C, C' \models N$, if $\varphi_{C'}(\spe{h}) = 1$, then $\varphi_C(\spe{h}) = 1$. \end{center}
 
 For the base case, assume $|N| = 1$, and let $\spe{h} \in \spe{H}_N$ such that $\varphi(\spe{h}) = 1$, Given any compositions $C \leq C' \models N$, we see that $C = C'$, as both have only one block. Moreover $\varphi_{C}(\spe{h}) = \varphi(\spe{h}) = 1$, so the condition is trivially satisfied.
 
 Suppose $|N| > 1.$
 Let $\spe{h} \in \spe{H}_N$ with $\varphi(\spe{h}) = 1$. Let $C \leq C' \models N$ be two compositions with $\varphi_{C'}(\spe{h}) = 1$. Write $C = C_1|\cdots|C_k$. By properties of coproduct, it follows that 
\[\Delta_{C}(\spe{h}) = (\id_{C_1} \otimes \Delta_{C/C_1}) \circ \Delta_{C_1, N \setminus C_1}(\spe{h}) \neq 0.\] Thus $\Delta_{C_1, N \setminus C_1}(\spe{h}) \neq 0$, and so $\varphi(\spe{h}|_{C_1}) = \varphi(\spe{h}/C_1) = 1$. We see that $C/C_1 \leq C'/C_1 \models N \setminus C_1$. Moreover, we have $\varphi_{C'/C_1}(\spe{h}/C_1) = \varphi(\spe{h}/C_1) = 1$, because $\varphi_{C'}(\spe{h}) = 1$. Since $|N \setminus C_1| < |N|$, it follows from induction that $\varphi_{C/C_1}(\spe{h}/C_1) = 1$. Then $\varphi_{C}(\spe{h}) = 1$. Thus, $\varphi$ is a convex character. \end{proof}
Recall that $\zeta$ is the character defined in Definition \ref{def:zeta} and $\chi$ is the character defined in Definition \ref{def:chromatic}.
\begin{theorem}
Let $\mathbb{K}\spe{H}$ be a linearized Hopf monoid. Then $\zeta$ and $\chi$ are convex characters. Suppose that, for every finite set $N$, and every $\spe{h} \in \spe{H}_N$, there exists non-empty $S \subsetneq N$ such that $\Delta_{S, N \setminus S}(\spe{h}) \neq 0.$ Then $\zeta$ and $\chi$ are balanced convex characters.
\label{thm:chiconvex}
\end{theorem}
\begin{proof}
First, we show that $\zeta$ and $\chi$ are convex characters. Let $N$ be a finite set, and let $\spe{h} \in \spe{H}_N$, such that $\varphi(\spe{h}) = 1$, and let $C \leq C' \models N$ such that $\zeta_{C'}(\spe{h}) \neq 0$. Then $\Delta_{C'}(\spe{h}) \neq 0$. This implies that $\Delta_C(\spe{h}) \neq 0$. By definition of $\zeta$, it follows that $\zeta_C(\spe{h}) = 1$. Thus, $\zeta$ is a convex character.

Let $\spe{h} \in \spe{H}_N$ such that $\chi(\spe{h}) = 1$. Let $C \leq C' \models N$ such that $\chi_{C'}(\spe{h}) = 1$. Then there exists a linear order $\ell$ on $N$, and elements $\spe{h}_i \in \spe{H}_{\ell_i}$ such that $\spe{h} = \spe{h}_1 \cdots \spe{h}_{|N|}$. Let $F_{\bdot} = F(C)$. We know $F_i \setminus F_{i-1} = \{\ell_{i_1}, \ldots, \ell_{i_k} \}$ for some $i_1, \ldots, i_k$. Then by compatibility of products and coproducts we have \[\spe{h}|_{F_i} / F_{i-1} = \spe{h}_{1}|_{F_i}/F_{i-1} \cdot \spe{h}_{2}|_{F_i}/F_{i-1} \cdots \spe{h}_{|N|}|_{F_i}/F_{i-1}.\] We see that \[\spe{h}_j|_{F_i}/F_{i-1} = \begin{cases} 1 & \ell_j \not\in F_i \setminus F_{i-1} \\ \spe{h}_i & \ell_j \in F_i \setminus F_{i-1} \end{cases}\] for a given $i$ and $j$. 
Thus $\spe{h}|_{F_i}/F_{i-1} = \spe{h}_{\ell_{i_1}} \spe{h}_{\ell_{i_2}} \cdots \spe{h}_{\ell_{i_k}}$. Hence $\chi(\spe{h}|_{F_i}/F_{i-1}) = 1$, and $\chi_C(\spe{h}) = 1$. Therefore $\chi$ is a convex character.

Let $N$ be a finite set, and $\spe{h} \in \spe{H}_N$. If $|N| = 1$, then we see that $\zeta(\spe{h}) = \chi(\spe{h}) = 1$. Thus $\zeta$ and $\chi$ both satisfy the first condition to be a balanced convex character. We see that the second condition is the hypothesis of our theorem.

So suppose $|N| > 1$, and let $S \subsetneq N$ be such that $\Delta_{S, N \setminus S}(\spe{h}) \neq 0$. 
We see that $\zeta(\spe{h}) = \zeta(\spe{h}|_S) = \zeta(\spe{h}/S) = 1$, and thus $\zeta$ is a balanced convex character.

Finally, suppose that $\chi(\spe{h}) = 1$.
Then there exists a linear order $\ell$ on $N$, and objects $\spe{h}_i \in \spe{H}_{\{\ell_i \}}$ such that $\spe{h} = \spe{h}_1 \cdot \spe{h}_2 \cdots \spe{h}_{|N|}$. By the same argument as above, $\spe{h}|_S = \prod_{i \in [|N|]: \ell_i \in S} \spe{h}_i$ and $\spe{h}/S = \prod_{i \in [|N|]: \ell_i \in N \setminus S} \spe{h}_i$. Hence $\chi(\spe{h}|_S) = \chi(\spe{h}/S) = 1$. Thus $\chi$ is a balanced convex character. \end{proof}

\begin{example}
For the linearized Hopf monoids of graphs $\mathbb{K}\spe{G}^{\bdot}$, posets $\mathbb{K}\spe{P}^{\bdot}$, and matroids $\mathbb{K}\spe{M}^{\bdot}$, the conditions of Theorem \ref{thm:chiconvex} hold. Thus $\chi$ and $\zeta$ are balanced convex characters on these linearized Hopf monoids.
\end{example}

\subsection{The geometric realization}
\label{subsec:geometric}

Now we give the geometric realization. 
\begin{definition}
Let $(\mathbb{K}\spe{H},\mathbb{K}(\varphi))$ be a linearized combinatorial Hopf monoid. Let $N$ be a finite set, and let $V$ be the collection of all proper subsets of $N$, and define $\kappa: V \to [|N|-1]$ by $\kappa(S) = |S|$. Given $\spe{h} \in \spe{H}_N$, define the $\varphi$-coloring complex of $\spe{h}$ by
\begin{align*} & \Sigma_{\varphi}(\spe{h})  = F(\mathcal{A}_{\varphi}(\spe{h})) \\ &  = \left\{ S_1 \subset S_2 \subset \cdots \subset S_k \subset N: \left(\bigotimes_{i=1}^k \varphi_{S_i \setminus S_{i-1}} \right) \circ \Delta_{S_1 | S_2 \setminus S_1 | \cdots | N \setminus S_k} (\spe{h}) = 1 \right\}.\end{align*}
\label{def:coloringcomplex}
\end{definition}
Given $S_1 \subset S_2 \subset \cdots \subset S_k \subset N$ in $\Sigma_{\varphi}(\spe{h})$, we have $\varphi(\spe{h}|_{S_i}/S_{i-1}) = 1$ for $1 \leq i \leq k+1$, where we define $S_0 = \emptyset$ and $S_{k+1} = N$. In our proofs, we often use this fact, and say $\varphi(\spe{h}|_{S_i}/S_{i-1}) = 1$ for all $i$, with the understanding that we are including the exceptional cases $\varphi(\spe{h}/S_k) = \varphi(\spe{h}|_{S_1}) = 1$.

\begin{lemma}
Let $(\mathbb{K}\spe{H},\mathbb{K}(\varphi))$ be a linearized combinatorial Hopf monoid. Let $N$ be a finite set, and define $\kappa: 2^N \setminus \{\emptyset, N \} \to [|N|-1]$ by $\kappa(S) = |S|$. Given $\spe{h} \in \spe{H}_N$, let $\Sigma_{\varphi}(\spe{h})$ be the $\varphi$-coloring complex. If $\varphi$ is a convex character, then $\Sigma_{\varphi}(\spe{h})$ is a colored relative simplicial complex.
\label{lem:defsigma}
\end{lemma}

\begin{proof}
Since $\varphi$ is a convex character, it follows from Proposition \ref{prop:convexequiv} that $\mathcal{A}_{\varphi}(\spe{h})$ is a convex album for all $\spe{h}$. By Proposition \ref{prop:convex}  $F(\mathcal{A}_{\varphi}(\spe{h}))$ is a relative simplicial complex. Given a flag $F_{\bdot} : F_1 \subset F_2 \subset \cdots \subset F_k \in \Sigma_{\varphi}(\spe{h})$, we see that $\kappa(F_{\bdot}) = \{|F_1|, |F_2|, \ldots, |F_k| \}$. Hence $\Sigma_{\varphi}(\spe{h})$ is a colored relative simplicial complex. \end{proof}

For a finite set $N$, let $\spe{K}_N$ denote the set of all colored relative simplicial complexes $(\Phi, \kappa)$ such that $\Phi$ is a subcomplex of the Coxeter complex of type $A$ and $\kappa$ is the geometric coloring where $\kappa(S) = |S|$. Then $\spe{K}$ is a set species. Given a linearized combinatorial Hopf monoid $(\mathbb{K}\spe{H}, \mathbb{K}(\varphi))$, where $\varphi$ is a convex character, and $\spe{h} \in \spe{H}_N$, we see that $(\Sigma_{\varphi}(\spe{h}), \kappa) \in \spe{K}_N$. Thus we have a natural transformation $\Sigma_{\varphi}: \spe{H} \to \spe{K}$.
\begin{theorem}
\label{thm:geometricrealization}
Let $(\mathbb{K}\spe{H},\mathbb{K}(\varphi))$ be a linearized combinatorial Hopf monoid. If $\varphi$ is a convex character, then $(\spe{K}, \Sigma_{\varphi})$ is a geometric realization of $\spe{H}$.

Thus, for every $\spe{H}$-structure $\spe{h}$, we have 
\[\Psi_{\varphi}(\spe{h}, \mathbf{x}) = \dHilb(\Sigma_{\varphi}(\spe{h}), \mathbf{x}). \]
\end{theorem}
\begin{proof}

Let $\spe{h} \in \spe{H}_N$ for some finite set $N$.
We know from Lemma \ref{lem:defsigma} that $\Sigma_{\varphi}(\spe{h}) \in \spe{K}_N$. The reader can check that $\Sigma_{\varphi}: \spe{H} \to \spe{K}$ is a natural transformation.

We show that $\Psi_{\varphi}(\spe{h}, \mathbf{x}) = \dHilb(\Sigma_{\varphi}(\spe{h}), \mathbf{x})$ for every $\spe{h}$-structure. We see that \[\Psi_{\varphi}(\spe{h}, \mathbf{x}) = \sum_{C \in \mathcal{A}_{\varphi}(\spe{h})} M_{\alpha(C)} = \sum_{F_{\bdot} \in F(\mathcal{A}_{\varphi}(\spe{h}))} M_{\alpha(F_{\bdot})}.\] where the first equality comes from Theorem \ref{thm:coloring}, and the second equality is due to the correspondence between flags and set compositions.  Since $\alpha(F_{\bdot}) = \alpha(\kappa(F_{\bdot}))$, and $\Sigma_{\varphi}(\spe{h}) = F(\mathcal{A}_{\varphi}(\spe{h}))$, we have \[\Psi_{\varphi}(\spe{h}, \mathbf{x}) = \sum_{F_{\bdot} \in \Sigma_{\varphi}(\spe{h})} M_{\alpha(\kappa(F_{\bdot}))} = \dHilb(\Sigma_{\varphi}(\spe{h}), \mathbf{x}).\] \end{proof}
\begin{example}
Let $\mathbb{K}\spe{G}^{\bdot}$ be the combinatorial Hopf monoid of graphs. Then the character $\chi$ is balanced convex. We see that $\mathcal{A}_{\chi}(\spe{g})$ is the collection of set compositions where each part is an independent set. Then the corresponding relative simplicial complex $\Sigma_{\chi}(\spe{g})$ consists of chains $S_1 \subset S_2 \subset \cdots \subset S_k \subset N$ such that $S_i \setminus S_{i-1}$ is an independent set for all $i$. One presentation for $\Phi(\spe{g})$ is as the pair $(\Sigma, \Gamma(\spe{g}))$, where $\Sigma$ is the Coxeter complex of type $A$, and $\Gamma(\spe{g})$ is collection of flags $S_1 \subset S_2 \subset \cdots \subset S_k \subset N$ where $S_{i+1} \setminus S_i$ must contain an edge for some $i$. The subcomplex $\Gamma(\spe{g})$ is the coloring complex, as introduced by Steingr\'imsson \cite{steingrimsson}.
\end{example}

\begin{example}
Let $\mathbb{K}\spe{P}^{\bdot}$ be the combinatorial Hopf monoid of posets. Then the character $\zeta$ is balanced convex. We see that $\mathcal{A}_{\zeta}(\spe{p})$ is the collection of set compositions $C_1|C_2|\ldots|C_k$ where $C_1 \cup C_2 \cup \cdots \cup C_i$ is an order ideal for all $i$. If we let $J(\spe{p})$ be the set of order ideals of $\spe{p}$ ordered by inclusion, and let $\overline{J(\spe{p})} = J(\spe{p}) \setminus \{\emptyset, \spe{p} \}$, then $\Sigma_{\zeta}(\spe{p}) = \Delta(\overline{J(\spe{p})})$, the order complex of $\overline{J(\spe{p})}$. Thus $\Psi_{\zeta}(\spe{p}, \mathbf{x}) = \Hilb(\Delta(\overline{J(\spe{p})}, \mathbf{x}).$ Since $J(\spe{p})$ is a graded poset with a maximum and minimum element, $\Hilb(\Delta(\overline{J(\spe{p})}), \mathbf{x})$ is the $F$-quasisymmetric function of the graded poset $J(\spe{p})$. As mentioned in Example \ref{ex:ehrenborg}, the $F$-quasisymmetric function of a graded poset is a generalization of the $p$-partition enumerator.
\label{ex:complexposet}
\end{example}

\begin{example}
Let $\mathbb{K}\spe{M}^{\bdot}$ be the combinatorial Hopf monoid of matroids. Then the character $\chi$ is balanced convex. Given a matroid $\spe{m}$, the corresponding relative simplicial complex $\Sigma_{\chi}(\spe{m})$ consists of chains $S_1 \subset S_2 \subset \cdots \subset S_k \subset N$ such that $\spe{m}|_{S_i} /  S_{i-1}$ has a unique basis for all $i$.
\end{example}

\subsection{Balanced realizations}

We also discuss a condition which ensures all the geometric realizations are \emph{balanced} complexes.
\begin{theorem}

Let $\mathbb{K}\spe{H}$ be a linearized combinatorial Hopf monoid with convex character $\mathbb{K}(\varphi)$. Then $\Sigma_{\varphi}(\spe{h})$ is balanced of dimension $|N|-2$ for every non-zero $\spe{H}$-structure $\spe{h} \in \spe{H}_N$ with $|N| \geq 1$ if and only if $\mathbb{K}(\varphi)$ is a balanced convex character.
\label{thm:geometricbalanced}
\end{theorem}
\begin{proof}
Suppose that $\Sigma_{\varphi}(\spe{h})$ is balanced for every $\spe{H}$-structure $\spe{h}$. Let $|N| = 1$, and let $\spe{h} \in \spe{H}_N$. Since $\Sigma_{\varphi}(\spe{h})$ is balanced of dimension $|N|-2 = -1$, we have $\Sigma_{\varphi}(\spe{h}) = \{\emptyset\}$, which is a simplicial complex. Thus, $\varphi(\spe{h}) = 1$, satisfying the first condition for being a balanced convex character. 

Now let $|N| \geq 2$, and let $\spe{h} \in \spe{H}_N$. Since $\Sigma_{\varphi}(\spe{h})$ is balanced of dimension $|N|-2$, there exists a facet $F_1 \subset F_2 \subset \cdots \subset F_{|N|-1}$ in $\Sigma_{\varphi}(\spe{h})$. This implies that $\Delta_{F_1, N \setminus F_1}(\spe{h}) \neq 0$. Thus we have the second condition for being a balanced convex character.

Now suppose that $\varphi(\spe{h}) = 1$, and let $S \subsetneq N$ such that $\Delta_{S, N \setminus S}(\spe{h}) \neq 0$. We know that there are faces $E_1 \subset \cdots \subset E_{\ell} \in \Sigma_{\varphi}(\spe{h}|_S)$ and $F_1 \subset \cdots \subset F_{k} \in \Sigma_{\varphi}(\spe{h}/S)$. Then $\varphi((\spe{h}|_S)|_{E_i} / E_{i-1}) = 1$ for all $i$, and $\varphi((\spe{h}/S)|_{F_i} / F_{i-1}) = 1$ for all $i$. However, this is equivalent to stating that $\varphi(\spe{h}|_{E_i}/E_{i-1}) = 1$ and $\varphi(\spe{h}|_{F_i \cup S} / (F_{i-1} \cup S)) = 1$ for all $i$.
Thus \[G_{\bdot} := E_1 \subset \cdots \subset E_{\ell} \subset S \subset S \cup F_1 \subset \cdots \subset S \cup F_{k} \] is an element of $\Sigma_{\varphi}(\spe{h})$. Since $\varphi(\spe{h}) = 1$, we know $\emptyset \in \Sigma_{\varphi}(\spe{h})$. Since $\emptyset \subseteq \{S\} \subseteq G_{\bdot}$, and $\Sigma_{\varphi}(\spe{h})$ is a relative simplicial complex, it follows that $ \{S\} \in \Sigma_{\varphi}(\spe{h})$. By definition of $\Sigma_{\varphi}(\spe{h})$, we conclude that $\varphi_{S| N \setminus S}(\spe{h}) = 1$. Thus $\varphi$ is a balanced convex character.

Now suppose that $\varphi$ is a balanced convex character. Then $\varphi$ is convex, so $\Sigma_{\varphi}(\spe{h})$ is a colored relative simplicial complex with geometric coloring given by $\kappa(S) = |S|$. Thus we only need to show that $\Sigma_{\varphi}(\spe{h})$ is pure of dimension $|N|-2$, when $N \neq \emptyset$.
 If $|N| = 1$, then $\varphi(\spe{h}) = 1$, and thus $\Sigma_{\varphi}(\spe{h}) = \{\emptyset \}$, which is pure of dimension $-1$. 
 
 Next we show that $\Sigma_{\varphi}(\spe{h})$ has dimension $|N|-2$ for all non-zero $\spe{h}$ and $|N| \geq 2$. The proof is by induction on $N$, with the base case $|N| = 1$ already established.
 Since $|N| \geq 2$, there exists $S \subsetneq N$ such that $\Delta_{S, N \setminus S}(\spe{h}) \neq 0$. By induction, we know that there are faces $E_1 \subset \cdots E_{|S|-1} \in \Sigma_{\varphi}(\spe{h}|_S)$ and $F_1 \subset \cdots \subset F_{|N \setminus S|-1} \in \Sigma_{\varphi}(\spe{h}/S)$. Note that $|E_i \setminus E_{i-1}| = 1$ for all $i$, and $|F_i \setminus F_{i-1}| = 1$ for all $i$.
Then \[ G_{\bdot} := E_1 \subset \cdots \subset E_{|S|-1} \subset S \subset S \cup F_1 \subset \cdots \subset S \cup F_{|N \setminus S|-1} \] is a flag such that $\Delta_{C(G_{\bdot})}(\spe{h}) \neq 0$. Since $|G_i \setminus G_{i-1}| = 1$ for all $i$, and $\varphi(\spe{h}') = 1$ whenever $\spe{h}' \in \spe{H}_{M}$ and $|M| = 1$, it follows that $\varphi_{C(G_{\bdot})}(\spe{h}) = 1$. Thus, $G_{\bdot} \in \Sigma_{\varphi}(\spe{h})$ is of maximum dimension. Hence $\Sigma_{\varphi}(\spe{h})$ has dimension $|N| - 2$.

To show that $\Sigma_{\varphi}(\spe{h})$ is pure, let $F_{\bdot} \in \Sigma_{\varphi}(\spe{h})$. Suppose that $F_{\bdot}$ does not have dimension $|N|-2$. Then there exists $i$ such that $|F_i \setminus F_{i-1}| \geq 2.$ Let $N' = F_i \setminus F_{i-1}$, and $\spe{h}' = \spe{h}|_{F_i}/F_{i-1}$. Since $|N'| \geq 2$, there exists $M \subseteq N'$ such that $\Delta_{M, N' \setminus M}(\spe{h}') \neq 0$. Since $F_{\bdot} \in \Sigma_{\varphi}(\spe{h})$, we see that $\varphi(\spe{h}') = \varphi(\spe{h}|_{F_i}/F_{i-1}) = 1$. Since $\varphi$ is a balanced convex character, it follows that $\varphi(\spe{h}'|_M) = \varphi(\spe{h'}/M) = 1$. Then $\varphi(\spe{h}|_{F_i \cup M} / F_{i-1}) = 1$ and $\varphi(\spe{h}|_{F_i} / (F_{i-1} \cup M)) = 1$. Thus 
\[ F_1 \subset \cdots \subset F_i \subset F_i \cup M \subset F_{i+1} \subset \cdots \subset F_i \] is a face of $\Sigma_{\varphi}(\spe{h})$ that contains $F_{\bdot}$. Therefore, a face is maximal only if it has dimension $|N|-2$, which means that $\Sigma_{\varphi}(\spe{h})$ is pure. \end{proof}

\subsection{Some consequences}
\label{subsec:consequences}

Now we derive inequalities for the coefficients of $\Psi_{\varphi}(\spe{h}, \mathbf{x})$ and $\chi_{\varphi}(\spe{h}, x)$ with respect to the basis of monomial quasisymmetric functions and binomial coefficients, respectively.

\begin{theorem}
Let $(\mathbb{K}\spe{H},  \mathbb{K}(\varphi))$ be a linearized combinatorial Hopf monoid. Suppose that $\varphi$ is a balanced convex character. Let $N$ be a finite set and  $\spe{h} \in \spe{H}_N$ be a non-zero $\spe{H}$-structure.
Then $\Psi_{\varphi}(\spe{h}, \mathbf{x})$ is $M$-increasing. 
Also, if we write $\chi_{\varphi}(\spe{h}, x) = \sum_{i=0}^n f_{i-2}(\spe{h}) \binom{x}{i}$, then $f_{-2}(\spe{h}) = 0$ and $(f_{-1}, \ldots, f_{|N|-1})$ is super flawless.
\label{thm:consequence}
\end{theorem}
\begin{proof}
Since $\varphi$ is convex, we know that $\Psi_{\varphi}(\spe{h}, \mathbf{x}) = \dHilb(\Sigma_{\varphi}(\spe{h}), \mathbf{x})$, where $\Sigma_{\varphi}(\spe{h})$ is the $\varphi$-proper coloring complex of $\spe{h}$. Since $\varphi$ is a balanced convex character, $\Sigma_{\varphi}(\spe{h})$ is balanced. We know by Proposition \ref{prop:increasing} that $\dHilb(\Sigma_{\varphi}(\spe{h}), \mathbf{x})$ is $M$-increasing for any balanced relative simplicial complex.

We write $\chi_{\varphi}(\spe{h}, n) = \sum_{i=0}^n f_{i-2}(\spe{h}) \binom{n}{i}$ for integers $f_{i-2}(\spe{h})$. We have 
\begin{align*}
\chi_{\varphi}(\spe{h}, n) & = \ps_n \Psi_{\varphi}(\spe{h}, \mathbf{x}) \\ & =  \sum_{S \subseteq [n]} f_S(\Sigma_{\varphi}(\spe{h})) \binom{n}{|S|+1} \\ & = \sum_{i=1}^{|N|} f_{i-2}(\Sigma_{\varphi}(\spe{h})) \binom{n}{i}. \end{align*} Thus, $f_{-2}(\spe{h}) = 0$, and $f_i(\spe{h}) = f_i(\Sigma_{\varphi}(\spe{h}))$ for $i \geq -1$. Since $\Sigma_{\varphi}(\spe{h})$ is pure of dimension $|N|-2$, the second result follows from Proposition \ref{prop:flawless}. \end{proof}
Consider a polynomial $p(x) = \sum_{i=0}^n f_i \binom{x}{i}$ of degree $n$. We say $p(x)$ is super flawless if $f_0 = 0$ and $(f_1, \ldots, f_n)$ is super flawless.

Since we know that $\chi$ and $\zeta$ are balanced convex characters for the linearized Hopf monoids of graphs $\mathbb{K}\spe{G}^{\bdot}$, posets $\mathbb{K}\spe{P}^{\bdot}$, and matroids $\mathbb{K}\spe{M}^{\bdot}$, we can apply Theorem \ref{thm:consequence} to these cases.
\begin{corollary}
Let $G$ be a graph on a finite set $N$. Then $X(G, \mathbf{x}) = \dHilb(\Phi(G), \mathbf{x})$, where $\Phi(G)$ is the relative coloring complex of $G$.

Thus $X(G, \mathbf{x})$ is $M$-increasing, and $\chi(G, x)$ is super flawless.
\label{cor:graphs}
\end{corollary}
\begin{corollary}
Let $\spe{p}$ be a poset on a finite set $N$. Let $J(P)$ be the poset of order ideals of $\spe{p}$, ordered by inclusion. Then $\Psi_{\zeta}(\spe{p}, \mathbf{x}) = \dHilb(\Delta(J(\spe{p})), \mathbf{x})$, where $\Delta(J(\spe{p}))$ is the order complex of $J(\spe{p})$.

Thus $\Psi_{\zeta}(\spe{p}, \mathbf{x})$ is $M$-increasing, and $\chi_{\zeta}(\spe{p}, x)$ is super flawless.
\label{cor:posets}
\end{corollary}

\begin{corollary}
Let $\spe{m}$ be a matroid on a finite set $N$. Then $\Psi_{\chi}(\spe{m}, \mathbf{x})$ is $M$-increasing, and $\chi_{\chi}(\spe{m},x)$ is super flawless.
\label{cor:matroids}
\end{corollary}
We will give a lot more examples of linearized Hopf monoids in Section \ref{sec:applications}, along with balanced convex characters for those Hopf monoids. Theorem \ref{thm:consequence} applies to each example.

\section{Coloring complexes}
\label{sec:coloringcomplex}

A natural question is to determine what types of relative simplicial complexes arise from our geometric construction. The aim of this section is to characterize such complexes completely. 

Now we define \emph{abstract} relative coloring complexes. Let $N$ be a finite set.
Given flags
 $E_{\bdot} :  E_1 \subset E_2 \subset \cdots \subset E_k \subseteq N$ and $F_{\bdot} : F_1 \subset F_2 \subset \cdots \subset F_m \subseteq N$, such that $E_a = F_c$ for some $a$ and $c$, then we define their \emph{exchange} with respect to $E_a$ to be the flag \[E_{\bdot} \cdot_a F_{\bdot} : E_1 \subset \cdots \subset E_a \subset F_{c+1} \subset F_{c+2} \subset \cdots \subset F_m. \]

\begin{definition} 
Let $\mathcal{F}$ be a collection of flags of subsets of $N$. We say that $\mathcal{F}$ satisfies the \emph{flag exchange condition} if whenever we have two flags $E_{\bdot}, F_{\bdot} \in \mathcal{F}$ such that $E_i = F_j$ for some $i$ and $j$, then $E_{\bdot} \cdot_i F_{\bdot} \in \mathcal{F}$. 
\label{def:flagexchange}
\end{definition}

\begin{figure}
\begin{center}
\begin{tabular}{cc}
\begin{tikzpicture}
\draw[color=white, fill=gray!20] (30:2cm) -- (330:2cm) -- (210:2cm) -- (150:2cm) -- cycle;
  \node[circle, draw=red, fill=white, dashed, thick] (b) at (150:2cm) {$abc$};
  \node[circle, draw=red, fill=white, dashed, thick] (a) at (210:2cm) {$c$};
  \node[circle, draw=red, fill=white,dashed, thick] (e) at (330:2cm) {$acd$};
  \node[circle, draw=red, fill=white,dashed, thick] (d) at (30:2cm) {$a$};
  \node[circle, draw=black, fill=white,thick] (or) at (0:0cm) {$ac$};
  \draw[dashed, red, thick] (d) -- (e);
  \draw[dashed, red, thick] (a) -- (b);
  \draw[thick] (b) -- (d);
  \draw[thick] (a) -- (e);
  \draw[thick] (b) -- (or) -- (d);
  \draw[thick] (a) -- (or) -- (e);

\end{tikzpicture}
&

\begin{tikzpicture}
\draw[color=white, fill=gray!20] (30:2cm) -- (0:0cm) -- (150:2cm) -- cycle;
\draw[color=white, fill=gray!20] (330:2cm) -- (210:2cm) -- (0:0cm) -- cycle;
  \node[circle, draw=black, fill=white, thick] (b) at (150:2cm) {$abc$};
  \node[circle, draw=black, fill=white, thick] (a) at (210:2cm) {$c$};
  \node[circle, draw=black, fill=white, thick] (e) at (330:2cm) {$acd$};
  \node[circle, draw=black, fill=white, thick] (d) at (30:2cm) {$a$};
  \node[circle, draw=black, fill=white,thick] (or) at (0:0cm) {$ac$};

  \draw[thick] (b) -- (d);
  \draw[thick] (a) -- (e);
  \draw[thick] (b) -- (or) -- (d);
  \draw[thick] (a) -- (or) -- (e);

\end{tikzpicture}
\end{tabular}
\end{center}
\caption{Two examples of $F(\mathcal{A})$.}
\label{fig:coloringcomplexes}
\end{figure}

\begin{definition}
Recall that the Coxeter Complex $\Sigma_N$ is the set of all flags in $N$. It is a balanced simplicial complex, with coloring given by $\kappa(S) = |S|$ for all vertices $S$ in $\Sigma_N$. Then $\Sigma_N$ is a simplicial complex. A relative simplicial complex $\Phi$ is an \emph{abstract relative coloring complex} if $\Phi \subseteq \Sigma_N$, and where $\Phi$ satisfies the flag exchange condition. An abstract relative coloring complex is \emph{balanced} if it is balanced with respect to the coloring $\kappa$.
\label{def:abstractcoloring}
\end{definition}
The fact that we require $\Phi$ to be balanced with respect to the coloring $\kappa$ implies that $\Phi$ is pure of dimension $|N|-2$. This also implies that any flag $E_{\bdot} \in \Phi$ can be refined to a \emph{complete} flag $F_{\bdot}$. That is, there is a flag $F_{\bdot}$ with $E_{\bdot} \subseteq F_{\bdot}$, and such that $|F_i \setminus F_{i-1}| = 1$ for all $i$.

The relative simplicial complex on the left in Figure \ref{fig:coloringcomplexes} is a balanced abstract relative coloring complex. The relative simplicial complex $\Phi$ on the right in Figure \ref{fig:coloringcomplexes} is not. For example, $\{a\} \subset \{a,c\} \subset \{a,b,c\}$ and $\{c\} \subset \{a,c\} \subset \{a,c,d\}$ are two faces of $\Phi$ that share $\{a,c\}$, but $\{a\} \subset \{a,c\} \subset \{a,c,d\}$ is not a face of $\Phi$.

\begin{proposition}
Let $(\mathbb{K}\spe{H},  \mathbb{K}(\varphi))$ be a linearized combinatorial Hopf monoid, and suppose that $\varphi$ is a convex character. Let $\spe{h}$ be an $\spe{H}$-structure. Then $\Sigma_{\varphi}(\spe{h})$ is an abstract relative coloring complex.
\end{proposition}
\begin{proof}
Let $E_{\bdot}, F_{\bdot} \in \Sigma_{\varphi}(\spe{h})$ such that there exists $i$ with $E_i = F_k$ for some $k$. Then $\varphi(\spe{h}|_{E_{j}}/E_{j-1}) = 1$ for all $j \leq i$. Similarly, $\varphi(\spe{h}|_{F_{j}}/F_{j-1}) = 1$ for all $j \geq k$. This implies that $E_{\bdot} \cdot_i F_{\bdot} \in \Sigma_{\varphi}(\spe{h})$. \end{proof}
Thus, every $\varphi$-coloring complex is an abstract relative coloring complex. We show that \emph{every} abstract relative coloring complex is a $\varphi$-coloring complex for some combinatorial Hopf monoid $\spe{H}$. We do this by showing that the species of abstract relative coloring complexes is itself a linearized combinatorial Hopf monoid.
\subsection{Constructions on coloring complexes}

Let $F_{\bdot}$ be a flag of subsets $F_1 \subset F_2 \subset \cdots \subset F_k \subset N$ of $N$. For a subset $S \subsetneq N$, we let $F_{\bdot} \cap S = F_1 \cap S \subseteq F_2 \cap S \subseteq \cdots \subseteq F_k \cap S \subseteq S$. We view this as a proper flag of $S$ by removing duplicates.

Given an abstract relative coloring complex $\Sigma$ on $N$, and a nonempty subset $S \subsetneq N$, we let \[\Phi|_S = \{ F_{\bdot} \cap S: F_{\bdot} \in \Phi, S \in F_{\bdot} \} \] 
and \[\Phi/S = \{ F_{\bdot} \cap (N \setminus S): F_{\bdot} \in \Phi, S \in F_{\bdot} \}. \]
These two relative complexes are the \emph{restriction} and \emph{contraction} with respect to $S$. We define $\Phi|_N = \Phi/\{\emptyset\} = \Phi$, and, when $\Phi \neq \emptyset$, we define $\Phi|_{\emptyset} = \Phi/N = 1_{\emptyset}$. If there is no $F_{\bdot} \in \Phi$ with $S \in F_{\bdot}$, then $\Phi|_S = \Phi/S = 0_S$. Thus, it is possible for the restriction or contraction of a non-empty complex to be empty.

\begin{example}
Let $\Phi$ be the relative simplicial complex from the middle of Figure \ref{fig:productcomplex}. We observe that $\Phi$ is an abstract relative coloring complex. We see that $\Phi|_{\{a,d \}} = \emptyset$, since $\{a,d\}$ is not a vertex of $\Phi$.
On the other hand, $A(\Phi|_{\{a,b,c \}}) = \{ abc, a|b|c, a|c|b, c|a|b, ab|c, a|bc, ac|b, c|ab \}$.
\end{example}

\begin{proposition}
Let $\Phi$ be a non-empty abstract relative coloring complex on a finite set $N$, with $|N| \geq 2$ and let $S \subsetneq N$ be nonempty. Suppose that there exists $H_{\bdot} \in \Phi$ with $S \in H_{\bdot}$. Then $\Phi|_S$ is a non-empty abstract relative coloring complex on $S$, and $\Phi/S$ is a non-empty abstract relative coloring complex on $N \setminus S$. Moreover, if $\Phi$ is balanced, then $\Phi|_S$ and $\Phi/S$ are also balanced.
\label{prop:coproduct}
\end{proposition}
\begin{proof}
Since there exists $H_{\bdot}$ such that $S \in H_{\bdot}$, we see that $\Phi|_S$ and $\Phi/S$ are non-empty.

Let $E_{\bdot}, F_{\bdot} \in \Phi|_S$ such that there exists $i$ with $E_i \in F_{\bdot}$. Then there exists $D_{\bdot} \in \Phi$ such that $S \in D_{\bdot}$ and $D_{\bdot} \cap S = E_{\bdot}$. Similarly, there exists $G_{\bdot} \in \Phi$ such that $S \in G_{\bdot}$ and $G_{\bdot} \cap S = F_{\bdot}$. We see that $E_i \in D_{\bdot}$ and $E_i \in G_{\bdot}$. Since $\Phi$ is an abstract relative coloring complex, we have $D_{\bdot} \cdot_i G_{\bdot} \in \Phi$. Since $S \in D_{\bdot} \cdot_i G_{\bdot}$, we see that $E_{\bdot} \cdot_i F_{\bdot} = (D_{\bdot} \cdot_i G_{\bdot}) \cap S \in \Phi|_S.$ Thus $\Phi|_S$ is an abstract relative coloring complex.

Let $E_{\bdot}, F_{\bdot} \in \Phi/S$ such that there exists $i$ with $E_i \in F_{\bdot}$. Then there exists $D_{\bdot} \in \Phi$ such that $S \in D_{\bdot}$ and $D_{\bdot} \cap (N \setminus S) = E_{\bdot}$. Similarly, there exists $G_{\bdot} \in \Phi$ such that $S \in G_{\bdot}$ and $G_{\bdot} \cap (N \setminus S) = F_{\bdot}$. We see that $E_i \cup S \in D_{\bdot}$ and $E_i \cup S \in G_{\bdot}$. Choose $j$ such that $D_j = E_i \cup S $. Since $\Phi$ is an abstract relative coloring complex, we have $D_{\bdot} \cdot_j G_{\bdot} \in \Phi$. Since $S \in D_{\bdot} \cdot_j G_{\bdot}$, we see that $E_{\bdot} \cdot_i F_{\bdot} = (D_{\bdot} \cdot_j G_{\bdot}) \cap (N \setminus S) \in \Phi/S.$ Thus $\Phi/S$ is an abstract relative coloring complex.

Suppose that $\Phi$ is balanced. Let $E_{\bdot} \in \Phi|_S$. Then there exists a flag $F_{\bdot} \in \Phi$ such that $S \in F_{\bdot}$ and $ F_{\bdot} \cap S = E_{\bdot}$. Since $\Phi$ is balanced, there exists a refinement $G_{\bdot}$ such that $F_{\bdot} \subset G_{\bdot}$ and $G_i \setminus G_{i-1}$ is a singleton for all $i$. Then we see that $E_{\bdot} \subset G_{\bdot} \cap S \in \Phi|_S$. Hence $\Phi|_S$ is balanced. A similar argument shows that $\Phi/S$ is also balanced. \end{proof}
Given an abstract relative coloring complex $\Sigma$ on a finite set $M$, and an abstract relative coloring complex $\Gamma$ on a finite set $N$, with $M \cap N = \emptyset$, we define the product complex $\Sigma \cdot \Gamma$ to be the relative simplicial complex on $M \cup N$ given by \[\Sigma \cdot \Gamma = \{F_{\bdot}: F_{\bdot} \cap M \in \Sigma \mbox{ and } F_{\bdot} \cap N \in \Gamma \}. \] 
We see that if $\Sigma = 0_M$ or $\Gamma = 0_N$, then $\Sigma \cdot \Gamma = 0_{M \sqcup N}$. Note that the sets $M$ and $N$ are relevant for the definition of the product: if we pick $\Sigma = \{\emptyset \} \in \spe{C}_{\{a\}}$, and $\Gamma = \{\emptyset\} \in \spe{C}_{\{b\}}$, then $\Sigma \cdot \Gamma = \{ \emptyset, \{\{a\} \}, \{\{b\} \} \}$.

\begin{example}
Let $\Phi$ be given by $A(\Phi) = \{a|b, ab \}$, and $\Sigma$ be given by $A(\Sigma) = \{c|d\}$. Then $\Phi \cdot \Sigma$ appears in Figure \ref{fig:productcomplex}. We write $\Phi = (\Gamma, \Sigma)$, where $\Sigma$ consists of all the vertices, edges, and triangles in the figure, while $\Gamma$ consists of the dashed vertices and dashed edges.
\end{example}
\begin{figure}
\begin{center}
\begin{tikzpicture}
\draw[color=white, fill=gray!20] (90:2cm) -- (30:2cm) -- (330:2cm) -- (270:2cm) -- (210:2cm) -- (150:2cm) -- cycle;
\node[circle, draw=red, fill=white, dashed, thick] (c) at (90:2cm) {$ab$};
  \node[circle, draw=black, fill=white, thick] (b) at (150:2cm) {$abc$};
  \node[circle, draw=black, fill=white, thick] (a) at (210:2cm) {$c$};
  \node[circle, draw=red, fill=white,dashed, thick] (f) at (270:2cm) {$cd$};
  \node[circle, draw=red, fill=white,dashed, thick] (e) at (330:2cm) {$acd$};
  \node[circle, draw=red, fill=white,dashed, thick] (d) at (30:2cm) {$a$};
  \node[circle, draw=black, fill=white,thick] (or) at (0:0cm) {$ac$};
  \draw[dashed, red, thick] (c) -- (d) -- (e) -- (f);
  \draw[thick] (f) -- (a) -- (b) -- (c);
  \draw[thick] (b) -- (d);
  \draw[thick] (a) -- (e);
  \draw[thick] (b) -- (or) -- (d);
  \draw[thick] (a) -- (or) -- (e);

\end{tikzpicture}
\end{center}
\caption{A product of two coloring complexes.}
\label{fig:productcomplex}
\end{figure}

\begin{proposition}
Given an abstract relative coloring complex $\Sigma$ on a finite set $M$, and an abstract relative coloring complex $\Gamma$ on a finite set $N$, with $M \cap N = \emptyset$, then $\Sigma \cdot \Gamma$ is an abstract relative coloring complex on $N \sqcup M$. Moreover, if $\Sigma$ and $\Gamma$ are balanced, then so is $\Sigma \cdot \Gamma$.
\label{prop:product}
\end{proposition}
\begin{proof}
 Let $E_{\bdot}, F_{\bdot} \in \Sigma \cdot \Gamma$ such that there exists $i$ with $E_i \in F_{\bdot}$. We let $C_{\bdot} = E_{\bdot} \cap M$ and $D_{\bdot} = E_{\bdot} \cap N$. Similarly, we define $G_{\bdot} = F_{\bdot} \cap M$ and $H_{\bdot} = F_{\bdot} \cap N$. 
 We see that there exists an integer $j$ such that $C_j = E_i \cap M$. Then $C_j \in G_{\bdot}$, and $C_{\bdot} \cdot_j G_{\bdot} \in \Sigma$. There also exists an integer $k$ such that $D_k = E_i \cap N$. Then $D_k \in H_{\bdot}$, and $D_{\bdot} \cdot_k H_{\bdot} \in \Gamma$. Finally, we observe that $(E_{\bdot} \cdot_i F_{\bdot}) \cap M = C_{\bdot} \cdot_j G_{\bdot}$ and $(E_{\bdot} \cdot_i F_{\bdot}) \cap N = D_{\bdot} \cdot_k H_{\bdot}$. Thus $E_{\bdot} \cdot_i F_{\bdot} \in \Sigma \cdot \Gamma$. Therefore $\Sigma \cdot \Gamma$ is an abstract relative coloring complex.
 
 Now suppose that $\Sigma$ and $\Gamma$ are balanced, and let $F_{\bdot} \in \Sigma \cdot \Gamma$. Write $F_{\bdot} :  F_1 \subset F_2 \subset \cdots \subset F_k$. Suppose that $F_{\bdot}$ is maximal. If $|F_i \setminus F_{i-1}| = 1$ for all $i$, then $F_{\bdot}$ has dimension $(\dim \Sigma + \dim \Gamma)$ as required. So suppose $|F_i \setminus F_{i-1}| \geq 2$ for some $i$. If $F_{i-1} \setminus F_i$ contains elements from both $M$ and $N$, then it is possible to refine $F_{\bdot}$, contradicting the fact that it is maximal. Hence, without loss of generality, we may assume that $F_i \setminus F_{i-1} \subset M$. Let $E_{\bdot} = F_{\bdot} \cap M$. Choose $j$ such that $E_j = F_i \cap M$. Then $|E_j \setminus E_{j-1}| \geq 2$. Since $\Sigma$ is balanced, $E_{\bdot}$ is not a facet, and it is possible to refine it further, by adding a new set $E$ between $E_{j-1}$ and $E_j$. However, then we see that we can add $E \cup F_{i-1}$ between $F_{i-1}$ and $F_i$, contradicting the assumption that $F_{\bdot}$ was maximal. Therefore $\Sigma \cdot \Gamma$ must be balanced.  \end{proof}

\subsection{Hopf monoid structure}
For a finite set $N$, let $\spe{C}_N$ be the set of coloring complexes on $N$. We view $\spe{C}_N$ as a pointed set by making the void complex the base point. Thus $\spe{C}$ is a pointed set species.

Now we define the multiplication and comultiplication for $\mathbb{K}\spe{C}$.
Given finite sets $M$ and $N$, and coloring complexes $\Sigma \in \spe{C}_M$ and $\Gamma \in \spe{C}_N$, we see from Proposition \ref{prop:product} that $\Sigma \cdot \Gamma \in \spe{C}_{M \sqcup N}$. We define $\mu_{M,N}$ by $\mu_{M,N}( \Sigma \otimes \Gamma) = \Sigma \cdot \Gamma$. 

Given finite sets $M$ and $N$, and an abstract relative coloring complex $\Phi \in \spe{C}_{M \sqcup N}$, let \[\Delta_{M,N}(\Phi) = \Phi|_M \otimes \Phi/M. \]
Note that $\Delta_{M,N}(\Phi) = 0_M \otimes 0_N$ if there is no $F_{\bdot} \in \Phi$ with $M \in F_{\bdot}$.

For a finite set $N$, with $|N| \geq 1$, we let $\spe{BC}_N$ be the set of all balanced coloring complexes. We let $\spe{BC}_{\emptyset}$ consist of the empty complex $1_{\emptyset}$. Then $\spe{BC}^{\bdot}$ is a subspecies of $\spe{C}$, where we are identifying the $0 \in \spe{BC}^{\bdot}_N$ with the void complex of $\spe{C}_N$.
\begin{theorem}
The species $\mathbb{K}\spe{C}$ of coloring complexes is a linearized Hopf monoid, and the species $\mathbb{K}\spe{BC}^{\bdot}$ of balanced coloring complexes is a Hopf submonoid.
\end{theorem}
\begin{proof}

 Given finite sets $L, M$ and $N$, and abstract relative coloring complexes $\Pi \in \spe{C}_L$, $\Sigma \in \spe{C}_M$, and $\Gamma \in \spe{C}_N$, we see that $(\Pi \cdot \Sigma) \cdot \Gamma = \Pi \cdot (\Sigma \cdot \Gamma) = \{ E_{\bdot}: E_{\bdot} \cap L \in \Pi, E_{\bdot} \cap M \in \Sigma, E_{\bdot} \cap N \in \Gamma \}$. Extending by linearity, we see that multiplication is associative.
 
 Note here that $\spe{C}_{\emptyset}$ has one non-zero element, the empty complex $1_{\emptyset}$. We observe that $1_{\emptyset} \cdot \Sigma = \Sigma$ for all $\Sigma \in \spe{C}_N$. Hence $1_{\emptyset}$ in $\spe{C}_{\emptyset}$ is the multiplicative identity. 
 
 Let $L, M, N$ be finite sets, and let $\Phi \in \spe{C}_{L \sqcup M \sqcup N}$. Define $\Phi_1 = \{F_{\bdot} \cap L: F_{\bdot} \in \Phi, L, L \cup M \in F_{\bdot} \},$ $\Phi_2 = \{F_{\bdot} \cap M: F_{\bdot} \in \Phi, L, L \cup M \in F_{\bdot} \}$, and $\Phi_3 = \{F_{\bdot} \cap N: F_{\bdot} \in \Phi, L, L \cup M \in F_{\bdot} \}$. 
 Then we observe that 
 \begin{align*}  ((\Delta_{L,M} \otimes \id_N) \circ \Delta_{L \cup M, N})(\Phi) & = \Phi_1 \otimes \Phi_2 \otimes \Phi_3  \\
 & = ((\id_L \otimes \Delta_{M,N}) \circ \Delta_{L, M \cup N})(\Phi). \end{align*} Thus, the coproduct is coassociative. We observe that the counit axiom is also satisfied.
 
 Let $A, B, C$ and $D$ be finite sets, and let $\Sigma \in \spe{C}_{A \sqcup B}$ and $\Gamma \in \spe{C}_{C \sqcup D}.$ We show that $(\Sigma \cdot \Gamma)|_{A \sqcup C} = \Sigma|_A \cdot \Gamma|_C$. 
 Let $F_{\bdot} \in (\Sigma \cdot \Gamma)|_{A \sqcup C}$. Then there exists $E_{\bdot} \in \Sigma \cdot \Gamma$ such that $A \sqcup C \in E_{\bdot}$ and $E_{\bdot} \cap (A \sqcup C) = F_{\bdot}.$ Since $E_{\bdot} \in \Sigma \cdot \Gamma,$ we know that $E_{\bdot} \cap (A \sqcup B) \in \Sigma.$ Since $A \sqcup C \in E_{\bdot}$, we see that $A \in E_{\bdot} \cap (A \sqcup B).$ This, together with the fact that $E_{\bdot} \cap (A \sqcup B) \in \Sigma$ shows that $E_{\bdot} \cap A = (E_{\bdot} \cap (A \sqcup B)) \cap A \in \Sigma|_A.$ Hence $E_{\bdot} \cap A = E_{\bdot} \cap (A \sqcup C) \cap A = F_{\bdot} \cap A \in \Sigma|_A.$ By a similar argument, $F_{\bdot} \cap C \in \Gamma|_C$. Thus, we see that $F_{\bdot} \in \Sigma|_A \cdot \Gamma|_C.$
 
 Let $F_{\bdot} \in \Sigma|_A \cdot \Gamma|_C.$ Then $F_{\bdot} \cap A \in \Sigma|_A$. Thus, there exists $D_{\bdot} \in \Sigma$ with $A \in D_{\bdot}$ and $D_{\bdot} \cap A = F_{\bdot} \cap A.$ Similarly, there exists $E_{\bdot} \in \Gamma$ with $C \in E_{\bdot}$ and $E_{\bdot} \cap C = F_{\bdot} \cap C.$ We see that there exists an $i$ such that $D_i = A$, and a $j$ such that $E_j = C$. Let \[G_{\bdot} = F_{\bdot} \cup \{D_k \cup A \cup C: k > i \} \cup \{E_k \cup A \cup B \cup C: k > j \}.\] By construction, $G_{\bdot} \cap (A \cup B) = D_{\bdot}$, and $G_{\bdot} \cap (C \cup D) = E_{\bdot}$. Hence $G_{\bdot} \in \Sigma \cdot \Gamma$. Moreover, $G_{\bdot} \cap (A \cup C) = F_{\bdot}$ and $A \cup C \in G_{\bdot}$. Thus $F_{\bdot} \in (\Sigma \cdot \Gamma)|_{A \cup C}$. Hence we have shown that $(\Sigma \cdot \Gamma)|_{A \sqcup C} = (\Sigma|_A) \cdot (\Gamma|_C)$.
 
 A similar proof shows that $(\Sigma \cdot \Gamma)/(A \sqcup C) = (\Sigma / A) \cdot (\Gamma / C)$. Thus, the multiplication and comultiplication are compatible.

 For the final claim, it is enough to notice that the set of balanced abstract relative coloring complexes $\spe{BC}_N$ is a subset of $\spe{C}_N$, and that this inclusion induces an inclusion $\spe{BC} \subset \spe{C}$ of species. We have already shown that the product and coproduct of balanced abstract relative coloring complexes remains balanced, and hence the result follows. \end{proof}

\begin{remark}
If we let $\spe{K}_N$ be the collection of all subcomplexes of the Coxeter complex of type $A$, then $\mathbb{K}\spe{K}_N$ is closed under the same product and coproduct operations as $\spe{C}$. However, the coproduct is no longer coassociative. For example, if we let $\Sigma$ be the simplicial complex on the right in Figure \ref{fig:coloringcomplexes}, then $\Sigma|_{\{a,c\}} = \{\emptyset, \{a\}, \{c\} \}$, while $(\Sigma|_{\{a,b,c\}})|_{\{a,c\}} = \{\emptyset, \{a\}\}.$ Hence \[(\Delta_{\{a,c \}, \{b\}} \otimes \id_{\{d\}}) \circ \Delta_{\{a,b,c\}, \{d\}} \neq (\id_{\{a,c\}} \otimes \Delta_{\{b,d\}}) \circ \Delta_{\{a,c\}, \{b,d\}}.\]
\end{remark}

Let $N$ be a finite set, and let $\sigma: \spe{C}_N \to \spe{E}_N$ be given by 
\[\sigma(\Sigma) = \begin{cases}
1 & \text{if $\emptyset \in \Sigma$} \\
0 & \text{otherwise.}
\end{cases}
\]
The reader may check that this is a multiplicative function.

\begin{lemma}
Let $\Phi$ be an abstract relative coloring complex on $N$. Let $C \models N$. Then $\sigma_C(\Phi) = 1$ if and only if $F(C) \in \Phi$.
\label{lem:coproducttolink}
\end{lemma}
\begin{proof}
Let $C \models N.$ Let $F_{\bdot} = F(C)$. Suppose that $F_{\bdot} \in \Phi$. Recall that $\Phi|_{F_i}/F_{i-1} = \{E_{\bdot} \cap C_{i}: F_{i-1}, F_{i} \in E_{\bdot} \}$. We see that $F_{\bdot} \cap C_{i} = \{C_{i} \}$, which implies that $\emptyset \in \Phi|_{F_i}/F_{i-1}$.  Hence $\sigma_C(\Phi) = 1$.

Now suppose that $\sigma_C(\Phi) = 1$ with $|C| = k$. Then we see that $\emptyset \in \Phi|_{F_i}/F_{i-1}$ for all $i$. For each $i$, there exists $E_{\bdot}^i$ and an integer $j_i$ such that $E^i_{j_i} = F_{i-1}$ and $E^i_{j_i+1} = F_{i}$. Since $\Phi$ is closed under exchange operations, we can define $G^i_{\bdot}$ recursively by letting $G^0_{\bdot} = E^k_{\bdot}$, and $G^i_{\bdot} = E^{k-i}_{\bdot} \cdot_{j_{k-i}} G^{i-1}_{\bdot}$. One can check, using induction, that $G^i_{k-j} = F_{k-j}$ for $j \leq i$. Thus $G^k_{\bdot} = F_{\bdot}$, and $F_{\bdot} \in \Phi$. \end{proof}

\begin{proposition}
The morphism $\sigma: \mathbb{K}\spe{C} \to \mathbb{K}\spe{E}$ is a convex character, and the morphism $\sigma: \mathbb{K}\spe{BC}^{\bdot} \to \mathbb{K}\spe{E}$ is a balanced convex character.
\end{proposition}
\begin{proof}
Let $\Phi$ be an abstract relative coloring complex on $N$ and suppose $\sigma(\Phi) = 1$. Let $C \leq C' \models N$, and suppose that $\sigma_{C'}(\Phi) = 1$. Then we see that $\Phi$ is a simplicial complex, since $\emptyset \in \Phi$. Moreover, by Lemma \ref{lem:coproducttolink}, $F(C')$ is a face of $\Phi$. However, since $F(C) \subseteq F(C')$, it follows that $F(C) \in \Phi$. By Lemma \ref{lem:coproducttolink}, $\sigma_C(\Phi) = 1$. Therefore $\sigma$ is a convex character.

We see that $\sigma$ restricts to a convex character on $\spe{BC}$. Thus there is a geometric realization $\Sigma_{\sigma}: \spe{BC} \to \spe{K}$. Based upon the definition of $\Sigma_{\sigma}$, and Lemma \ref{lem:coproducttolink}, we see that $\Sigma_{\sigma}(\Phi) = \Phi$. In particular, the geometric realization of a balanced abstract relative coloring complex of dimension $|N|-2$ is itself, which is a balanced simplicial complex of dimension $|N|-2$. By Theorem \ref{thm:geometricbalanced}, it follows that $\sigma$ is a balanced convex character on $\spe{BC}$. \end{proof}

\begin{theorem}
The pair $(\mathbb{K}\spe{C}, \mathbb{K}(\sigma))$ is the terminal object in the category of combinatorial Hopf monoids with convex characters. In particular, for every Hopf monoid $\spe{H}$ with convex character $\varphi$, the map $\Sigma_{\varphi}: \spe{H} \to \spe{C}$ given by $\Sigma_{\varphi}(\spe{h}) = F(\mathcal{A}_{\varphi}(\spe{h}))$ is a morphism of Hopf monoids. 

Similarly, the pair $(\mathbb{K}\spe{BC}^{\bdot}, \mathbb{K}(\varphi))$ is the terminal object in the category of combinatorial Hopf monoids with balanced convex characters.
\label{thm:terminal}
\end{theorem}
\begin{proof}
Let $(\mathbb{K}\spe{H},  \mathbb{K}(\varphi))$ be a linearized combinatorial Hopf monoid such that $\varphi$ is a convex character. Given an $\spe{H}$-structure $\spe{h}$, we see that $\Sigma_{\varphi}(\spe{h})$ is an abstract relative coloring complex. Thus we obtain a map $\Sigma_{\varphi, N}: \spe{H}_N \to \spe{C}_N$. The reader can verify that this map is natural in $N$, and hence we obtain a natural transformation $\Sigma_{\varphi}: \spe{H} \to \spe{C}$. We show that the linearization $\mathbb{K}(\Sigma_{\varphi})$ is a morphism of Hopf monoids.

Let $\spe{h}$ and $\spe{h}'$ be $\spe{H}$-structures with $\spe{h} \in \spe{H}_M$ and $\spe{h}' \in \spe{H}_N$. Let $F_{\bdot} :  F_1 \subset F_2 \subset \cdots \subset F_k$. Then $F_{\bdot} \in \Sigma_{\varphi}(\spe{h} \cdot \spe{h}')$ if and only if $\varphi_{C(F_{\bdot})}(\spe{h} \cdot \spe{h}') = 1.$  Then $\varphi( (\spe{h} \cdot \spe{h}')|_{F_i}/F_{i-1}) = 1$ for all $i$ if and only if $\varphi(\spe{h}|_{F_{i} \cap M}/F_{i-1} \cap M) = \varphi(\spe{h}'|_{F_{i} \cap N}/F_{i-1} \cap N) = 1$ for all $i$. Thus $F_{\bdot} \cap M \in \Sigma_{\varphi}(\spe{h})$ and $F_{\bdot} \cap N \in \Sigma_{\varphi}(\spe{h}')$, which is equivalent to $F_{\bdot} \in \Sigma_{\varphi}(\spe{h}) \cdot \Sigma_{\varphi}(\spe{h}')$. Thus we have $\Sigma_{\varphi}(\spe{h} \cdot \spe{h}') = \Sigma_{\varphi}(\spe{h}) \cdot \Sigma_{\varphi}(\spe{h}').$ Hence the morphism preserves multiplication.

We focus on comultiplication.
Let $F_{\bdot} \in \Sigma_{\varphi}(\spe{h})|_S$. Then there exists $E_{\bdot} \in \Sigma_{\varphi}(\spe{h})$ such that $S \in E_{\bdot}$ and $E_{\bdot} \cap S = F_{\bdot}$. We see that $(\spe{h}|_S)|_{F_i}/F_{i-1} = \spe{h}|_{E_i}/E_{i-1}$ for $i \leq \ell(F_{\bdot})$. Then $\varphi(\spe{h}|_{E_i}/E_{i-1}) = 1$ for all $i$, which implies that $\varphi((\spe{h}|_S)|_{F_i}/F_{i-1}) = 1$ for all $i$. Thus $F_{\bdot} \in \Sigma_{\varphi}(\spe{h}|_S).$ Thus we conclude that $\Sigma_{\varphi}(\spe{h})|_S \subseteq \Sigma_{\varphi}(\spe{h}|_S).$
By a similar argument, $\Sigma_{\varphi}(\spe{h})/S \subseteq \Sigma_{\varphi}(\spe{h}/S).$

If $\Sigma_{\varphi}(\spe{h}/S) = 0_{N \setminus S}$ or $\Sigma_{\varphi}(\spe{h}|_S) = 0_S$, then we have $\Delta_{S, T}(\Sigma_{\varphi}(\spe{h})) = 0 = \Sigma_{\varphi}(\spe{h}|_S) \otimes \Sigma_{\varphi}(\spe{h}/S)$. So assume both relative complexes are not the void complex. Let $E_{\bdot} \in \Sigma_{\varphi}(\spe{h}|_S)$ and let $F_{\bdot} \in \Sigma_{\varphi}(\spe{h}/S).$ Define $G_{\bdot}$ by \[G_{\bdot} :  E_1 \subset E_2 \subset \cdots \subset E_k \subset S \subset S \cup F_1 \subset S \cup F_2 \subset \cdots \subset S \cup F_r.\] Since $E_{\bdot} \in \Sigma_{\varphi}(\spe{h}|_S)$, we have $1 = \varphi((\spe{h}|_S)|_{E_i}/E_{i-1}) = \varphi(\spe{h}|_{E_i}/E_{i-1})$ for all $i$. Since $F_{\bdot} \in \Sigma_{\varphi}(\spe{h}/S),$ we have $1 = \varphi((\spe{h}/S)|_{F_i}/F_{i-1}) = \varphi(\spe{h}|_{S \cup F_{i}} / S \cup F_{i-1})$ for all $i$. We see then that $\varphi(\spe{h}|_{G_i}/G_{i-1}) = 1$ for all $i$. Thus $G_{\bdot} \in \Sigma_{\varphi}(\spe{h})$ with $S \in G_{\bdot}.$ Then $E_{\bdot} = G_{\bdot} \cap S \in \Sigma_{\varphi}(\spe{h})|_S$, and $F_{\bdot} = G_{\bdot} \cap (N \setminus S) \in \Sigma_{\varphi}(\spe{h})/S$. Hence $\Sigma_{\varphi}(\spe{h}|_S) = \Sigma_{\varphi}(\spe{h})|_S$ and $\Sigma_{\varphi}(\spe{h}/S) = \Sigma_{\varphi}(\spe{h})/S$.
Therefore $\Delta_{S,T}(\Sigma_{\varphi}) = (\Sigma_{\varphi} \otimes \Sigma_{\varphi}) \circ \Delta_{S,T}$, and $\mathbb{K}(\Sigma_{\varphi})$ is a morphism of Hopf monoids.

Now we show that $\mathbb{K}(\Sigma_{\varphi}): \mathbb{K}\spe{H} \to \mathbb{K}\spe{C}$ is the unique Hopf monoid homomorphism with the property that $\sigma(\Sigma_{\varphi}(\spe{h})) = \varphi(\spe{h})$ for every $\spe{H}$-structure $\spe{h}$. Let $\Gamma$ be another such linearized Hopf monoid homomorphism. Let $\spe{h}$ be an $\spe{H}$-structure. Let $F_{\bdot}$ be a flag. Then by Lemma \ref{lem:coproducttolink}, $F_{\bdot} \in \Sigma_{\varphi}(\spe{h})$ if and only if $\varphi_{C(F_{\bdot})}(\spe{h}) = 1$. Since $\Gamma$ is a morphism of Hopf monoids which satisfies $\sigma \circ \Gamma = \varphi$, then $\sigma_{C(F_{\bdot})}(\Gamma(\spe{h})) = 1$. By Lemma \ref{lem:coproducttolink}, this holds if and only if $F \in \Gamma(\spe{h})$. Thus, if $\varphi_{C(F_{\bdot})}(\spe{h}) = 1,$ then $F_{\bdot}$ is a face of both $\Sigma_{\varphi}(\spe{h})$ and $\Gamma(\spe{h})$. On the other hand, if $\varphi_{C(F_{\bdot})}(\spe{h}) = 0,$ then $\sigma_{C(F_{\bdot})}(\Gamma(\spe{h})) = 0,$ and $F_{\bdot}$ is not a face of either of complex. Therefore $\Sigma_{\varphi}(\spe{h}) = \Gamma(\spe{h}).$

Finally, let $\mathbb{K}\spe{H}$ be a linearized combinatorial Hopf monoid with a balanced convex character $\mathbb{K}(\varphi).$ Then we have a morphism $\Sigma_{\varphi}: \spe{H} \to \spe{C}$. Moreover, by Theorem \ref{thm:geometricbalanced}, $\Sigma_{\varphi}(\spe{h})$ is balanced for all $\spe{h}$, which means that $\Sigma_{\varphi}(\spe{H}_N) \subseteq \spe{BC}_N$ for all $N$. Thus, we have a morphism $\Sigma_{\varphi}: \spe{H} \to \spe{BC}$, and $\mathbb{K}(\Sigma_{\varphi})$ is a morphism of Hopf monoids. \end{proof}

\section{Applications}
\label{sec:applications}
\subsection{Graphs}
\label{subsec:graphs}
We already discussed the standard character that yields the chromatic symmetric function. Here we detail one more interesting character that was first studied by Aval, Bergeron, and Machacek \cite{aval-bergeron-machacek}. Let 
\[\psi(\spe{g}) = 
\begin{cases} 
1 & \text{if $\spe{g}$ is a perfect matching} \\ 
0 & \text{otherwise.} 
\end{cases} \] 
It is clear that $\psi$ is a character. The resulting quasisymmetric function, $\Psi_{\psi}(\spe{g}, \mathbf{x})$, enumerates functions $f: V \to \mathbb{N}$ such that the induced subgraph on each color class is a perfect matching. 

We show that $\psi$ is a convex character. Let $\psi(\spe{g}) = 1$. Then $\spe{g}$ is a perfect matching. Let $C \leq C'$ and suppose that $\psi_{C'}(\spe{g}) = 1$. Hence, $\spe{g}$ restricted to each block of $C'$ is also a perfect matching. This implies that $\spe{g}|_{C_i}$ is also a perfect matching for each block of $C$. Hence $\psi_C(\spe{g}) = 1$. Therefore $\psi$ is a convex character.

We see that $\Sigma_{\psi}(\spe{g}) \neq \emptyset$ only when $\spe{g}$ has a perfect matching. Thus $\psi$ is not a balanced convex character. However, it does possess similar properties. If $\spe{g}$ has a perfect matching, then the facets of $\Sigma_{\psi}(\spe{g})$ are associated to set compositions where each block has size two, corresponding to taking a matching of $\spe{g}$, and putting the edges of the matching in a linear order. Hence, $\Sigma_{\psi}(\spe{g})$ is pure, and is balanced with the coloring given by $\kappa'(S) = |S|/2$. 

Given an integer composition $\alpha,$ we let $2 \alpha = (2\alpha_1, 2\alpha_2, \ldots, 2\alpha_k).$
Given $C \in \mathcal{A}_{\psi}(\spe{g})$ with $\alpha(C) = \alpha$, and $F_{\bdot} = F(C)$, we have $|F_i \setminus F_{i-1}| = \alpha_i$ for all $i$. However, $\kappa'(F_{\bdot}) = \{|F_1|/2, |F_2|/2, \ldots, |N \setminus F_k| / 2 \}$. Thus $2\alpha(\kappa'(F_{\bdot})) = \alpha(C)$. Since this implies that every block of $C$ has even size, we see that $[M_{\alpha}] \Psi_{\psi}(\spe{g}, \mathbf{x}) = 0$ if $\alpha$ contains an odd part. Otherwise, we see that, given $S \subseteq [|V|/2 - 1]$, we have $f_S(\Sigma_{\psi}(\spe{g})) = [M_{2\alpha(S)}] \Psi_{\psi}(\spe{g}, \mathbf{x}).$
 Since $\Sigma_{\psi}(\spe{g})$ is balanced with respect to $\kappa'$, we see that $f_S(\Sigma_{\psi}(\spe{g})) \leq f_T(\Sigma_{\psi}(\spe{g}))$ whenever $S \subseteq T \subseteq [|V|/2 - 1]$. Thus if $\alpha \leq \beta$ and $\beta$ contains only even parts, it follows that $[M_{\alpha}] \Psi_{\psi}(\spe{g}, \mathbf{x}) \leq [M_{\beta}] \Psi_{\psi}(\spe{g}, \mathbf{x})$.
 
Moreover, if we write $\Psi_{\psi}(\spe{g}, x) = \sum_{i=0}^n f_i \binom{x}{i}$, then $f_0 = 0$, and $f_j = 0$ for $j > n/2$. Also, the sequence $(f_1, \ldots, f_{n/2})$ is strongly flawless, and the polynomial has degree $|N|/2$. Aval, Bergeron, and Machacek \cite{aval-bergeron-machacek} showed that the polynomial can be computed from chromatic polynomials.

They also give a large class of invariants coming from linearized characters. Let $\mathcal{C}$ be a collection of isomorphism classes of connected graphs. Then we can define $\varphi_{\mathcal{C}}$ by
\[\varphi_{\mathcal{C}}(\spe{g}) = 
\begin{cases} 
1 & \text{if every component of $\spe{g}$ is isomorphic to a graph in $\mathcal{C}$} \\ 
0 & \text{otherwise.} 
\end{cases} \]  In particular, $\chi = \varphi_{\{K_1\}}$ and $\psi = \varphi_{\{K_2\}}$, where $K_n$ is the complete graph on $n$ vertices.

We are able to classify the balanced convex characters on graphs.
\begin{theorem}
Let $\varphi$ be a linearized character on $\mathbb{K}\spe{G}^{\bdot}$. Then there exists a collection $\mathcal{C}$ of isomorphism classes of connected graphs such that $\varphi = \varphi_{\mathcal{C}}$. Moreover, $\varphi_{\mathcal{C}}$ is a balanced convex character if and only if $\mathcal{C}$ satisfies the following two properties:
\begin{enumerate}
    \item For every $\spe{g} \in \mathcal{C}$, and every induced connected subgraph $\spe{h}$ of $\spe{g}$, we have $\spe{h} \in \mathcal{C}$.
    \item The collection $\mathcal{C}$ contains $K_1$.
\end{enumerate}
\end{theorem}
\begin{proof}

Let $\mathbb{K}(\varphi)$ be a linearized character on $\mathbb{K}\spe{G}^{\bdot}$. Let $\mathcal{C}$ be the isomorphism classes of connected graphs $\spe{g}$ for which $\varphi(\spe{g}) = 1$. We claim that $\varphi = \varphi_{\mathcal{C}}$. Let $\spe{g} \in \spe{G}_N$. Then there are connected graphs $\spe{g}_1, \ldots, \spe{g}_k$ such that $\spe{g} = \spe{g}_1 \cdots \spe{g}_k$. By definition, we have $\varphi_{\mathcal{C}}(\spe{g}_i) = \varphi(\spe{g}_i)$ for all $i$. Thus $\varphi(\spe{g}) = \varphi(\spe{g}_1) \cdots \varphi(\spe{g}_k) = \varphi_{\mathcal{C}}(\spe{g}_1) \cdots \varphi_{\mathcal{C}}(\spe{g}_k) = \varphi_{\mathcal{C}}(\spe{g}).$

Let $\varphi_{\mathcal{C}}$ be a balanced convex character. We denote $\varphi_{\mathcal{C}}$ by $\varphi$. Suppose that $K_1 \not\in \mathcal{C}$, and let $\spe{g} \in \mathcal{C}$ be a graph of minimum size. Then $\varphi(\spe{g}) = 1$. Since $\varphi$ is balanced, there must be a nontrivial composition $C$ such that $\varphi_C(\spe{g}) = 1$, which implies that there is a smaller graph in $\mathcal{C}$, a contradiction. Thus $K_1 \in \mathcal{C}.$

Let $\spe{g} \in \mathcal{C}$ be a graph on $N$, and let $\spe{h}$ be an induced connected graph on vertex set $S$. We know $\Delta_{S, N \setminus S}(\spe{g}) \neq 0$. Since $\varphi$ is balanced, it follows that $\varphi(\spe{h}) = \varphi(\spe{g}|_S) = 1$. Thus $\spe{h} \in \mathcal{C}$. Hence $\mathcal{C}$ is closed under inclusion of induced connected subgraphs.

Conversely, suppose that $\mathcal{C}$ is closed under inclusion of induced connected subgraphs, and contains $K_1$. Let $\spe{g} \in \spe{G}_N$ be a graph for which $\varphi_{\mathcal{C}}(\spe{g}) = 1$. Let $N = S \sqcup T$. Then we see that $\spe{g}|_S$ is a disjoint union of connected components $\spe{h}_i$. Moreover, each is an induced subgraph of a component of $\spe{g}$. Thus $\spe{h}_i$ is isomorphic to a graph in $\mathcal{C}$ for all $i$. Therefore $\varphi(\spe{g}|_S) = 1$. Similarly, $\varphi(\spe{g}/S) = 1$. By Proposition \ref{prop:criterion}, the character $\varphi_{\mathcal{C}}$ is balanced and convex. \end{proof}

As an example, let $d$ be an integer, and let $\mathcal{C}_d$ denote the set of isomorphism classes of connected graphs of maximum degree at most $d$. Then $\varphi_{\mathcal{C}_d}$ is balanced convex. Let us denote the character as $\varphi_d$. Then $\Psi_{\varphi_d}(\spe{g}, \mathbf{x})$ enumerates functions $f:N \to \mathbb{N}$ such that the induced subgraph on each color class has degree at most $d$. Such colorings are called \emph{defective} \cite{cowen-cowen-woodall}.

As another example, let $\mathcal{K} = \{K_1, K_2, \ldots \}$. Then $\varphi_{\mathcal{K}}$ is balanced convex, and $\Psi_{\varphi_{\mathcal{K}}}(\spe{g}, \mathbf{x})$ enumerates functions $f: N \to \mathbb{N}$ such that the induced subgraph on each color class is a disjoint union of cliques. These functions are called \emph{subcolorings} \cite{subchromatic}.

\subsection{Acyclic mixed graphs}
\label{subsec:mixed}
  Given a finite set $N$, a mixed graph is a triple $(N, U, D)$, where $U$ is a set of undirected edges, and $D$ is a set of directed edges.  A mixed graph is \emph{acyclic} if it does not contain a directed cycle.

There are two well-known polynomial invariants associated to acyclic mixed graphs: the weak and strong chromatic polynomial, both introduced in \cite{beck-et-al-2}, motivated by work in \cite{beck-et-al}. Given an acyclic mixed graph $\spe{g},$ a \emph{weak coloring} is a function $f:N \to \mathbb{N}$, subject to:
\begin{enumerate}
    \item For every $uv \in U$, we have $f(u) \neq f(v)$.
    \item For every $(u,v) \in D,$ we have $f(u) \leq f(v)$.
\end{enumerate} 
If we replace the second condition with strict inequalities, then we obtain the definition of \emph{strong coloring}.
the \emph{weak chromatic polynomial} $\chi(\spe{g}, k)$ counts the number of weak colorings $f$ such that $f(N) \subseteq [k].$ Similarly,
the \emph{strong} chromatic polynomial $\bar{\chi}(\spe{g}, k)$ counts the number of strong colorings $f$ with $f(N) \subseteq [k]$. Naturally there are quasisymmetric function generalizations of both these invariants, which we introduce. Moreover, both invariants can be understood from the theory of combinatorial Hopf monoids, and come from convex characters, as we now demonstrate.

Let $\spe{MG}_N$ be the set of all acyclic mixed graphs with vertex set $N$. 
Given a bijection $\sigma: N \to M$, and an acyclic mixed graph $\spe{g} = (N,U,D) \in \spe{MG}_N$, we define $\spe{MG}_{\sigma}(\spe{g}) = (M, U', D')$, where $U' = \{ \{\sigma(u), \sigma(v) \}: \{u, v \} \in U \}$ and $D' = \{ (\sigma(u), \sigma(v)): (u,v) \in D \}.$
We see that we get a map $\spe{MG}_{\sigma}: \spe{MG}_N \to \spe{MG}_M.$ Thus $\spe{MG}$ is a set species. 

Now we turn $\mathbb{K}\spe{MG}^{\bdot}$ into a linearized Hopf monoid. Given disjoint sets $M$ and $N$, and acyclic mixed graphs $\spe{g} \in \spe{MG}_M$ and $\spe{d}\in \spe{MG}_N$, we define $\spe{g} \cdot \spe{d}$ by taking disjoint unions of edge sets.

Likewise, given an acyclic mixed graph $\spe{g} = (M \sqcup N, U, D) \in \spe{MG}_{M \sqcup N}$, we define $\Delta_{M,N}(\spe{g})$ as follows. First, if there exists $(n,m) \in D$ with $m \in M, n \in N$, then $\Delta_{M,N}(\spe{g}) = 0$. Otherwise, we let $\spe{g}|_M$ be the induced subgraph on $M$, and $\spe{g}/M$ be the induced subgraph on $N$, and define $\Delta_{M,N}(\spe{g}) = \spe{g}|_M \otimes \spe{g}/M$.

There are two interesting examples of characters we will be interested in. The first is $\chi$. We see that $\chi(\spe{g}) = 1$ if and only if $\spe{g}$ has no edges. 
We define 
\[\varphi(\spe{g}) = 
\begin{cases} 
0 & \text{if $\spe{g}$ has at least one undirected edge} \\ 
1 & \text{otherwise.} 
\end{cases} \]
The reader may check that this is a multiplicative function.

\begin{theorem}
The pairs $(\mathbb{K}\spe{MG}^{\bdot}, \mathbb{K}(\varphi))$ and $(\mathbb{K}\spe{MG}^{\bdot}, \mathbb{K}(\chi))$ are linearized combinatorial Hopf monoids. Moreover, both $\varphi$ and $\chi$ are balanced convex characters.
\end{theorem}
\begin{proof}
Let $L, M$ and $N$ be finite sets. Let $\spe{g} \in \spe{MG}_L, \spe{h} \in \spe{MG}_M$, and $\spe{k} \in \spe{MG}_N$. Then $(\spe{g}\cdot \spe{h}) \cdot \spe{k} = \spe{g}\cdot (\spe{h} \cdot \spe{k})$, as the result in both cases is an acyclic mixed graph that is the disjoint union of the three individual graphs.

Let $\spe{f} \in \spe{MG}_{L \sqcup M \sqcup N}$. If there is a directed edge from a vertex in $L$ to a vertex in $M$, or from a vertex in $L$ to a vertex in $N$, or from a vertex in $M$ to a vertex in $N$, then 
\begin{align*}((\Delta_{L, M} \otimes \id_N) \circ \Delta_{L \sqcup M, N})(\spe{f}) & = 0 \\ & = ((\id_L \otimes \Delta_{M,N}) \circ \Delta_{L, M \sqcup N})(\spe{f}). \end{align*}

If no such directed edges exist, then 
\begin{align*}((\Delta_{L, M} \otimes \id_N) \circ \Delta_{L \sqcup M, N})(\spe{f}) & = \spe{f}|_L \otimes \spe{f}|_M \otimes \spe{f}|_N \\ & = ((\id_L \otimes \Delta_{M,N}) \circ \Delta_{L, M \sqcup N})(\spe{f}), \end{align*} and hence we have coassociativity.

Now, we show compatibility.
Let $A, B, C$ and $D$ be finite sets, and let $\spe{g} \in \spe{MG}_{A \sqcup B}$ and $\spe{h} \in \spe{MG}_{C \sqcup D}$. If there is a directed edge in $\spe{g}$ from a vertex in $A$ to a vertex in $B$, or there is a directed edge in $\spe{h}$ from a vertex in $C$ to a vertex in $D$, then 
\begin{align*} & \Delta_{A \sqcup C, B \sqcup D}(\spe{g} \cdot \spe{h}) = 0 \\ & = ((\mu_{A,C} \otimes \mu_{B,C}) \circ (\id_A \otimes \tau_{B,C} \otimes \id_D) \circ (\Delta_{A,B} \otimes \Delta_{C,D}))(\spe{g} \otimes \spe{h}). \end{align*}
Otherwise, 
\begin{align*} & \Delta_{A \sqcup C, B \sqcup D}(\spe{g} \cdot \spe{h}) = \spe{g}|_A \cdot \spe{h}|_C \otimes \spe{g}|_B \cdot \spe{h}|_D \\ & = ((\mu_{A,C} \otimes \mu_{B,C}) \circ (\id_A \otimes \tau_{B,C} \otimes \id_D) \circ (\Delta_{A,B} \otimes \Delta_{C,D}))(\spe{g} \otimes \spe{h}). \end{align*} Thus we see that the product and coproduct are compatible.

Now we show that $\varphi$ and $\chi$ are balanced convex characters. It is not hard to see that $\varphi$ is multiplicative. Given a mixed graph $\spe{g} \in \spe{MG}_N$, with $|N| = 1$, then $\spe{g}$ is just a vertex, and we have $\varphi(\spe{g}) = 1$.

Given $\spe{g} \in \spe{MG}_N$ with $|N| \geq 2$, let $v$ be a sink of $\spe{g}$ with respect to the directed edges. Then $\Delta_{\{v\}, N \setminus \{v\}}(\spe{g}) \neq 0$.

Now suppose that $\spe{g}$ is an acyclic mixed graph where $\varphi(\spe{g}) = 1$. Then $\spe{g}$ has no undirected edges. Let $V = S \sqcup T$, and suppose $\Delta_{S,T}(\spe{g}) \neq 0$. Then $\Delta_{S,T}(\spe{g}) = \spe{g}|_S \otimes \spe{g} / S$. Both $\spe{g}|_S$ and $\spe{g}/S$ are induced subgraphs, and hence have no undirected edges. Thus $\varphi(\spe{g}|_S) = \varphi(\spe{g}/S) = 1$. Thus by Proposition \ref{prop:criterion}, the character $\varphi$ is balanced convex.

The fact that $\chi$ is balanced convex follows from Theorem \ref{thm:chiconvex}. \end{proof}
Given an acyclic mixed graph $\spe{g}$ with vertex set $V$, we see that $\Psi_{\varphi}(\spe{g}, \mathbf{x})$ is a summation over all weak colorings of $\spe{g}$. 
Similarly, $\Psi_{\chi}(\spe{g}, \mathbf{x})$ is a summation over strong colorings of $\spe{g}.$

As an example, if we let $\spe{g}$ be the graph on the left in Figure \ref{fig:mixedgraph}, then 
\begin{align*} \Psi_{\varphi}(\spe{g}, \mathbf{x}) & = 3M_{2,2} + 4M_{1,1,2} + 2M_{1,2,1}+4M_{2,1,1} + 6M_{1,1,1,1} \\ 
& = 3F_{2,2} + F_{1,1,2} + 2F_{1,2,1}+F_{2,1,1}-F_{1,1,1,1}. \end{align*}
We have that \[\Psi_{\chi}(\spe{g}, \mathbf{x}) = M_{2,2} + 2M_{1,1,2} + 2M_{2,1,1} + 6M_{1,1,1,1}. \] Note that $\Psi_{\varphi}(\spe{g}, \mathbf{x})$ is not $F$-positive. 
The relative simplicial complex in the middle of Figure \ref{fig:mixedgraph} is the coloring complex of $\spe{g}$ with respect to $\varphi$, while the relative simplicial complex on the right is the coloring complex of $\spe{g}$ with respect to $\chi$.

\begin{corollary}
For every acyclic mixed graph $\spe{g}$, the quasisymmetric functions $\Psi_{\chi}(\spe{g}, \mathbf{x})$ and $\Psi_{\varphi}(\spe{g}, \mathbf{x})$ are both $M$-increasing. Moreover, the corresponding weak and strong chromatic polynomials are both super flawless.
\end{corollary}

We will show in a future paper that $\Psi_{\chi}(\spe{g})$ is $F$-positive by showing that the relative simplicial complex is relatively shellable. It is an interesting problem to determine necessary and sufficient criteria for when $\Sigma_{\varphi}(\spe{g})$ is relatively shellable.

\begin{figure}
\begin{center}
\begin{tabular}{ccc}
\begin{tikzpicture}
  \node[circle, draw=black, fill=white] (b) at (0,2) {$b$};
  \node[circle, draw=black, fill=white] (a) at (0,0) {$a$};
  \node[circle, draw=black, fill=white] (c) at (1.5,0) {$c$};
  \node[circle, draw=black, fill=white] (d) at (1.5,2) {$d$};
 \draw[-Latex] (b) -- (a);
  \draw[-Latex] (d) edge (c);
  \draw[-] (b) -- (c);
    \draw[-] (d) -- (a);

\end{tikzpicture}
&

\begin{tikzpicture}
\draw[color=white, fill=gray!20] (90:1.5cm) -- (30:1.5cm) -- (330:1.5cm) -- (270:1.5cm) -- (210:1.5cm) -- (150:1.5cm) -- cycle;
\node[circle, draw=black, fill=white, thick] (c) at (90:1.5cm) {$ab$};
  \node[circle, draw=red, fill=white, dashed, thick] (b) at (150:1.5cm) {$abc$};
  \node[circle, draw=red, fill=white, dashed, thick] (a) at (210:1.5cm) {$c$};
  \node[circle, draw=black, fill=white, thick] (f) at (270:1.5cm) {$cd$};
  \node[circle, draw=red, fill=white, dashed, thick] (e) at (330:1.5cm) {$acd$};
  \node[circle, draw=red, fill=white,dashed, thick] (d) at (30:1.5cm) {$a$};
  \node[circle, draw=black, fill=white,thick] (or) at (0:0cm) {$ac$};
  \draw[thick] (a) -- (b) -- (c) -- (d) -- (e) -- (f) -- (a);
  \draw[dashed, red, thick] (b) -- (d);
  \draw[dashed, red, thick] (a) -- (e);
  \draw[thick] (b) -- (or) -- (d);
  \draw[thick] (a) -- (or) -- (e);

\end{tikzpicture}

&

\begin{tikzpicture}
\draw[color=white, fill=gray!20] (90:1.5cm) -- (30:1.5cm) -- (330:1.5cm) -- (270:1.5cm) -- (210:1.5cm) -- (150:1.5cm) -- cycle;
\node[circle, draw=red, fill=white,dashed, thick] (c) at (90:1.5cm) {$ab$};
  \node[circle, draw=red, fill=white, dashed, thick] (b) at (150:1.5cm) {$abc$};
  \node[circle, draw=red, fill=white, dashed, thick] (a) at (210:1.5cm) {$c$};
  \node[circle, draw=red, fill=white,dashed, thick] (f) at (270:1.5cm) {$cd$};
  \node[circle, draw=red, fill=white,dashed, thick] (e) at (330:1.5cm) {$acd$};
  \node[circle, draw=red, fill=white,dashed, thick] (d) at (30:1.5cm) {$a$};
  \node[circle, draw=black, fill=white,thick] (or) at (0:0cm) {$ac$};
  \draw[dashed, red, thick] (a) -- (b) -- (c) -- (d) -- (e) -- (f) -- (a);
  \draw[dashed, red, thick] (b) -- (d);
  \draw[dashed, red, thick] (a) -- (e);
  \draw[thick] (b) -- (or) -- (d);
  \draw[thick] (a) -- (or) -- (e);

\end{tikzpicture}

\end{tabular}
\end{center}
\caption{an acyclic mixed graph, and its two coloring complexes with respect to different characters.}
\label{fig:mixedgraph}
\end{figure}

\subsection{Rooted Connected Graphs}
\label{subsec:rooted}
  Given a finite set $N$, a rooted connected graph is a connected graph on $N \sqcup \{r \}$. We view $r$ as an extra vertex, which we call the root.

Let $\spe{RG}_N$ be the set of all rooted connected graphs with vertex set $N \sqcup \{r\}$. 
Given a bijection $\sigma: N \to M$, we extend to a bijection $\sigma: N \sqcup \{r\} \to M \sqcup \{r\}$ by defining $\sigma(r) = r$. Given a rooted graph $\spe{g} \in \spe{RG}_N$, we define $\spe{RG}_{\sigma}(\spe{g})$ by $E(\spe{RG}_{\sigma}(\spe{g})) = \{ \{\sigma(u), \sigma(v) \}: \{u,v\} \in E(\spe{g}) \}$.
We see that we get a map $\spe{RG}_{\sigma}: \spe{RG}_N \to \spe{RG}_M.$ Thus $\spe{RG}$ is a set species. 

Now we turn $\mathbb{K}\spe{RG}^{\bdot}$ into a linearized Hopf monoid. Given disjoint sets $M$ and $N$, and rooted connected graphs $\spe{g} \in \spe{RG}_M$ and $\spe{h}\in \spe{RG}_N$, we define $\spe{g} \cdot \spe{h}$ by taking disjoint unions of edge sets. In this case $\spe{g} \cdot \spe{h}$ has vertex set $M \sqcup N \sqcup \{r\}$, and we are identifying the roots of $\spe{g}$ and $\spe{h}$. As a result, the product is a rooted connected graph.

Likewise, let $\spe{g}  \in \spe{RG}_{M \sqcup N}$ such that $\spe{g}|_{M \cup r}$, the induced subgraph on $M \cup r$, is connected. Let $W = \{v \in N: uv \in E(\spe{g}) \mbox{ for some } u \in M \}$. We define a graph $\spe{g}/M$ on $N \cup \{r\}$ by saying $uv \in E(\spe{g}/M)$ if and only if $u, v \in N$ and $uv \in E(\spe{g})$ or $u = r$ and $v \in W$. Equivalently, $\spe{g}/M$ is the graph obtained from $\spe{g}$ by contracting all the edges of $\spe{g}|_{M \cup \{r\}}$. We claim that $\spe{g}/M \in \spe{RG}_N$. To see this, let $u \in N$. Then there is a path $u, v_1, \ldots, v_k, r$. Let $i$ be the first index where $v_{i+1} \in M \cup \{r\},$ if such an index exists. Otherwise, let $i=r.$ Then we see that $u, v_1, \ldots, v_i, r$ is a path in $\spe{g}/M$, and thus every $\spe{g}/M$ is connected, and rooted at $r.$

Given $\spe{g} \in \spe{RG}_{M \sqcup N}$, we define 
\[\Delta_{M,N}(\spe{g}) = \begin{cases} \spe{g}|_{M \cup \{r\}} \otimes \spe{g}/M & \spe{g}|_{M \cup \{r\}} \mbox{ is connected } \\ 0 & \mbox{ otherwise } \end{cases} \]

Given a rooted connected graph $\spe{g}$, we see that $\chi(\spe{g}) = 1$ if and only if $\spe{g}$ is a star graph rooted at the only nonleaf vertex.

\begin{theorem}
The pair $(\mathbb{K}\spe{RG}, \mathbb{K}(\chi))$ is a linearized combinatorial Hopf monoid. Moreover, $\chi$ is a balanced convex character.
\end{theorem}
\begin{proof}
Let $L, M$ and $N$ be finite sets. Let $\spe{g} \in \spe{RG}_L, \spe{h} \in \spe{RG}_M$, and $\spe{k} \in \spe{RG}_N$. Then $(\spe{g}\cdot \spe{h}) \cdot \spe{k} = \spe{g}\cdot (\spe{h} \cdot \spe{k})$, as the result in both cases is a rooted connected graph that is the union of the three individual graphs, with the corresponding roots identified.

Let $\spe{g} \in \spe{RG}_{L \sqcup M \sqcup N}$. If $\spe{g}|_{L \cup \{r\}}$ is disconnected, or $\spe{g}|_{L \sqcup M \sqcup \{r\}}$ is disconnected, then 
\begin{align*}((\Delta_{L, M} \otimes \id_N) \circ \Delta_{L \sqcup M, N})(\spe{g}) & = 0 \\ & = ((\id_L \otimes \Delta_{M,N}) \circ \Delta_{L, M \sqcup N})(\spe{g}). \end{align*}
Otherwise, we can see that 
\begin{align*} ((\Delta_{L, M} \otimes \id_N) \circ \Delta_{L \sqcup M, N})(\spe{g}) & = \spe{g}|_{L \sqcup \{r\}} \otimes (\spe{g}|_{L \sqcup M \sqcup \{r\}}) / L \otimes \spe{g} / (L \sqcup M) .\end{align*}  
Similarly, \[((\id_L \otimes \Delta_{M,N}) \circ \Delta_{L, M \sqcup N})(\spe{g}) = \spe{g}|_{L \sqcup \{r\}} \otimes (\spe{g} / L)|_{M \sqcup \{r\}} \otimes \spe{g} / (L \cup M). \] We see that $(\spe{g}|_{L \sqcup M \sqcup \{r\}}) / L = (\spe{g} / L)|_{M \sqcup \{r\}}$, as both graphs consist of restricting to $M$ and contracting the vertices in $L \sqcup \{r\}$. Hence we have coassociativity.

Now we show compatibility.
Let $A, B, C$ and $D$ be finite sets, and let $\spe{g} \in \spe{RG}_{A \sqcup B}$ and $\spe{h} \in \spe{RG}_{C \sqcup D}$. If $\spe{g}|_{A \cup \{r\}}$ or $\spe{h}|_{C \cup \{r\}}$ is disconnected, then 
\begin{align*} & \Delta_{A \sqcup C, B \sqcup D}(\spe{g} \cdot \spe{h}) = 0 \\ &  = ((\mu_{A,C} \otimes \mu_{B,C}) \circ (\id_A \otimes \tau_{B,C} \otimes \id_D) \circ (\Delta_{A,B} \otimes \Delta_{C,D}))(\spe{g} \otimes \spe{h}). \end{align*}
Otherwise, we see that 
\begin{align*} & \Delta_{A \sqcup C, B \sqcup D}(\spe{g} \cdot \spe{h}) = \spe{g}|_{A \cup \{r\}} \cdot \spe{h}|_{C \cup \{r\}} \otimes \spe{g}/A \cdot \spe{h}/C \\ & = (\mu_{A,C} \otimes \mu_{B,C}) \circ (\id_A \otimes \tau_{B,C} \otimes \id_D) \circ (\Delta_{A,B} \otimes \Delta_{C,D}))(\spe{g} \otimes \spe{h}) .\end{align*} Thus we see that the product and coproduct are compatible.

Now we show that $\chi$ is a balanced convex character. For $\spe{g} \in \spe{RG}_N$ with $|N| \geq 2$, let $v$ be a vertex that is adjacent to $r$. Then $\Delta_{\{v \}, N \setminus \{v\}}(\spe{g}) \neq 0$.
By Theorem \ref{thm:chiconvex}, $\chi$ is a balanced convex character. \end{proof}
Given a rooted connected graph $\spe{g}$ with vertex set $V$, we see that $\Psi_{\chi}(\spe{g}, \mathbf{x})$ is a summation over proper colorings that have the additional feature that, for every vertex $v$, there exists a path $r, v_1, \ldots, v_k, v$ with the property that $f(v_1) < f(v_2) < \cdots < f(v_k) < f(v)$. If there is an edge between $r$ and every vertex of $N$, then $\Psi_{\chi}(\spe{g}, \mathbf{x})$ is just the chromatic symmetric function of the unrooted graph $\spe{g} \setminus \{r\}$. However, for other types of rooted connected graphs, we obtain a new quasisymmetric function.

As an example, if we let $\spe{g}$ be the graph on the left in Figure \ref{fig:rootedgraph}, then \[\Psi_{\chi}(\spe{g}, \mathbf{x}) = M_{2,1} + 4M_{1,1,1}. \]  
The relative simplicial complex on the right of Figure \ref{fig:rootedgraph} is the coloring complex of $\spe{g}$ with respect to $\chi$.
\begin{corollary}
For every rooted connected graph $\spe{g}$, the quasisymmetric function $\Psi_{\chi}(\spe{g}, \mathbf{x})$ is $M$-increasing. Moreover, the corresponding chromatic polynomial is super flawless.
\end{corollary}

\begin{figure}
\begin{center}
\begin{tabular}{cc}

\begin{tikzpicture}

  \node[circle, draw=black, fill=white, thick] (a) at (210:2cm) {$a$};
  \node[circle, draw=black, fill=white, thick] (f) at (270:2cm) {$r$};
  \node[circle, draw=black, fill=white, thick] (e) at (330:2cm) {$c$};
  \node[circle, draw=black, fill=white,thick] (or) at (0:0cm) {$b$};
  \draw[thick] (a) -- (or) -- (e) -- (f) -- (a);

\end{tikzpicture}

&

\begin{tikzpicture}

  \node[circle, draw=red, fill=white, dashed, thick] (b) at (150:2cm) {$a, b$};
  \node[circle, draw=red, fill=white, dashed, thick] (a) at (210:2cm) {$a$};

  \node[circle, draw=red, fill=white, dashed, thick] (e) at (330:2cm) {$c$};
  \node[circle, draw=red, fill=white, dashed, thick] (d) at (30:2cm) {$b,c$};
  \node[circle, draw=black, fill=white,thick] (or) at (0:0cm) {$a, c$};
 
  \draw[thick] (b) -- (a) -- (or) -- (e) -- (d);

\end{tikzpicture}

\end{tabular}
\end{center}
\caption{A rooted connected graph, and its corresponding coloring complex.}
\label{fig:rootedgraph}
\end{figure}

\subsection{Double Posets}
\label{subsec:double}
Now we will discuss double posets. The Hopf algebra of double posets was introduced by Malvenuto and Reutenauer \cite{malvenuto}. However, we show that double posets also form a linearized combinatorial Hopf monoid. The associated quasisymmetric function is a generalization of Gessel's $(P, \omega)$-partition enumerator for labeled posets $(P, \omega)$. The corresponding quasisymmetric function for double posets is studied extensively by Grinberg \cite{grinberg}.

  Given a finite set $N$, a double poset on $N$ is a triple $(N, \leq_1, \leq_2)$ where $\leq_1$ and $\leq_2$ are both partial orders on $N$.
Let $\spe{DP}_N$ be the set of all double posets with vertex set $N$. 
Given a bijection $\sigma: N \to M$, and a double poset $\spe{d} \in \spe{DP}_N$, we define $\spe{DP}_{\sigma}(\spe{d})$ to be the double poset given by declaring $x \leq_i y$ if and only if $\sigma^{-1}(x) \leq_i \sigma^{-1}(y)$ for all $x,y \in N, i \in \{1,2 \}$.
We see that we get a map $\spe{DP}_{\sigma}: \spe{DP}_N \to \spe{DP}_M.$ Thus $\spe{DP}$ is a set species. 

Now we turn $\mathbb{K}\spe{DP}^{\bdot}$ into a linearized Hopf monoid. Given disjoint sets $M$ and $N$, and double posets $\spe{d} \in \spe{DP}_M$ and $\spe{h}\in \spe{DP}_N$, we define two partial orders on $M \sqcup N$.
For $x, y \in M \sqcup N$, we say $x \leq_1 y$ if one of the following holds:
\begin{enumerate}
    \item $x, y \in M$ and $x \leq_1 y$ in $\spe{d}$.
    \item $x,y \in N$ and $x \leq_1 y$ in $\spe{h}$.
\end{enumerate}
We observe that $\leq_1$ is the \emph{disjoint union} of two partial orders.
For $x, y \in M \sqcup N$, we say $x \leq_2 y$ if one of the following holds:
\begin{enumerate}
    \item $x, y \in M$ and $x \leq_2 y$ in $\spe{d}$.
    \item $x,y \in N$ and $x \leq_2 y$ in $\spe{h}$.
        \item $x \in M$ and $y \in N$.
\end{enumerate}
We observe that $\leq_2$ is the \emph{ordinal sum} of two partial orders.
Then we let $\spe{d} \cdot \spe{h} = (M \sqcup N, \leq_1, \leq_2).$ Given $\spe{d} \in \spe{DP}_M$ and $\spe{h} \in \spe{DP}_N$, where $M$ and $N$ are disjoint sets, we define a double poset $\spe{d} \cdot \spe{h} = (M \sqcup N, \leq_1, \leq_2)$. 
Hence we have a multiplication operation for $\spe{D}$.

Now we define the comultiplication.
Given a double poset $\spe{d} \in \spe{DP}_{M \sqcup N}$, we define $\spe{d}|_M$ by: for $x, y \in M$, we say $x \leq_i y$ in $\spe{d}|_M$ if and only if $x \leq_i y$ in $\spe{d}$.
We define \[\Delta_{M,N}(\spe{d}) = \begin{cases} \spe{d}|_M \otimes \spe{d}|_N& \text{ if $M$ is an order ideal of $\leq_1$} \\ 0 & \text{ otherwise. } \end{cases}\]
Since we are going to prove that $\Delta$ is coassociative, it is worth noting that, given a poset $\spe{p} \in \spe{P}_{M \sqcup N}$, then $M$ is a $\leq_1$-order ideal if and only if $N$ is a $\leq_1$-order filter.

Given a double poset $\spe{d}$, a pair $(m,m') \in M$ is an \emph{inversion} if $m <_1 m'$ and $m' <_2 m$.
Finally, given a double poset $\spe{d} \in \spe{DP}_N$, we define 
\[\varphi(\spe{d}) = 
\begin{cases} 
0 & \text{if $\spe{d}$ has an inversion} \\ 
1 & \text{otherwise.} 
\end{cases} \]
The reader may check that this is a multiplicative function.
\begin{theorem}
The pairs $(\mathbb{K}\spe{DP}^{\bdot}, \mathbb{K}(\varphi))$ and $(\mathbb{K}\spe{DP}^{\bdot}, \mathbb{K}(\chi))$ are linearized combinatorial Hopf monoids. Moreover, $\varphi$ and $\chi$ are balanced convex characters.
\end{theorem}
\begin{proof}
Let $L, M$ and $N$ be finite sets. Let $\spe{d} \in \spe{DP}_L, \spe{f} \in \spe{DP}_M$, and $\spe{g} \in \spe{DP}_N$. Then $(\spe{d}\cdot \spe{f}) \cdot \spe{g} = \spe{d}\cdot (\spe{f} \cdot \spe{g})$, as the result in both cases is a double poset such that $\leq_1$ is the disjoint union of the first partial orders for $\spe{d}$, $\spe{f}$ and $\spe{g}$, and $\leq_2$ is the ordinal sum of the second partial orders of $\spe{d}, \spe{f}$ and $\spe{g}$.

Let $\spe{d} \in \spe{DP}_{L \sqcup M \sqcup N}$. If $L$ is not a $\leq_1$-order ideal of $\spe{d}$, or $N$ is not a $\leq_1$-order filter of $\spe{d}$, then \begin{align*} ((\Delta_{L, M} \otimes \id_N) \circ \Delta_{L \sqcup M, N})(\spe{f}) & = 0  \\ & = ((\id_L \otimes \Delta_{M,N}) \circ \Delta_{L, M \sqcup N})(\spe{f}). \end{align*}
Otherwise, we see that $M$ is a $\leq_1$-convex subset of $\spe{d}$, and we have
\begin{align*} ((\Delta_{L, M} \otimes \id_N) \circ \Delta_{L \sqcup M, N})(\spe{f}) & = \spe{f}|_L \otimes \spe{f}|_M \otimes \spe{f}|_N  \\ & = ((\id_L \otimes \Delta_{M,N}) \circ \Delta_{L, M \sqcup N})(\spe{f}). \end{align*}
Hence we have coassociativity.

Now we show compatibility.
Let $A, B, C$ and $D$ be finite sets, and let $\spe{d} \in \spe{DP}_{A \sqcup B}$ and $\spe{f} \in \spe{DP}_{C \sqcup D}$. If $A$ is not a $\leq_1$-order ideal of $\spe{d}$ or $C$ is not a $\leq_1$-order ideal of $\spe{f}$, then $A \cup C$ is not a $\leq_1$-order ideal of $\spe{d} \cdot \spe{f}$. In that case, we have
\begin{align*} & \Delta_{A \sqcup C, B \sqcup D}(\spe{d} \cdot \spe{f}) = 0 \\ & = ((\mu_{A,C} \otimes \mu_{B,C}) \circ (\id_A \otimes \tau_{B,C} \otimes \id_D) \circ (\Delta_{A,B} \otimes \Delta_{C,D}))(\spe{d} \otimes \spe{f}). \end{align*} 
If $A$ and $C$ are $\leq_1$-order ideals of their corresponding posets, then $A \cup C$ is a $\leq_1$-order ideal of $\spe{d} \cdot \spe{f}$.
Write $\Delta_{A \sqcup C, B \sqcup D}(\spe{d} \cdot \spe{f}) = \spe{h} \otimes \spe{k}$. Then $(\spe{h}, \leq_1)$ is the disjoint union of $(\spe{d}|_A, \leq_1)$ and $(\spe{f}|_C, \leq_1)$, while $(\spe{h}, \leq_2)$ is the ordinal sum of $(\spe{d}|_A, \leq_2)$ and $(\spe{f}|_C, \leq_2)$. Similarly, $(\spe{k}, \leq_1)$ is the disjoint union of $(\spe{d}|_B, \leq_1)$ and $(\spe{f}|_D, \leq_1)$, while $(\spe{h}, \leq_2)$ is the ordinal sum of $(\spe{d}|_B, \leq_2)$ and $(\spe{f}|_D, \leq_2)$. It follows that 
\[\spe{h} \otimes \spe{k} = ((\mu_{A,C} \otimes \mu_{B,C}) \circ (\id_A \otimes \tau_{B,C} \otimes \id_D) \circ (\Delta_{A,B} \otimes \Delta_{C,D}))(\spe{d} \otimes \spe{f}). \]

Now we show that the characters $\varphi$ and $\chi$ are balanced convex. Let $v \in N$ be a minimal element with respect to $\leq_1$. Then $\Delta_{\{v\}, N \setminus \{v\}}(\spe{d}) \neq 0$. By Theorem \ref{thm:chiconvex}, it follows that $\chi$ is balanced convex.

Also, given $\spe{d} \in \spe{DP}_N$, where $|N| = 1$, we see that $\spe{d}$ consists of a single vertex, and thus $\varphi(\spe{d}) = 1$. 

Let $\spe{d}$ be a double poset such that $\varphi(\spe{d}) = 1$. Let $N = S \sqcup T$ such that $\Delta_{S,T}(\spe{d}) \neq 0$. Then $\Delta_{S,T}(\spe{d}) = \spe{d}|_S \otimes \spe{d}|_T$, where each factor is the induced double poset on the corresponding subset. We see that if $\spe{d}$ has no inversions, then every induced subposet of $\spe{d}$ also has no inversions. Hence $\varphi_{S|T}(\spe{d}) = 1$. Thus, by Proposition \ref{prop:criterion}, $\varphi$ is a balanced convex character. \end{proof}
The corresponding quasisymmetric function $\Psi_{\varphi}(\spe{d}, \mathbf{x})$ enumerates double poset partitions. A double poset partition is a function $\sigma: N \to \mathbb{N}$ subject to:
\begin{enumerate}
    \item for $x, y \in N$, if $x \leq_1 y$, then $\sigma(x) \leq \sigma(y)$.
    \item for $x, y \in N$, if $x \leq_1 y$ and $y <_2 x$, then $\sigma(x) < \sigma(y)$.
\end{enumerate}

It is not hard to show that $\Psi_{\varphi}(\spe{d}, \mathbf{x}) = \sum_{\sigma} \prod_{m \in M} x_{\sigma(m)}$. This amounts to showing that a double poset partition is a $\varphi$-proper function.

As an example, if we let $\spe{d}$ be the double poset in Figure \ref{fig:doubleposet}, where the Hasse diagram on the left is for $\leq_1$ and the Hasse diagram in the middle is for $\leq_2$. Then 
\begin{align*} \Psi_{\varphi}(\spe{d}, \mathbf{x}) & = M_{2,2} + 2M_{1,1,2} + 2M_{1,2,1}+2M_{2,1,1} + 4M_{1,1,1,1} \\ & = F_{2,2}+F_{1,1,2}+2F_{1,2,1}+F_{2,1,1}-F_{1,1,1,1}. \end{align*} This quasisymmetric function is not $F$-positive.
The relative simplicial complex on the right of Figure \ref{fig:doubleposet} is the coloring complex of $\spe{d}$ with respect to $\varphi$.
\begin{corollary}
For every double poset $\spe{d}$, the quasisymmetric functions $\Psi_{\varphi}(\spe{p}, \mathbf{x})$ is $M$-increasing. Moreover, the corresponding order polynomial is strongly flawless.
\end{corollary}
It is an interesting problem to determine necessary and sufficient criteria for when $\Sigma_{\varphi}(\spe{d})$ is relatively shellable.

We also briefly mention the character $\chi$. Given a double poset $\spe{d}$, define the mixed graph $\spe{g}(\spe{d})$ as follows:
\begin{enumerate}
    \item For $u, v \in N$, we have $(u,v) \in E(\spe{g}(\spe{d}))$ if and only if $u \leq_1 v$.
    \item For $u, v \in N$, we have $uv \in E(\spe{g}(\spe{d}))$ if and only if $u, v$ are $\leq_2$-incomparable.
\end{enumerate}
Given a finite set $N$, let $\psi: \spe{DP}_N \to \spe{MG}_N$ be given by $\psi(\spe{d}) = \spe{g}(\spe{d})$. One can show that $\mathbb{K}(\psi)$ is a morphism of combinatorial Hopf monoids with respect to the character $\chi$. Thus, while at first glance $\Psi_{\chi}(\spe{d}, \mathbf{x})$ appears to be a new invariant, it is enumerating strong colorings of an associated mixed graph.

\begin{figure}
\begin{center}
\begin{tabular}{ccc}
\begin{tikzpicture}
  \node[circle, draw=black, fill=white] (b) at (0,2) {$b$};
  \node[circle, draw=black, fill=white] (a) at (0,0) {$a$};
  \node[circle, draw=black, fill=white] (c) at (2,0) {$c$};
  \node[circle, draw=black, fill=white] (d) at (2,2) {$d$};
 \draw[-Latex] (b) -- (a);
  \draw[-Latex] (d) edge (c);
  \draw[-Latex] (b) -- (c);
    \draw[-Latex] (d) -- (a);

\end{tikzpicture}

& 

\begin{tikzpicture}
  \node[circle, draw=black, fill=white] (b) at (0,2) {$c$};
  \node[circle, draw=black, fill=white] (a) at (0,0) {$b$};
  \node[circle, draw=black, fill=white] (c) at (2,0) {$d$};
  \node[circle, draw=black, fill=white] (d) at (2,2) {$a$};
 \draw[-Latex] (b) -- (a);
  \draw[-Latex] (d) edge (c);

\end{tikzpicture}

&

\begin{tikzpicture}
\draw[color=white, fill=gray!20] (30:2cm) -- (330:2cm) -- (210:2cm) -- (150:2cm) -- cycle;
  \node[circle, draw=red, fill=white, dashed, thick] (b) at (150:2cm) {$abc$};
  \node[circle, draw=red, fill=white, dashed, thick] (a) at (210:2cm) {$c$};
  \node[circle, draw=red, fill=white, dashed, thick] (e) at (330:2cm) {$acd$};
  \node[circle, draw=red, fill=white,dashed, thick] (d) at (30:2cm) {$a$};
  \node[circle, draw=black, fill=white,thick] (or) at (0:0cm) {$ac$};
  \draw[thick] (a) -- (b);
  \draw[thick] (d) -- (e);

  \draw[dashed, red, thick] (b) -- (d);
  \draw[dashed, red, thick] (a) -- (e);
  \draw[thick] (b) -- (or) -- (d);
  \draw[thick] (a) -- (or) -- (e);

\end{tikzpicture}

\end{tabular}
\end{center}
\caption{A double poset, and its coloring complex.}
\label{fig:doubleposet}
\end{figure}

\subsection{Antimatroids}
\label{subsec:antimatroid}

  Given a finite set $N$, an \emph{antimatroid} is a non-empty collection $\spe{a}$ of subsets of $N$ that satisfies the following conditions:
\begin{enumerate}
    \item For every $S \in \spe{a}$, there exists $x \in S$ such that $S \setminus \{x\} \in \spe{a}.$
    \item For every $S, T \in \spe{a}$, we have $S \cup T \in \spe{a}$.
\end{enumerate}
Antimatroids were introduced by Jamison \cite{antimatroid}. Given a partially ordered set $P$, let $J(P)$ be the set of order ideals of $P$. Then $J(P)$ is an antimatroid. Another example of an antimatroid appears on the left in Figure \ref{fig:antimatroid}.

Let $\spe{A}_N$ be the set of all antimatroids with vertex set $N$. 
Given a bijection $\sigma: N \to M$, and an antimatroid $\spe{a} \in \spe{A}_N$, we define $\sigma(\spe{a}) = \{\sigma(S): S \in \spe{a} \}$.
We see that we get a map $\spe{A}_{\sigma}: \spe{A}_N \to \spe{A}_M.$ Thus $\spe{A}$ is a set species. 

Now we turn $\mathbb{K}\spe{A}^{\bdot}$ into a linearized Hopf monoid. Given disjoint sets $M$ and $N$, and antimatroids $\spe{a} \in \spe{A}_M$ and $\spe{b}\in \spe{A}_N$, we define $\spe{a} \cdot \spe{b} = \{X \cup Y: X \in \spe{a}, Y \in \spe{b} \}$.

Likewise, given an antimatroid $\spe{a} \in \spe{A}_{M \sqcup N}$ with $M \in \spe{a}$, we define $\spe{a}|_M = \{S: S \in \spe{a}, S \subseteq M \}$ and $\spe{a}/M = \{T: T \subseteq N, T \cup M \in \spe{a} \}$. These operations are called the \emph{restriction} and \emph{contraction}, respectively. It is easy to see that both the restriction and contraction are closed under unions, and that for $S \in \spe{a}|_M$, there exists $x \in S$ such that $S \setminus \{x\} \in \spe{a}|_M$. Let $S \in \spe{a}/M$. Then $S \cup M \in \spe{a}$. There exists $x_1 \in S \cup M$ such that $S \cup M \setminus \{x_1\} \in \spe{a}$. If $x_1 \in M$, then there exists $x_2 \in S \cup M \setminus \{x_1\}$ such that $ S \cup M \setminus \{x_1, x_2\}\in \spe{a}$. Continuing, there is a sequence of elements $x_1, \ldots, x_i$, with $x_1, \ldots, x_i \in M$ and $S \cup M \setminus \{x_1, \ldots, x_i \} \in \spe{a}$. Since $M$ is a finite set, we see that eventually we find $T \subseteq M$ and $x \in S$ such that $S \cup M \setminus (T \cup \{x\}) \in \spe{a}$. Since $\spe{a}$ is closed under unions, $M \cup ( S \cup M \setminus (T \cup \{x\})) \in \spe{a}$. Thus $(S \setminus \{x\}) \cup M \in \spe{a}$ and $S \setminus \{x\} \in \spe{a}/M$. Therefore, the restriction and contraction of an antimatroid is an antimatroid.

For $\spe{a} \in \spe{A}_{M \sqcup N}$, we define 
\[\Delta_{M,N}(\spe{a}) = 
\begin{cases}
\spe{a}|_M \otimes \spe{a}/M & \text{ if $M \in \spe{a}$} \\
0 & \text{ otherwise.}
\end{cases}
\]
Finally, given an antimatroid $\spe{a} \in \spe{A}_N$, we see that $\chi(\spe{a}) = 1$ if and only if $\spe{a} = 2^N$.

\begin{theorem}
The pair $(\mathbb{K}\spe{A}, \mathbb{K}(\chi))$ is a linearized combinatorial Hopf monoid. Moreover, $\chi$ is a balanced convex character.
\end{theorem}
\begin{proof}
Let $L, M$ and $N$ be finite sets. Let $\spe{a} \in \spe{A}_L, \spe{b} \in \spe{A}_M$, and $\spe{c} \in \spe{A}_N$. Then $(\spe{a}\cdot \spe{b}) \cdot \spe{c} = \spe{a}\cdot (\spe{b} \cdot \spe{c})$, as the result in both cases is the antimatroid \[\spe{a} \cdot \spe{b} \cdot \spe{c} = \{X \cup Y \cup Z: X \in \spe{a}, Y \in \spe{b}, Z \in \spe{c} \}. \]

Next, we show coassociativity. Let $\spe{a} \in \spe{A}_{L \sqcup M \sqcup N}$. If $L \not\in \spe{a}$ or $L \cup M \not\in \spe{a}$, then \[((\Delta_{L, M} \otimes \id_N) \circ \Delta_{L \sqcup M, N})(\spe{a}) = 0 = ((\id_L \otimes \Delta_{M,N}) \circ \Delta_{L, M \sqcup N})(\spe{a}). \] Otherwise, we have 
\begin{align*}((\Delta_{L, M} \otimes \id_N) \circ \Delta_{L \sqcup M, N})(\spe{a}) & = \spe{a}|_L \otimes \spe{a}_{L,M} \otimes \spe{a}/(L \cup M)  \\ & = ((\id_L \otimes \Delta_{M,N}) \circ \Delta_{L, M \sqcup N})(\spe{a}) \end{align*}
 where $\spe{a}_{L,M} = \{Y \subseteq M: Y \cup L \in \spe{a} \}$. Hence we have coassociativity.

Now we show compatibility.
Let $A, B, C$ and $D$ be finite sets, and let $\spe{a} \in \spe{A}_{A \sqcup B}$ and $\spe{b} \in \spe{A}_{C \sqcup D}$. If $A \not\in \spe{a}$ or $C \not\in \spe{b}$, then 
\begin{align*} & \Delta_{A \sqcup C, B \sqcup D}(\spe{a} \cdot \spe{b})  = 0 \\  & = ((\mu_{A,C} \otimes \mu_{B,C}) \circ (\id_A \otimes \tau_{B,C} \otimes \id_D) \circ (\Delta_{A,B} \otimes \Delta_{C,D}))(\spe{a} \otimes \spe{b}). \end{align*}
In all other cases, we have 
\begin{align*} & \Delta_{A \sqcup C, B \sqcup D}(\spe{a} \cdot \spe{b})  = \spe{a}|_A \cdot \spe{b}|_C \otimes \spe{a}/A \cdot \spe{b}/C \\ 
& = (\mu_{A,C} \otimes \mu_{B,C}) \circ (\id_A \otimes \tau_{B,C} \otimes \id_D) \circ (\Delta_{A,B} \otimes \Delta_{C,D}))(\spe{a} \otimes \spe{b}).\end{align*} Thus we see that the product and coproduct are compatible.

Now we show that $\chi$ is a balanced character. Given $\spe{a} \in \spe{A}_N$ with $|N| \geq 2$, there exists a vertex $v$ such that $\{v\} \in \spe{a}$. Then $\Delta_{\{v\}, N \setminus \{v\}}(\spe{a}) \neq 0$.
By Theorem \ref{thm:chiconvex}, $\chi$ is a balanced convex character. \end{proof}

As an example, in Figure \ref{fig:antimatroid}, we have an antimatroid $\spe{a}$ and the corresponding coloring complex $\Sigma_{\chi}(\spe{a})$. The resulting quasisymmetric function is \[\Psi_{\chi}(\spe{a}, \mathbf{x}) = 2M_{1,2} + M_{2,1} + 4M_{1,1,1}.\]

We see that $\Psi_{\chi}(\spe{a}, \mathbf{x})$ enumerate all functions $f: I \to \mathbb{N}$ such that $f^{-1}([i]) \in \spe{a}$ and $f^{-1}(i)$ is a boolean lattice, for all $i$. This appears to be a new invariant. 

\begin{remark}
Consider the function $J_N: \spe{P}_N \to \spe{A}_N$ given by $J(\spe{p}) = \{I \subseteq N: I \mbox{ is an ideal of } \spe{p} \}$. Then $J: \spe{P} \to \spe{A}$ is a morphism of Hopf monoids, and $\chi \circ J = \chi$. Thus $\Psi_{\chi}(J(\spe{p}), \mathbf{x}) = \Psi_{\chi}(\spe{p}, \mathbf{x})$. In this sense, our new invariant can be regarded as a generalization of the quasisymmetric function which enumerates strict $\spe{p}$-partitions.
\end{remark}

\begin{corollary}
Let $\spe{a}$ be an antimatroid. Then the invariant $\Psi_{\varphi}(\spe{a}, \mathbf{x})$ is $M$-increasing. Moreover, $\Psi_{\varphi}(\spe{a},x)$ is strongly flawless.
\end{corollary}

We did not discuss $\zeta$. However, given an antimatroid $\spe{a}$, it is partially ordered by inclusion, and $\Sigma_{\zeta}(\spe{a}) = \Delta(\spe{a} \setminus \{\emptyset, N \})$, the order complex of the corresponding poset. In particular, $\Psi_{\zeta}(\spe{a}, \mathbf{x})$ is a $F$-quasisymmetric function of a graded poset, as defined in Example \ref{ex:ehrenborg}.
\begin{figure}
\begin{center}
\begin{tabular}{cc}

\begin{tikzpicture}

\node[circle, draw=black, fill=white, thick] (c) at (90:2cm) {$\{a, b, c\}$};
  \node[circle, draw=black, fill=white, thick] (b) at (150:2cm) {$\{a, b\}$};
  \node[circle, draw=black, fill=white, thick] (a) at (210:2cm) {$\{a\}$};
  \node[circle, draw=black, fill=white, thick] (f) at (270:2cm) {$\emptyset$};
  \node[circle, draw=black, fill=white, thick] (e) at (330:2cm) {$\{c\}$};
  \node[circle, draw=black, fill=white, thick] (d) at (30:2cm) {$\{b,c\}$};
  \node[circle, draw=black, fill=white,thick] (or) at (0:0cm) {$\{a, c\}$};
  \draw[thick] (a) -- (b) -- (c) -- (d) -- (e) -- (f) -- (a);
  \draw[thick] (c) -- (or);
  \draw[thick] (a) -- (or) -- (e);

\end{tikzpicture}

&

\begin{tikzpicture}

  \node[circle, draw=red, fill=white, dashed, thick] (b) at (150:2cm) {$\{a, b\}$};
  \node[circle, draw=black, fill=white, thick] (a) at (210:2cm) {$\{a\}$};

  \node[circle, draw=black, fill=white, thick] (e) at (330:2cm) {$\{c\}$};
  \node[circle, draw=red, fill=white, dashed, thick] (d) at (30:2cm) {$\{b,c\}$};
  \node[circle, draw=black, fill=white,thick] (or) at (0:0cm) {$\{a, c\}$};
 
  \draw[thick] (b) -- (a) -- (or) -- (e) -- (d);

\end{tikzpicture}

\end{tabular}
\end{center}
\caption{An antimatroid, and its corresponding coloring complex.}
\label{fig:antimatroid}
\end{figure}

\subsection{Generalized Permutohedra}
\label{subsec:permutohedra}

  The Hopf monoid of generalized permutohedra can be described in terms of generalized permutohedra, or in terms of submodular functions. The first description is in Section 5 of \cite{aguiar-ardila}, while the description in terms of submodular functions is in Section 12. We will focus on the latter description.
A \emph{submodular function} is a function $\spe{z}: 2^N \to \mathbb{R}$ which satisfies:
\begin{enumerate}
    \item $\spe{z}(\emptyset) = 0$,
    \item For $A, B \subseteq N$, we have $\spe{z}(A \cup B) + \spe{z}(A \cap B) \leq \spe{z}(A) + \spe{z}(B)$.
\end{enumerate}
We say $\spe{z}$ is modular if we have $\spe{z}(A \cup B) + \spe{z}(A \cap B) = \spe{z}(A) + \spe{z}(B)$.

Let $\spe{SF}_N$ be the set of all submodular functions with vertex set $N$. 
Given a bijection $\sigma: N \to M$, and a submodular function $\spe{z} \in \spe{SF}_N$, we define $\sigma(\spe{z})$ by $\sigma(\spe{z}) = \spe{z} \circ \sigma^{-1}$.
We see that we get a map $\spe{SF}_{\sigma}: \spe{SF}_N \to \spe{SF}_M.$ Thus $\spe{SF}$ is a set species. 

Now we turn $\mathbb{K}\spe{SF}^{\bdot}$ into a linearized Hopf monoid. Given disjoint sets $M$ and $N$, and submodular functions $\spe{w} \in \spe{SF}_M$ and $\spe{z}\in \spe{SF}_N$, we define $\spe{w} \cdot \spe{z}$ by $(\spe{w} \cdot \spe{z})(X) = \spe{w}(X \cap M) + \spe{z}(X \cap N)$.

Likewise, given a submodular function $\spe{z} \in \spe{A}_{M \sqcup N}$, we define $\spe{z}|_M$ by $\spe{z}|_M(S) = \spe{z}(S)$ for $S \subseteq M$. We define $\spe{z}/M$ by $\spe{z}/M(X) = \spe{z}(M \cup X) - \spe{z}(M)$, for $S \subseteq N$. Then $\Delta_{M,N}(\spe{z}) = \spe{z}|_M \otimes \spe{z}/M$.

Finally, given a submodular function $\spe{z} \in \spe{SF}_N$, we see that $\chi(\spe{z}) = 1$ if and only if $\spe{z}$ is modular. This is the same character by Aguiar and Ardila in Section 17 of \cite{aguiar-ardila}, which they refer to as the basic invariant (they define it for generalized permutohedra instead).

\begin{theorem}
The pair $(\mathbb{K}\spe{SF}, \mathbb{K}(\chi))$ is a linearized combinatorial Hopf monoid. Moreover, $\chi$ is a balanced convex character. Hence, for every submodular function $\spe{z}$, the quasisymmetric function $\Psi_{\chi}(\spe{z}, \mathbf{x})$ is $M$-increasing, and the polynomial $\chi_{\chi}(\spe{z}, \mathbf{x})$ is strongly flawless.
\end{theorem}
\begin{proof}
The fact that $\mathbb{K}\spe{SF}$ is a Hopf monoid is Theorem 12.2 in \cite{aguiar-ardila}. The fact that $\chi$ is balanced and convex follows from Theorem \ref{thm:chiconvex}. \end{proof}

The polynomial $\chi_{\chi}(\spe{z}, \mathbf{x})$ was introduced by Aguiar and Ardila \cite{aguiar-ardila}. In the context of generalized permutohedra, it counts functions that are maximized on a unique vertex, while $\Psi_{\chi}(\spe{z}, \mathbf{x})$ is a natural quasisymmetric function generalization of the polynomial.

Finally, Aguiar and Ardila have shown that many Hopf monoids arise as Hopf submonoids of $\mathbb{K}\spe{SF}$, or a related Hopf monoid of extended generalized permutohedra. We show that $\mathbb{K}\spe{A}$ \emph{does not} correspond to a linearized Hopf submonoid. This demonstrates that there are examples of Hopf monoids in the literature that can be studied via coloring complexes but \emph{not} generalized permutohedra. An \emph{extended generalized permutohedron} corresponds to an \emph{extended submodular function} $\spe{z}: N \to \mathbb{R} \cup \{\infty \}$ which has the property that, for any $A, B \subset N$ with $\spe{z}(A) < \infty$ and $\spe{z}(B) < \infty$, we have $\spe{z}(A \cup B) + \spe{z}(A \cap B) \leq \spe{z}(A) + \spe{z}(B)$. For a finite set $N$, let $\spe{ESF}_N$ be the set of extended submodular functiions. The corresponding linear species $\mathbb{K}\spe{ESF}$ is also a linearized combinatorial Hopf monoid, with the rule that $\Delta_{S,T}(\spe{z}) = \spe{z}|_S \otimes \spe{z}/S$ if $\spe{z}(S) < \infty$, and $\Delta_{S,T}(\spe{z}) = 0$ otherwise. 

For many examples of combinatorial Hopf monoids $(\mathbb{K}\spe{H}^{\bdot}, \varphi)$, Aguiar and Ardila  \cite{aguiar-ardila} define a linearized morphism $\psi: \mathbb{K}\spe{H}^{\bdot} \to \mathbb{K}\spe{ESF}^{\bdot}$ such that $\varphi = \chi \circ \psi$.
\begin{proposition}
There is no linearized homomorphism $\mathbb{K}(\psi): \mathbb{K}\spe{A}^{\bdot} \to \mathbb{K}\spe{ESF}^{\bdot}$ such that $\chi \circ \psi = \chi$.

There is also no injective linearized homomorphism $\mathbb{K}(\varphi):\mathbb{K}\spe{A}^{\bdot} \to \mathbb{K}\spe{ESF}^{\bdot}.$
\label{prop:counterexample}
\end{proposition}
\begin{proof}
Let $\mathbb{K}(\psi): \mathbb{K}\spe{A}^{\bdot} \to \mathbb{K}\spe{ESF}^{\bdot}$ such that $\chi \circ \psi = \chi$.
Given an extended submodular function $\spe{z}$, and sets $S, T$ with $\spe{z}(S) < \infty$ and $\spe{z}(T) < \infty$, it follows that $\spe{z}(S \cap T) < \infty$. In particular, if we let $\mathcal{C}(\spe{z})$ be the collection of sets $S$ with $\spe{z}(S) < \infty$, then $\mathcal{C}(\spe{z})$ is closed under intersections. We know that $\mathcal{C}(\spe{z}) \setminus \{\emptyset, N \}$ is also the vertex set of $\Sigma_{\chi}(\spe{z})$. On the other hand, given an antimatroid $\spe{a}$, the vertex set of $\Sigma_{\chi}(\spe{a})$ are the subsets $S \in \spe{a}$. Since $\chi \circ \psi = \chi,$ it follows that $\Sigma_{\chi}(\psi(\spe{a})) = \Sigma_{\chi}(\spe{a}).$ This implies that the subsets of $\spe{a}$ are closed under intersection \emph{for every} antimatroid. However, this is not true: the antimatroid in Figure \ref{fig:antimatroid} is a counterexample. Thus we have obtained a contradiction.

 Suppose that $\mathbb{K}(\varphi): \mathbb{K}\spe{A} \to \mathbb{K}\spe{ESF}$ is an injective homomorphism. We show that $\chi \circ \varphi = \chi$. Then our result follows from the fact that no such linearized homomorphism exists. Clearly, given an antimatroid $\spe{a}$, if $\chi(\spe{a}) = 1$, then $\spe{a} = \spe{a}_1 \cdots \spe{a}_k$ for several antimatroids $\spe{a}_i$, where each $\spe{a}_i$ has a singleton vertex set. Then $\varphi(\spe{a}) = \varphi(\spe{a}_1) \cdots \varphi(\spe{a}_k)$. We see each $\varphi(\spe{a}_i)$ is modular, and the product of two modular functions is modular. Thus $\chi(\varphi(\spe{a})) = 1$. 

Suppose then that there exists a finite set $N$ and an antimatroid $\spe{a} \in \spe{A}[N]$ with $\chi(\varphi(\spe{a})) = 1$. There is a linear order $\ell$ on $N$, and modular functions $\spe{z}_1, \ldots, \spe{z}_n$ such that $\varphi(\spe{a}) = \spe{z}_1 \cdots \spe{z}_n$, where $\spe{z}_i \in \spe{ESF}_{\{\ell_i\}}$. For each $i$, we see that \[\Delta_{\{\ell_1, \ldots, \ell_{i-1} \}| \ell_i|\ell_{i+1},\ldots, \ell_n}(\varphi(\spe{a})) = \spe{z}_1\cdots\spe{z}_{i-1} \otimes \spe{z}_i \otimes \spe{z}_{i+1}\cdots \spe{z}_n \neq 0.\] Since $\varphi$ is a homomorphism, we see that $\Delta_{\{\ell_1, \ldots, \ell_{i-1} \}| \ell_i|\ell_{i+1},\ldots, \ell_n} (\spe{a}) \neq 0$. Thus $\spe{a}|_{ \{\ell_1, \ldots, \ell_i \}} / \{\ell_1, \ldots, \ell_{i-1} \}$ exists for all $i$.
Let $\spe{a}_i = \spe{a}|_{\{\ell_1, \ldots, \ell_i \}} / \{\ell_1, \ldots, \ell_{i-1} \}$. Applying coproducts, and the fact that $\varphi$ is a homomorphism, we see that $\varphi(\spe{a}_i) = \spe{z}_i$. Then $\varphi(\spe{a}) = \varphi(\spe{a}_1 \cdots \spe{a}_n)$. Since $\varphi$ is injective, $\spe{a} = \spe{a}_1 \cdots \spe{a}_n$, and hence $\chi(\spe{a}) = 1$. Hence $\chi \circ \varphi = \chi.$ \end{proof}
In particular, there are natural examples of linearized combinatorial Hopf monoids which are not linearized submonoids of the Hopf monoid of extended generalized permutohedra, but which may be viewed as submonoids of the Hopf monoid of coloring complexes. A result similar to Proposition \ref{prop:counterexample} also holds for rooted connected graphs, with a similar proof. That is, there is no injective linearized homomorphism $\mathbb{K}(\varphi): \mathbb{K}\spe{RG}^{\bdot} \to \mathbb{K}\spe{ESF}^{\bdot}.$

\section{Acknowledgments}
The author wishes to thank Marcelo Aguiar, Jeremy Martin, and Frank Sottile for helpful conversations related to this project. We also thank the anonymous referees for many helpful comments.


\newcommand{\etalchar}[1]{$^{#1}$}
\providecommand{\bysame}{\leavevmode\hbox to3em{\hrulefill}\thinspace}
\providecommand{\MR}{\relax\ifhmode\unskip\space\fi MR }
\providecommand{\MRhref}[2]{%
  \href{http://www.ams.org/mathscinet-getitem?mr=#1}{#2}
}
\providecommand{\href}[2]{#2}

\end{document}